\documentclass[10pt,a4paper,reqno]{amsart}

\usepackage{amsmath,amssymb,amsthm}
\usepackage{mathtools}

\usepackage[T1]{fontenc}
\usepackage{lmodern}
\usepackage{float}
\usepackage[colorlinks=true,
            linkcolor=blue,
            citecolor=blue,
            urlcolor=blue]{hyperref}

\usepackage{todonotes}

\if {

}{{\caption}}
 \fi

\def\argmin{\mathop{\rm argmin}}

\def\1B{{\bf  1}}

\newcommand\be{\begin{equation}}
\newcommand\ee{\end{equation}}
\newcommand{\benl}{\begin{equation*}}
\newcommand{\eenl}{\end{equation*}}
\newcommand\ba{\begin{array}}
\newcommand\ea{\end{array}}
\newcommand{\bean}{\begin{eqnarray*}}
\newcommand{\eean}{\end{eqnarray*}}

\usepackage{cancel}
\usepackage{amsmath,amssymb,latexsym,amsthm,xcolor}
\usepackage{mathrsfs}  
\usepackage{mathtools}
\usepackage[margin=2cm]{geometry}
\usepackage{upgreek}
\usepackage{graphicx}

\providecommand{\Jreg}{\mathcal J}
\providecommand{\HS}{\mathcal{HS}}
\providecommand{\Proj}{\operatorname{Proj}}
\DeclareMathOperator*{\argminop}{arg\,min}
\providecommand{\argmin}{\argminop}
\providecommand{\intinf}{\operatorname{int}_{L^\infty}}
\providecommand{\Bor}{\operatorname{Bor}}

\newtheorem{theorem}{Theorem}[section]
\newtheorem{lemma}[theorem]{Lemma}
\newtheorem{proposition}[theorem]{Proposition}
\newtheorem{corollary}[theorem]{Corollary}

\usepackage[normalem]{ulem}
\usepackage{cancel}

\usepackage{soul}

\theoremstyle{remark}
{
    \newtheorem{definition}[theorem]{\bf Definition}
    
    \newtheorem{remark}[theorem]{\bf Remark}
    \newtheorem{assumption}{\bf Assumption}
}

\newcommand{\Vspace}{L^2([0,T]; \mathbb{R}^m)}
\newcommand{\Uad}{\mathcal{U}_{\rm ad}}


\begin{document}

\title[Asymptotics of SAA in optimal control]{
Asymptotic analysis of sample average approximations for optimal control problems with uncertain dynamics
}

\author[M.S. Aronna]{M. Soledad Aronna}
\address{M.S. Aronna\\ Escola de Matem\'atica Aplicada, Fundação Getulio Vargas, Praia de Botafogo 190,   22250-900 Rio de Janeiro - RJ, Brazil}
\email{soledad.aronna@fgv.br}

\author[J.M. Machado]{João Miguel Machado}
\address{J.M. Machado \\ CREST, ENSAE Paris, Institut Polytechnique de Paris\\
5 avenue Henry Le Chatelier\\
91120 Palaiseau, France}
\email{joaomiguel.machado@ensae.fr}

\author[O. Sierra]{Oscar Sierra}
\address{O. Sierra\\ Escola de Matem\'atica Aplicada, Fundação Getulio Vargas, Praia de Botafogo 190,   22250-900 Rio de Janeiro - RJ, Brazil}
\email{oscar.fonseca@fgv.br}

\begin{abstract}
We study sample average approximations (SAA) of optimal control problems with
control-affine dynamics depending on an uncertain parameter. The objective is an
expected terminal cost plus a Tikhonov regularization term, and the controls are
square-integrable functions subject to pointwise convex constraints. We show that
the reduced objective is twice continuously differentiable in the
square-integrable topology, so that no two-norm discrepancy arises, and that its
Hessian is a Hilbert--Schmidt perturbation of a multiple of the identity. This
yields Lipschitz regularity of stationary controls and second-order conditions
without gap for polyhedral constraints. For the SAA scheme with sample size $M$,
we prove almost sure uniform consistency of the empirical minimizers,
convergence rates of order $M^{-1/2}$ for the optimal values, independent of the
dimension of the parameter space
. We then derive the asymptotic distributions. The
optimal values satisfy a central limit theorem, obtained from a functional
central limit theorem for the SAA objective. The empirical minimizers are
asymptotically Gaussian when the constraints are inactive; under polyhedral
constraints, their limit solves a quadratic problem on the critical cone with a
Gaussian linear term. We also derive confidence intervals for linear functionals
of the optimal control and illustrate the results on two numerical examples.
\end{abstract}

\maketitle
\footnotetext[1]{This work was supported by FAPERJ (Rio de Janeiro, Brazil) process codes SEI-260003/004623/2021 and SEI-260003/000175/2024; and CNPq (Brazil) process code 312407/2023-8.}

\footnotetext[2]{
The second author was supported by the ERC project ``Optinfinite'' and would like to thank the hospitality of EMAp-FGV, which greatly helped this project.  
}

\section{Introduction}\label{Introduction}

Many control problems arising in applications involve dynamics that are only
partially known: physical parameters, initial states or model coefficients are
uncertain and are described by a probability distribution rather than by a
single value. A natural way to design a control that performs well for the
whole family of possible dynamics is to minimize the expected value of the
cost with respect to this distribution. This leads to \emph{average} (or
\emph{risk-neutral}, or \emph{ensemble}) optimal control problems, in which a
single open-loop control steers a parametrized family of systems and the
objective is the average, with respect to the probability distribution of the
uncertain parameter, of the terminal cost of each system. Optimal control
problems in which a single control must act on a whole family of systems appear
under several names: averaged and simultaneous control
\cite{zuazua2014,Zuzua-Loheac}, ensemble control of quantum systems
\cite{LiKhaneja2006,RuthsLi2012}, Riemann--Stieltjes optimal control in
aerospace guidance \cite{Ross2015,RossKarpenko2015}, control under uncertain
initial positions \cite{MarigondaQuincampoix2018}, and the training of deep
residual networks and neural ODEs, which can be formulated as a mean-field
optimal control problem in which the random initial condition encodes the
distribution of the data \cite{WEINAN2018}. Such problems also arise in
mathematical oncology, where ensemble optimal control has been used to design
drug schedules that manage resistance under uncertain tumour dynamics
\cite{ScagliottiEtAl2026cancer}, and control-affine dynamics appear in
robotics, for instance in the mid-air motion planning of free-floating robots
\cite{FanLiuBelabbas2019}, where the initial momentum imparted at take-off is
a natural source of uncertainty.

For risk-neutral problems, necessary optimality conditions in the form of a
Pontryagin principle for the averaged Hamiltonian were derived in
\cite{BettiolKhalil,BettiolKhalil2024} and, for control-affine systems with
uncertain initial data, in \cite{AronnaMonteiroSierra2025}, including a
characterization of bang--bang and singular arcs; a dynamic programming
approach was developed in \cite{AronnaPalladinoSierra2026}. Ensembles of
control-affine systems were studied in \cite{Scagliotti2023}, and their
worst-case and risk-averse counterparts in
\cite{Vinter2005,Scagliotti2025minimax} and \cite{ScagliottiSurowiec2026risk},
respectively; see also \cite{BonalliBonnet2023,KouriSurowiec2016} for
risk-averse stochastic and PDE-constrained problems.

A second line of research concerns the dependence of these problems on the
underlying distribution $\mu$, which is rarely available in closed form, or is
unknown or intractable, and has to be replaced by an approximation $\mu_k$, for
instance an empirical measure built from data. In
\cite{PesarePalladinoFalcone2021b}, the convergence of the value functions of
problems with unknown dynamics is established when the approximating
distributions converge in the $1$-Wasserstein distance. The linear-quadratic
case with uncertainty in the drift is analyzed in
\cite{PesarePalladinoFalcone.LQR.2021}, and a related analysis based on the $2$-Wasserstein distance appears in the mean-field setting of
\cite{WEINAN2018}. For control-affine dynamics with a quadratic control cost, a qualitative notion of convergence is investigated in \cite{palladino2025convergence}, through the notion of $\Gamma$-convergence, or equivalently epi-convergence, with respect to the weak topology of $L^2$. In particular, this yields the convergence of optimal controls in this topology.

The \emph{Sample Average Approximation} (SAA) method, first proposed in the
literature of stochastic programming \cite{shapiro2021}, consists of solving the
problem in which $\mu$ is replaced by the empirical measure $\mu_M$ of an
i.i.d.\ sample of size $M$. For optimal control problems with uncertain
parameters, the consistency of SAA, that is, the almost sure convergence of
optimal values and of cluster points of empirical minimizers, follows from the
epi-convergence of the sample average functionals \cite{ArtsteinWets1995,AttouchWets1983}
and was established in \cite{Phelps2013,PhelpsRoysetGong2016,Phelps2015}; see
also \cite{machado2025nash} for similar results in the context of game theory.
Epi-convergence, however, is purely qualitative. In the broader setting of
infinite-dimensional stochastic optimization, quantitative stability,
consistency and sample size estimates were obtained in
\cite{HoffhuesRomischSurowiec2021,MilzSurowiec2024,Milz2025stationary,Milz2023SAA,milz2024sample,LewBonalliPavone2024};
in particular, a different approach based on the minimal information distance
between probability measures is developed in \cite{romisch2024asymptotic} for
elliptic PDE-constrained problems. For optimal control of uncertain dynamical
systems, convergence rates have recently been obtained through metric entropy
arguments. In \cite{melnikov2026convergence}, for control-affine systems with a
strongly convex regularization, the gradients of the reduced cost lie in a
bounded subset of $W^{1,2}
$, and all critical controls, true
and empirical, are its image under the nonexpansive proximal mapping of the
regularizer. The covering numbers of $W^{1,2}$-balls in $L^2$ then yield
dimension-free rates of order $M^{-1/2}$ for the optimal values and the critical
points. Such entropy arguments in the space of controls thus avoid the curse of
dimensionality that affects stability estimates based on Wasserstein distances,
which inherit the rate $M^{-1/d}$ of $W_1(\mu_M,\mu)$ for parameters in
$\mathbb R^d$, $d\ge3$ \cite{fournier2015rate}.

In finite-dimensional stochastic programming, the next step is the asymptotic
distribution of the estimators
\cite{DupacovaWets1988,Shapiro1993,KingRockafellar1993,shapiro2021}. The
optimal values satisfy a central limit theorem, while the optimal solutions
converge, after rescaling by $\sqrt M$, to the solution of a random quadratic
problem posed on the critical cone, which is Gaussian only when the constraints
are inactive. In infinite dimension, central limit theorems for the optimal
values have been proved for linear elliptic PDE-constrained problems
\cite{romisch2024asymptotic} and for discrete-time stochastic optimal control
\cite{MilzShapiro2026}. For continuous-time control-affine problems with Mayer
costs, \cite{JaveedMilz2026} recently obtained a functional central limit
theorem for the SAA objective and confidence intervals for the optimal values,
without any regularization of the control. To the best of our knowledge, no
limit distribution is available for the SAA \emph{solutions} of optimal control
problems, where the decision variable lives in an infinite-dimensional space
and the relevant topologies are those of $L^2$ and $L^\infty$.

In this paper, we consider a metric measure space $(\Omega, d_\Omega, \mu)$ and
we study the convergence of the SAA method for systems of the form
$$
\min_{u\in\Uad}\ \mathbb{E}_{\mu}\bigl[ g\bigl(x^u(T,\omega),\omega\bigr)\bigr]
+\frac\alpha2\|u\|_{L^2}^2,
\qquad
\dot x^u=f_0(x^u,\omega)+\sum_{i=1}^m u_i\,f_i(x^u,\omega),
$$
where the admissible controls take values in a compact convex set
$U_c\subset\mathbb R^m$, the initial condition may be random, and the
Lipschitz constants of the terminal cost are allowed to be random and only
square integrable. Two structural features play a central role. The first is
the affine dependence of the dynamics on the control. It ensures that the cost
functional is twice Fr\'echet differentiable with respect to the $L^2$-topology,
so that the two-norm discrepancy typical of nonlinear optimal control problems
\cite{Ioffe1979III,casas2012} does not arise and the whole second-order analysis
can be carried out in the Hilbert space $L^2([0,T];\mathbb R^m)$; in addition,
the Hessian of the uncertain cost is represented by a compact operator. The
second is the Tikhonov term $\tfrac\alpha2\|u\|_{L^2}^2$. In the control-affine
setting it is not a mere technical convenience: since a compact operator is
never coercive in infinite dimension, it is the regularization that makes
second-order sufficient conditions possible, and it turns the first-order
optimality condition into a projection formula, from which the regularity of
optimal controls and the uniform convergence of the empirical minimizers
follow. Larger values of $\alpha$ improve the coercivity and convergence constants, at the price of a larger deviation from the unregularized problem.
\subsection*{Main contributions}
While first-order optimality conditions for ensemble problems have been
extensively studied \cite{BettiolKhalil,AronnaMonteiroSierra2025,ScagliottiSurowiec2026risk},
the second-order analysis below is, to the best of our knowledge, new in this
setting. Our main results are the following.

\noindent$\bullet$ \emph{Structure of the Hessian.} The Hessian of the averaged
cost is a Hilbert--Schmidt integral operator with bounded kernel and a square
integrable random element of the space of Hilbert--Schmidt operators
(Proposition  \ref{hessian.HS}), which yields laws of large numbers and central
limit theorems for empirical Hessians.

\noindent$\bullet$ \emph{Consistency, rates and stability.} The empirical
minimizers admit a measurable selection and converge almost surely to the true
minimizer in the uniform norm (Lemma  \ref{measurable.selection},
Proposition  \ref{uniform.consistency.Linfty}). Under the mere square
integrability of the Lipschitz constant of the cost, we obtain dimension-free
rates of order $M^{-1/2}$ for the optimal values and the squared $L^2$-error of
the minimizers (Proposition  \ref{prop.conv_rates_energy},
Remark  \ref{remark.rates.controls}).

\noindent$\bullet$ \emph{Asymptotic distributions.} We prove a central limit
theorem for the optimal values (Theorem  \ref{asymp.values}), extending these obtained in 
\cite{JaveedMilz2026} to regularized problems
with parameter-dependent costs. For polyhedral constraints,
the rescaled error 
of the empirical minimizers converges in distribution in $L^2$ to the solution of
a quadratic problem on the critical cone with a Gaussian linear term
(Theorem  \ref{asymp.active}), which is Gaussian when the constraints are
inactive or under strict complementarity (Corollary  \ref{asymp.gaussian}). These
appear to be the first limit theorems for SAA solutions in infinite-dimensional
optimal control. As a consequence, we identify the second-order term of the
optimal values, which quantifies their downward bias
(Corollary  \ref{cor.second.order.values}). A consistent covariance estimator
yields confidence intervals for linear functionals of the optimal control
(Corollary  \ref{inference}).

\noindent$\bullet$ \emph{Numerical illustration.} On an ensemble of oscillators
with uncertain frequency and on a sterile insect technique model, the observed
errors, fluctuations and confidence intervals agree closely with the asymptotic
predictions, including the vanishing of the fluctuations on active arcs
(Section  \ref{sec:numerics}).

\subsection*{Strategy of the proofs and organization of the paper}
Our analysis follows the finite-dimensional theory of statistical inference for
SAA estimators \cite[Chapter 5]{shapiro2021}, for
$\min_{x\in X}\mathbb E[F(x,\xi)]$ with $X\subset\mathbb R^n$ compact. Each of its ingredients needs a substitute for our framework, given by some sections of the paper:

\noindent (a) \emph{Smooth integrand} (Sections  \ref{sec:statemap}
and  \ref{sec:diff}). We show that the reduced cost, obtained by composing the
terminal cost with the control-to-state mapping, is twice Fr\'echet
differentiable in $L^2$, with Lipschitz constants that are only square
integrable with respect to the uncertain parameter. Its Hessian is the sum of
a Hilbert--Schmidt operator and a multiple $\alpha$ of the identity
(Proposition  \ref{hessian.HS}). This structure is used throughout: it makes the
second variation a Legendre form, and it allows us to apply laws of large
numbers and central limit theorems to the empirical Hessians.

\noindent (b) \emph{Compactness of the feasible set} (Sections  \ref{sec:necessary}
and  \ref{sec:quantitative}). The set of admissible controls is not compact in
$L^2$. The Tikhonov term turns the first-order optimality condition into a
pointwise projection formula, which shows that, outside an event of probability
of order $1/M$, all stationary controls, true and empirical, are Lipschitz in
time with a common constant (Corollary  \ref{regularity.stationary}). They thus
lie in a compact subset of $L^2$ whose metric entropy we estimate
(Lemma  \ref{entropy.K}), which yields uniform consistency and dimension-free
convergence rates.

\noindent (c) \emph{Optimal values} (Theorem  \ref{asymp.values}). On this compact
set, a functional central limit theorem holds for the sample average cost, and
the Danskin theorem and the Delta method apply as in
\cite[Theorem 5.7]{shapiro2021}.

\noindent (d)
\emph{Minimizers} (Section  \ref{sec:sosc} and
Subsection  \ref{subsec:minimizers}). In infinite dimension, bounded sequences
converge only weakly, and no functional central limit theorem is available for
the gradients in a norm of $C^1$ type. We therefore pass to the limit in the
first-order optimality condition rather than in the optimal value. Three facts
replace finite dimension: second-order sufficient conditions, which are
possible only because $\alpha>0$; the polyhedricity of the constraint set,
which removes the curvature term of \cite[Theorem 5.8]{shapiro2021}; and the
Legendre property of the second variation, which upgrades weak convergence to
strong convergence (Lemma  \ref{key.lemma} and Theorem  \ref{asymp.active}). The
second-order asymptotics of the optimal values then follows from that of the
minimizers (Corollary  \ref{cor.second.order.values}), in the reverse order of
\cite{shapiro2021}.

Section  \ref{sec:preliminary} introduces the problem and the SAA scheme,
Section  \ref{sec:numerics} contains the numerical examples,
Section  \ref{sec:conclusions} some concluding remarks, 
and the appendices
collect technical proofs and the abstract second-order theory.

\section{Notation, problem formulation and standing assumptions}\label{sec:preliminary}

\subsection{Functional-analytic notation}
 Norms and inner products without a subscript on $\mathbb R^k$ are Euclidean, and in subscripts we abbreviate
$\Vspace$ as $L^2$; for instance, $\mathcal L(L^2)$ denotes the space of bounded
linear operators on $\Vspace$, with norm $\|\cdot\|_{\mathcal L(L^2)}$, and
$\mathcal B(L^2)$ and $\mathcal B(L^2;\mathbb R^n)$ the spaces of bounded bilinear
maps on $\Vspace\times\Vspace$ with values in $\mathbb R$ and $\mathbb R^n$,
respectively. For a metric space $S$, $\Bor(S)$ denotes its Borel
$\sigma$-algebra. For a set $A$, $\mathbf 1_A$ denotes its indicator function. For a normed space $\mathbb X$, $v\in \mathbb X$ and $\varepsilon>0$, $B_{\mathbb X}(v,\varepsilon)$
denotes the open ball of radius $\varepsilon$ centered at $v.$ 

\subsection{Probabilistic notation}\label{notation.prob}
Let $\mathbb H$ be a separable Hilbert space and $(\Omega',\mathcal F',\mathcal P)$ a probability space with
expectation $\mathbb E$; in our applications it is either $(\Omega,\Bor(\Omega),\mu)$ of
Subsection  \ref{subsec:problem}, with expectation denoted by $\mathbb E_\mu$, or
$(\Omega_{\mathbb P},\mathcal F,\mathbb P)$ of Subsection  \ref{SAA}. A measurable map
$\upchi:(\Omega',\mathcal F')\to(\mathbb H,\Bor(\mathbb H))$ is called a \emph{random element} of $\mathbb H$. If
$\mathbb E\bigl[\|\upchi\|_{\mathbb H}^{2}\bigr]<\infty$, its \emph{mean} and \emph{covariance operator} are the Bochner
integrals
$
\mathbb E[\upchi]:=\int_{\Omega'}\upchi\,d\mathcal P\in\mathbb H,$ $
\operatorname{Cov}(\upchi):=\int_{\Omega'}(\upchi-\mathbb E\upchi)\otimes(\upchi-\mathbb E\upchi)\,d\mathcal P
\in\mathcal L(\mathbb H),
$
where $(u\otimes v)h:=\langle v,h\rangle_{\mathbb H}\,u$ for $u,v,h\in\mathbb H$. A sequence $(\nu_k)$ of Borel probability measures on a metric space $S$
\emph{converges narrowly} (or \emph{weakly}) to $\nu$, written $\nu_k\xrightarrow{\rm n}\nu$, if
$\int_S\phi\,d\nu_k\to\int_S\phi\,d\nu$ for every bounded continuous
$\phi:S\to\mathbb R$. Random elements $\upchi_M,\upchi$ with values in $S$,
possibly defined on different probability spaces, \emph{converge in
distribution}, written $\upchi_M\xrightarrow[]{\mathfrak D}\upchi$, if their laws
converge narrowly.
A $\mathbb{H}$-valued random element $\upchi$ is said to be \emph{Gaussian} if,
for every $h\in\mathbb{H}$, the real-valued random variable
$\langle \upchi,h\rangle_{\mathbb{H}}$ has a (possibly degenerate) Gaussian distribution. In that
case, if $\mathfrak m:=\mathbb{E}[\upchi]$ and $\mathscr K :=\operatorname{Cov}(\upchi)$,
we write $\upchi\sim\mathcal N(\mathfrak m,\mathscr K)$; then, for $h\in\mathbb{H}$,
$\langle \upchi,h\rangle_{\mathbb{H}}\sim
\mathcal N(\langle \mathfrak m,h\rangle_{\mathbb{H}},\langle \mathscr Kh,h\rangle_{\mathbb{H}})$.

A sequence of random elements $\upchi_M$ with values in a metric space $(S,d)$
converges to $\upchi$ \emph{almost surely} (a.s.) if
$\mathcal P\bigl(d(\upchi_M,\upchi)\to0\bigr)=1$, and \emph{in probability} if
$\mathcal P\bigl(d(\upchi_M,\upchi)>\varepsilon\bigr)\to0$ for every
$\varepsilon>0$. For real random variables $Z_M$ and positive numbers $r_M$, we
write $Z_M=o_p(r_M)$ if $Z_M/r_M\to0$ in probability, and we say that $(Z_M)$
is \emph{bounded in probability} if
$\lim_{R\to\infty}\sup_M\mathcal P(|Z_M|>R)=0$; for random elements of a normed
space, these notions refer to their norms. 

\subsection{The problem and assumptions}\label{subsec:problem}
Let  $(\Omega,d_\Omega)$ be a {\it Polish space} (complete, separable metric space) endowed with its Borel
$\sigma$-algebra $\Bor(\Omega)$, and denote by $\mathscr P(\Omega)$ the set of Borel
probability measures on $\Omega$. For $\mu\in\mathscr P(\Omega)$ and $\alpha>0$
we consider the problem
\begin{equation}\label{Pmu}\tag{$P_\mu$}
\min_{u\in\Uad}\ \Jreg_\mu(u):=J_\mu(u)+\frac\alpha2\|u\|_{L^2}^2,
\qquad
J_\mu(u):=\int_\Omega g\bigl(x^u(T,\omega),\omega\bigr)\,d\mu(\omega),
\end{equation}
where, for each $\omega\in\Omega$, $x^u(\cdot,\omega)$ solves
\begin{equation}\label{dynamics}
\dot x(t)=f\bigl(x(t),u(t),\omega\bigr)\ \text{ for a.e. }t\in[0,T],
\qquad x(0)=\varphi(\omega),
\end{equation}
with a measurable initial condition $\varphi:\Omega\to K_0$, $K_0\subset\mathbb R^n$
compact, and the dynamics is affine in the control:
\begin{equation}\label{f.affine}
f(x,u,\omega)=f_0(x,\omega)+\sum_{i=1}^m u_i\,f_i(x,\omega)=f_0(x,\omega)+F(x,\omega)\,u,
\ \ F:=\big[f_1\,|\,\cdots\,|\,f_m\big].
\end{equation}
For $t\in[0,T]$ we set
$g_t(u,\omega):=g(x^u(t,\omega),\omega)$. 
We consider the following set of admissible controls 
$$
    \Uad:=\{u\in\Vspace:\ u(t)\in U_c\text{ for a.e. }t\in[0,T]\},    
$$
where $U_c\subset\mathbb R^m$ is compact and convex, and we denote by
$$
    V_\mu:=\inf_{\Uad}\Jreg_\mu
$$
the value of \eqref{Pmu} with data distribution $\mu$. When $\mu$ is fixed we write $J$, $\Jreg$ and $V$. We set $\rho:=\max_{v\in U_c}\|v\|$ and since $\Uad$ has empty interior in $\Vspace$, differentiability statements are
formulated on the open ball
\begin{equation}\label{O.def}
\mathcal O:=\{u\in\Vspace:\ \|u\|_{L^2}<r_\rho\},\qquad r_\rho:=\sqrt T\rho+1,
\end{equation}
which contains $\Uad$ together with a margin, so that derivatives can be taken
at points of $\Uad$ in arbitrary directions; the value of $r_\rho$ plays no role.

\begin{assumption}\label{A1}
For every $i=0,\dots,m$ and $\omega\in\Omega$, $f_i(\cdot,\omega)\in
C^2(\mathbb R^n)$, and $f_i(x,\cdot)$, $\partial_xf_i(x,\cdot)$,
$\partial^2_{xx}f_i(x,\cdot)$ are measurable for every $x$. There is $L\ge1$ such
that $$\|\partial_x^kf_i(x,\omega)-\partial_x^kf_i(y,\omega)\|\le L\|x-y\|$$ for
$k=0,1,2$ and $\|f_i(0,\omega)\|\le L$, for all $x,y$, $\omega$ and $i$.
\end{assumption}
\begin{remark}\label{A1.implications}
Assumption  \ref{A1} implies, for every $\omega\in\Omega$, $i=0,\dots,m$ and
$x\in\mathbb R^n$,
$
\|\partial_xf_i(x,\omega)\|,\ \|\partial^2_{xx}f_i(x,\omega)\|\le L,$ $
\|f_i(x,\omega)\|\le L(1+\|x\|),
$
since a $C^1$ map that is $L$-Lipschitz has derivative bounded by $L$, and
$\|f_i(x,\omega)\|\le\|f_i(0,\omega)\|+L\|x\|$.

Using
$\|u\|_1\le\sqrt m\,\|u\|$, the estimates above yield, for all
$(x,u)\in\mathbb R^n\times\mathbb R^m$ and every $\omega\in\Omega$,
\begin{equation}\label{growth.f}
\begin{aligned}
&\|f(x,u,\omega)\|\le L(1+\|x\|)\big(1+\sqrt m\,\|u\|\big),  &&\|f_u(x,\omega)\|\le L\sqrt m\,(1+\|x\|),
\\
&\|f_x(x,u,\omega)\|,\ \|f_{xx}(x,u,\omega)\|
\le L\big(1+\sqrt m\,\|u\|\big), &&\|f_{xu}(x,\omega)\|\le L\sqrt m.
\end{aligned}
\end{equation}
Moreover, for all $x,y\in\mathbb{R}^n,$ $u,v\in \mathbb{R}^m,$ and all $\omega\in \Omega,$
one obtains the Lipschitz estimates
\begin{equation}\label{Lipschitz.bounds}
\begin{aligned}
&\|f(x,u,\omega)-f(y,v,\omega)\|
\le
L\big(1+\sqrt m\,\|u\|\big)\,\|x-y\|
+
L\sqrt m\,(1+\|y\|)\,\|u-v\|, \\
&\|f_x(x,u,\omega)-f_x(y,v,\omega)\|, \,\|f_{xx}(x,u,\omega)-f_{xx}(y,v,\omega)\| \\
& \hspace{5.4cm}
\le
L\big(1+\sqrt m\,\|u\|\big)\,\|x-y\|
+
L\sqrt m\,\|u-v\|, \\
&\|f_u(x,\omega)-f_u(y,\omega)\|,\,\|f_{xu}(x,\omega)-f_{xu}(y,\omega)\|
\le
L\sqrt m\,\|x-y\|.
\end{aligned}
\end{equation}
\end{remark}

\begin{assumption}\label{A2}
For every $\omega\in\Omega$, $g(\cdot,\omega)\in C^2(\mathbb R^n)$, and
$g(x,\cdot)$, $\nabla_xg(x,\cdot)$, $\nabla_{xx}g(x,\cdot)$ are measurable for
every $x$. Moreover $g(0,\cdot)\in L^2(\mu,\Omega)$ and there is
$L_g:\Omega\to[1,\infty)$ with $L_g\in L^2(\mu,\Omega)$ such that $g(\cdot,\omega)$,
$\nabla_xg(\cdot,\omega)$ and $\nabla_{xx}g(\cdot,\omega)$ are
$L_g(\omega)$-Lipschitz for every $\omega$.
\end{assumption}

\begin{remark}\label{bounded.deriv.g}
Under Assumption  \ref{A2}, for all $x\in\mathbb R^n$ and all $\omega \in \Omega$,
$$
\|\nabla_x g(x,\omega)\|,
\,
\|\nabla_{xx}g(x,\omega)\|
\le
L_g(\omega),
\qquad
|g(x,\omega)|\le|g(0,\omega)|+L_g(\omega)\|x\|.
$$
\end{remark}

\subsection{Sample average approximation}\label{SAA}
Let $(\Omega_{\mathbb P},\mathcal F,\mathbb P)$ be a probability space supporting
a sequence $(X_i)_{i\in\mathbb N}$ of independent random elements
$X_i:\Omega_{\mathbb P}\to\Omega$ with common law $\mu$. For $M\in\mathbb N$,
define the empirical random measure $\mu_M:=\frac1M\sum_{i=1}^M\delta_{X_i}$.
The associated SAA functionals and optimal value are
\begin{equation}\label{JM}
J_M(u):=\frac1M\sum_{i=1}^Mg_T(u,X_i),\qquad
\Jreg_M(u):=J_M(u)+\frac\alpha2\|u\|_{L^2}^2,\qquad
V_M:=\inf_{u\in\Uad}\Jreg_M(u).
\end{equation}
We write $\mathbb E$ for the expectation with respect to $\mathbb P$, so that $\mathbb P$ and $\mathbb E$ refer to the
randomness of the sample, whereas $\mu$ is the law of a single observation. Expectations, variances and covariances with
respect to $\mu$ are denoted by $\mathbb E_\mu$, $\operatorname{Var}_\mu$ and $\operatorname{Cov}_\mu$; thus
$\mathbb E_\mu[Z]=\int_\Omega Z\,d\mu$ for every $\mu$-integrable $Z$ on $\Omega$ (a Bochner integral if $Z$ is
vector-valued), and in particular $\mathbb E[Z(X_1)]=\mathbb E_\mu[Z]$.

\begin{remark}\label{JM.notation}
For every sample-space outcome, the empirical measure $\mu_M$ belongs to
$\mathscr P(\Omega)$, and the SAA functionals are the objective functionals associated with it:
$$
J_M=J_{\mu_M},\qquad \Jreg_M=\Jreg_{\mu_M}.
$$
Assumptions  \ref{A1}-- \ref{A2} hold with $\mu$ replaced by $\mu_M$, since the
integrability conditions in Assumption  \ref{A2} are trivially satisfied by a
measure with finite support. Hence every deterministic statement proved for a
generic $\mu\in\mathscr P(\Omega)$ applies, outcome by outcome, to the SAA
problem, with constants that depend on the sample only through
$\|L_g\|_{L^1(\mu_M,\Omega)}$ and $\|L_g\|_{L^2(\mu_M,\Omega)}$; by the strong
law of large numbers, these converge almost surely to the corresponding
quantities for $\mu$. We are mainly interested in the asymptotic behavior, as $M\to\infty$, of
\emph{empirical minimizers}, that is, random elements $\hat u_M$ such that
$\hat u_M\in\argminop_{\Uad}\Jreg_M$ for every sample-space outcome; their
existence follows from Theorem  \ref{existence}, and their measurability from
Lemma  \ref{measurable.selection}.
\end{remark}

\section{The control-to-state mapping and the cost functional}\label{sec:statemap} In this section we collect the well-posedness
and Lipschitz continuity of both the state equation \eqref{dynamics} and its
associated adjoint equation, as well as the differentiability of the control-to-state
mapping. The proofs of these results are provided in the Appendix 
 \ref{appendixA}.

\subsection{Well-posedness of the state-equation and Lipschitz continuity of the control-to-state operator}

\begin{proposition}[Well-posedness of the state equation]
\label{wellposedness}
Suppose that Assumption  \ref{A1} holds, and for $u \in \Vspace, $ set
$
  \Lambda_u := L\big(T + \sqrt{mT}\,\|u\|_{L^2}\big).
$
Then, for every $u \in \Vspace$ and $\omega \in \Omega$:
\begin{enumerate}
  \item[(i)] \emph{Existence and uniqueness.} System \eqref{dynamics} admits a unique solution
        $x^u(\cdot,\omega) \in W^{1,1}([0,T];\mathbb{R}^n)$.
  \item[(ii)] A {\it priori} bound. For all $t \in [0,T]$,
        $\|x^u(t,\omega)\| \le \big(\|\varphi(\omega)\| + \Lambda_u\big)\, e^{\Lambda_u}.
        $
  \item[(iii)] \emph{Uniform bound on $\mathcal{O}$.} There exists a compact set
        $\mathcal{X}_0 \subset \mathbb{R}^n$, independent of $u \in \mathcal{O}$, such that,
        for every $u \in \mathcal{O}$,
        $x^u(t,\omega) \in \mathcal{X}_0 $ for all $ (t,\omega) \in [0,T]\times\Omega .
        $
        Consequently, $x^u(\cdot,\omega) \in W^{1,2}([0,T];\mathbb{R}^n)$.
\end{enumerate}
\end{proposition}

\begin{proof}
    See Appedndix  \ref{App.wellposedness}
\end{proof}

\begin{remark}\label{bounded.deriv.f}
Let $\mathcal{X}_0:=\overline{B(0,R)}$ be the compact set given in Proposition  \ref{wellposedness}. If $x,y\in \mathcal{X}_0$ and $u,v\in U_c$, the estimates in \eqref{Lipschitz.bounds} become
$$
\|D^\gamma f(x,u,\omega)-D^\gamma f(y,v,\omega)\|\le L_f\big(\|x-y\|+\|u-v\|\big), \ \ |\gamma|\le2,
$$
where
$
L_f:=L\max\big\{1+\sqrt m\,\rho,\ \sqrt m\,(1+R)\big\}.
$
Moreover, by \eqref{growth.f}, $\|F(x,\omega)\|=\|f_u(x,\omega)\|\le L_f$ and
$\|f_{xu}(x,\omega)\|\le L_f$ for all $x\in \mathcal{X}_0$, and $\|f_x(x,u,\omega)\|\le L_f$ for
all $x\in \mathcal{X}_0$ and $u\in U_c$. We stress that the last bound, and hence every
estimate below in which $\|f_x\|$ or $\|f_{xx}\|$ is bounded by $L_f$, requires
$u(t)\in U_c$; for controls $u\in\mathcal O$ we only have
$\|f_x(x^u(t,\omega),u(t),\omega)\|\le L(1+\sqrt m\|u(t)\|)$, whose integral over
$[0,T]$ is bounded by $\Lambda_{r_\rho}$.
\end{remark}

\begin{remark}[Local version of the assumptions]\label{local.assumptions}
Assumptions  \ref{A1}-- \ref{A2} are needed only along trajectories. Hence all
the results remain valid, with $\mathcal{X}_0$ replaced by $\mathcal{X}_0'$, if there are a compact
set $\mathcal{X}_0'$ and a bounded open set $\mathcal N\supset \mathcal{X}_0'$ such that the
following holds. For every $u\in\mathcal O$ and $\omega\in\Omega$,
\eqref{dynamics} has a solution on $[0,T]$ with values in $\mathcal{X}_0'$, and the
conditions of Assumptions  \ref{A1}-- \ref{A2} hold with $x,y$ restricted to
$\mathcal N$. Indeed, let $\psi\in C_c^\infty(\mathcal N)$ with $\psi=1$ near
$\mathcal{X}_0'$. Then $\psi f_i$ and $\psi g$ satisfy Assumptions  \ref{A1}-- \ref{A2}
globally. By uniqueness, they generate the same trajectories, costs and
derivatives on $\mathcal O$ as the original data.
\end{remark}

\begin{proposition}[Lipschitz continuity of the control-to-state mapping] \label{traj.lipschitz} Suppose that Assumption  \ref{A1} holds. Then, for all $\omega\in\Omega,$ the control-to-state mapping from $\mathcal O$ to $C([0,T];\mathbb{R}^n)$ given by $ u \mapsto x^u(\cdot,\omega) $ is Lipschitz continuous, with a Lipschitz constant $L_x$ depending only on $T,L,m,K_0$ and $r_\rho$. \end{proposition}
\begin{proof}
    See Appendix  \ref{App.traj.lipschitz}
\end{proof}

\subsection{The adjoint equation}
Let $u\in\mathcal O$ and $\Phi^u(t,s,\omega)$ denote the fundamental matrix associated with the
linearized system along the trajectory $x^u(\cdot,\omega)$ i.e., $\Phi^u(t,s,\omega)$ is the fundamental matrix of
$\dot y=f_x(x^u(t,\omega),u(t),\omega)y.$
By Proposition  \ref{wellposedness}, $x^u(\cdot,\omega)\in \mathcal{X}_0$, and hence
\eqref{growth.f} yields, for a.e. $t\in[0,T]$,
$
\bigl\|f_x(x^u(t,\omega),u(t),\omega)\bigr\|
\leq
L\bigl(1+\sqrt{m}\,\|u(t)\|\bigr).
$
Since $u\in\mathcal{O}$, the right-hand side belongs to $L^1([0,T])$.
Consequently, the classical existence and uniqueness theory for linear
ordinary differential equations with integrable coefficients
(see, e.g., \cite[Section 18.3]{Kurzweil1986}) yields the existence of a unique absolutely
continuous solution $\Phi^u(\cdot,s,\omega)$ 
. Moreover,
Gr\"onwall's inequality gives
\begin{equation}\label{f.EDO.bound}
\|\Phi^u(t,s,\omega)\|
\leq
e^{
\int_s^t
L(1+\sqrt{m}\,\|u(\tau)\|)\,d\tau}
\leq
e^{\,\Lambda_u}
\leq
e^{\Lambda_{r_\rho}},
\end{equation}
uniformly for $u\in\mathcal{O}$, $0\le s\le t\le T$ and $\omega\in \Omega$. For $t\in[0,T],$ $u\in\mathcal O$ and all $\omega \in\Omega$ consider the adjoint equation
\begin{equation}\label{adjoint.t}
\begin{cases}
\dot p(s,\omega)
=
-f_x\bigl(x^u(s,\omega),u(s),\omega\bigr)^\top p(s,\omega),
& s\in[0,t],\\
p(t,\omega)
=
\nabla_x g\bigl(x^u(t,\omega),\omega\bigr).
\end{cases}
\end{equation}
Using the standard identity
$
\frac{d}{d s}\Phi^u(t,s,\omega)
=
-\Phi^u(t,s,\omega)
f_x\bigl(x^u(s,\omega),u(s),\omega\bigr),
$
one checks that
$$
p(s,\omega)
=
\Phi^u(t,s,\omega)^\top
\nabla_x g\bigl(x^u(t,\omega),\omega\bigr), \ \
s\in[0,t],
$$
solves \eqref{adjoint.t}, and by uniqueness for linear equations with integrable
coefficients it is the unique solution; in particular
$
p(\cdot,\omega)\in W^{1,1}([0,t];\mathbb R^n).
$

\begin{lemma}[Properties of the adjoint state]
\label{p.lc.w.t.u}
Suppose that Assumptions  \ref{A1}-- \ref{A2} hold, let $t\in[0,T]$ and, for
$u\in \mathcal{O}$, let $p^u(\cdot,\omega)$ denote the solution of the adjoint
equation \eqref{adjoint.t} associated with the control $u$. Then:
\begin{itemize}
    \item [(i)] For every $u\in\mathcal O$ and $\omega\in\Omega$,
$
\sup_{s\in[0,t]}\|p^u(s,\omega)\|\le e^{\Lambda_{r_\rho}}L_g(\omega);
$
in particular $p^u(s,\cdot)\in L^2(\mu,\Omega;\mathbb{R}^n)$ for all $s\in[0,t]$.
Moreover, for all $u,v\in\mathcal O$ and $\omega \in \Omega,$
$\sup_{s\in[0,t]}\|p^u(s,\omega)-p^v(s,\omega)\|\leq L_{p}(\omega)\|u-v\|_{L^2},$
for some $L_p\in L^{2}(\mu,\Omega;\mathbb{R}^{+})$ depending on $\Lambda_{r_\rho}$ and $L_g.$
\item[(ii)] If $u\in\Uad$, then $p^u(\cdot,\omega)\in W^{1,\infty}([0,t];\mathbb{R}^n)$ and
$\|\dot p^u(\cdot,\omega)\|_{L^\infty}\le L_fe^{\Lambda_{r_\rho}}L_g(\omega)$ for all $\omega \in\Omega.$
\end{itemize}
\end{lemma}
\begin{proof}
  See Appendix  \ref{App.p.lc.w.t.u}
\end{proof}

\subsection{Fr\'echet differentiability of the control-to-state mapping} As a consequence of the existence and uniform boundedness
properties of the fundamental matrix, the first variational equation is
well-posed and admits uniform estimates with respect to the perturbations
in $\Vspace$.

In the following results, we consider the first- and second-order variations as functions of both $t$ and $\omega$, that is,
$
\delta x=\delta x(t,\omega),
$ $
\delta^2x=\delta^2x(t,\omega).
$
For convenience, we drop the dependence on $\omega$ from the notation whenever there is no risk of ambiguity.

\begin{proposition}[Properties of the first variation]
\label{Diff.x}
Let Assumption  \ref{A1} hold and, for every $\omega\in\Omega$, 
$u\in\mathcal O$, and $\delta u\in\Vspace$, consider the variational
equation
\begin{equation}\label{variational}
\begin{cases}
\dot{\delta x}(s)
=
f_x\bigl(x^u(s,\omega),u(s),\omega\bigr)\delta x(s)
+
f_u\bigl(x^u(s,\omega),\omega\bigr)\delta u(s),
& s\in[0,T],\\
\delta x(0)=0.
\end{cases}
\end{equation}
Then, the following statements hold.
\begin{itemize}
\item[(i)] The variational equation \eqref{variational} admits a unique solution in $C([0,T];\mathbb R^n)$.
Moreover, for all $t\in[0,T],$  the mapping
$
\delta u\mapsto \delta x(t)
$
is linear and bounded from $\Vspace$ into $\mathbb R^n.$  In particular, there exists a constant
$C_1>0$, depending only on $T,L,m,K_0$ and $r_\rho$, such that
$
\|\delta x\|_{C([0,T];\mathbb R^n)}
\le
C_1\|\delta u\|_{L^2}.
$
\item[(ii)] For all $t\in[0,T],$  the map from $\mathcal O$ to $\mathcal L(L^2;\mathbb R^n)$ given by $u\mapsto \delta x^u(t)$ is Lipschitz continuous, where $\delta x^u$ is the solution of \eqref{variational} with control $u.$
\end{itemize}
\end{proposition}
\begin{proof}
    See Appendix  \ref{App.variational}
\end{proof}

\begin{theorem}[Second-order differentiability of the control-to-state mapping]\label{Diff.x2}
Let Assumption  \ref{A1} hold.
Then, for every  $\omega\in\Omega$, the following statements hold.
\begin{itemize}
    \item[(i)] The control-to-state mapping
is twice Fr\'echet
differentiable in $\mathcal O.$ Moreover, for every $t\in[0,T]$ and $\delta u_1,\delta u_2\in\Vspace,$ the second Fr\'echet differential is given by
$
D^2x^u(t,\omega)[\delta u_1,\delta u_2]=
\delta^2 x_{12} (t),
$
where $\delta^2x_{12}$ is the solution of the second
variational equation
\begin{equation}\label{variational2}\begin{cases}
\frac{d}{dt}{\delta^2x}_{12}
=
f_x(s)\,\delta^2x_{12}
+
f_{xx}(s)[\delta x_1,\delta x_2]
+
f_{xu}(s)[\delta x_1,\delta u_2]
+
f_{ux}(s)[\delta u_1,\delta x_2],
&s\in[0,t],
\\
\delta^2x_{12}(0)=0,
\end{cases}
\end{equation}
where $f_x(s):=f_x(x^u(s),u(s),\omega),$ $f_{ij}(s):=f_{ij}(x^u(s),u(s),\omega)$ for $i,j\in \{x,u\},$  $x^u(s):=x^u(s,\omega) $ and $\delta x_i$ denotes the solution of the variational equation
\eqref{variational} associated with the perturbation $\delta u_i$,
$i=1,2.$
 \item[(ii)]
The map $
u\mapsto
D^2x^u(t,\omega)
$
is Lipschitz continuous from
$\mathcal O$ to
${\mathcal B}(L^2;\mathbb R^n)$.
\end{itemize}
\end{theorem}
\begin{proof}
    See Appendix  \ref{App.Diff.x2}
\end{proof}

\begin{remark}\label{no.L2.diff}
The control-affine structure of the dynamics is what allows us to work with the
$L^2$-topology in Theorem  \ref{Diff.x2}. For a general dynamics of class $C^2$
with non-vanishing bounded and Lipschitz continuous $f_{uu}$, the second variational equation
contains the additional term $f_{uu}(s)[\delta u,\delta u]$, and the expansion of
$f$ with respect to the control contains the remainder
$$
\int_0^1(1-\theta)\,\Big(f_{uu}\bigl(x^u(s),u(s)+\theta\delta u(s),\omega\bigr)
-f_{uu}\bigl(x^u(s),u(s),\omega\bigr)\Big)
[\delta u(s),\delta u(s)]\,d\theta ,
$$
which can only be estimated by $C\|\delta u(s)\|^3$. Its contribution to
$\|r(\cdot,\omega)\|_{L^\infty}$ is therefore controlled by
$\int_0^T\|\delta u(s)\|^3\,ds=\|\delta u\|_{L^3}^3$, a quantity which is not
dominated by $\|\delta u\|_{L^2}^3$ and may even be infinite for
$\delta u\in\Vspace$. 
\end{remark}

\begin{corollary}\label{Diff.x2Linfty}
Suppose Assumption  \ref{A1} holds. Then, for every $\omega\in\Omega$, the control-to-state mapping
$
u\mapsto x^u(\cdot,\omega)
$
is twice continuously Fr\'echet differentiable on $\mathcal O\cap L^\infty([0,T];\mathbb{R}^m)$ with respect to the $L^\infty$-norm. Its first and second Fr\'echet differentials are given by the solutions of
\eqref{variational} and \eqref{variational2}, respectively.
\end{corollary}
\begin{proof}
    See Appendix  \ref{App.Diff.x2Linfty}
\end{proof}

\section{Differentiability of the cost function}\label{sec:diff}

To obtain the Fr\'echet derivatives of the objective functional $J_\mu$ (and
hence, those of $\Jreg_\mu$), we begin by studying
the Fr\'echet differentiability properties of the operator $g_t$. The derivation follows the framework developed in \cite{polak2012}. As explained in
Section  \ref{sec:preliminary}, all the statements below are formulated on the
open set $\mathcal O,$ introduced in \eqref{O.def}.

\begin{lemma}\label{gTlipschitz}
Let $\mu\in\mathscr P(\Omega)$ and suppose that Assumptions  \ref{A1}-- \ref{A2} hold. Then, for
every $t\in[0,T]$:
\begin{itemize}
\item[(i)] For every $\omega\in\Omega$ and all $u,v\in\mathcal O$,
$
|g_t(u,\omega)-g_t(v,\omega)|\le L_T(\omega)\,\|u-v\|_{L^2},
$ with $
L_T(\omega):=L_g(\omega)\,L_x,
$
where $L_x\le\sqrt T\,L_f\,e^{\Lambda_{r_\rho}}$ is the constant of Proposition  \ref{traj.lipschitz}. Thus
$L_T\in L^2(\mu,\Omega)$. Moreover, the map $(u,\omega)\mapsto g_t(u,\omega)$ is
$\Bor(\Vspace)\otimes\Bor(\Omega)$-measurable.
\item[(ii)] For every $\omega\in\Omega$ and every $u\in\mathcal O$,
$
|g_t(u,\omega)|\le a(\omega):=|g(0,\omega)|+L_g(\omega)\,R,
$ and $a\in L^2(\mu,\Omega),
$
where $R$ is given by Proposition  \ref{wellposedness}.
\item[(iii)] Consequently, $g_t(u,\cdot)\in L^2(\mu,\Omega)$ for every
$u\in\mathcal O$, so that $J_\mu$ is real valued on $\mathcal O$, and
$
|J_\mu(u)-J_\mu(v)|\le\|L_T\|_{L^1(\mu,\Omega)}\,\|u-v\|_{L^2}
$ for all $  u,v\in\mathcal O .
$
In particular, $\Jreg_\mu$ is real valued and Lipschitz continuous on
$\mathcal O$, and it is bounded from below on $\Uad$ by
$-\|a\|_{L^1(\mu,\Omega)}$.
\end{itemize}
\end{lemma}

\begin{proof}
(i) By Proposition  \ref{traj.lipschitz} and the Lipschitz property of
$g(\cdot,\omega)$ in Assumption  \ref{A2},
$$
\begin{aligned}
|g_t(u,\omega)-g_t(v,\omega)|
&=\big|g(x^u(t,\omega),\omega)-g(x^v(t,\omega),\omega)\big|\\
&\le L_g(\omega)\,\|x^u(t,\omega)-x^v(t,\omega)\|
 \le L_g(\omega)L_x\|u-v\|_{L^2}.
 \end{aligned}
$$
Since $L_x$ is deterministic and $L_g\in L^2(\mu,\Omega)$, we get
$L_T\in L^2(\mu,\Omega)$; the bound $L_x\le\sqrt T L_fe^{\Lambda_{r_\rho}}$
follows from $L\sqrt m(1+R)\le L_f$. For the measurability, fix $u\in\Vspace$.
The Picard iterates in the proof of Lemma  \ref{lem:picard}, started at the
constant function $\varphi(\omega)$, are measurable in $\omega$ for each $t$, since
$f$ is continuous in $x$ and measurable in $\omega$ and $\varphi$ is measurable;
their pointwise limit $\omega\mapsto x^u(t,\omega)$ is therefore measurable, and so
is $\omega\mapsto g(x^u(t,\omega),\omega)$, $g$ being continuous in $x$ and
measurable in $\omega$. Hence $g_t$ is a Carath\'eodory function, continuous in
$u$ by the first part of (i) and measurable in $\omega$, and it is jointly
measurable since $\Vspace$ is separable and metrizable
\cite[Lemma 4.51]{aliprantis2006infinite}. (ii) By Remark  \ref{bounded.deriv.g} and
$\|x^u(t,\omega)\|\le R$ (Proposition  \ref{wellposedness}),
$$
|g_t(u,\omega)|\le|g(0,\omega)|+L_g(\omega)\|x^u(t,\omega)\|
\le a(\omega).
$$
Since $g(0,\cdot)\in L^2(\mu,\Omega)$ and
$L_g\in L^2(\mu,\Omega)$, we have $a\in L^2(\mu,\Omega)$. (iii) The first assertion is immediate from (i) and (ii); the estimate follows
by integrating (i) with respect to $\mu$. The statements about $\Jreg_\mu$
follow since on the bounded set $\mathcal O$ the map
$u\mapsto\tfrac\alpha2\|u\|_{L^2}^2$ is Lipschitz continuous and nonnegative.
\end{proof}

\subsection{First order Fr\'echet derivatives}
From now on, we fix a probability measure $\mu\in \mathscr{P}(\Omega)$ and, for notational simplicity, we write $J$ in place of $J_\mu$ and $\Jreg$ in place of $\Jreg_\mu$.
\begin{proposition}\label{1order}
Suppose that Assumptions  \ref{A1}-- \ref{A2} hold and let $t\in[0,T]$. Then:
\begin{itemize}
\item[(i)] For all $\omega\in\Omega$ the map $g_t(\cdot,\omega)$ is Fr\'echet
differentiable on $\mathcal O$, with
$$
Dg_t(u,\omega)[\delta u]=\langle\nabla g_t(u,\omega),\delta u\rangle_{L^2},
\quad
\nabla g_t(u,\omega)(s)=\mathbf 1_{[0,t]}(s)\,
F\big(x^u(s,\omega),\omega\big)^{\!\top}p(s,\omega),
$$
where $p(\cdot,\omega)$ solves the adjoint equation \eqref{adjoint.t}. Moreover,
$u\mapsto\nabla g_t(u,\omega)$ is Lipschitz on $\mathcal O$, with a Lipschitz
constant $\ell(\omega)$ such that $\ell\in L^2(\mu,\Omega)$.
\item[(ii)] For every $u\in\mathcal O$ one has
$\nabla g_t(u,\cdot)\in L^2(\mu,\Omega;\Vspace)$, with
$$
\|\nabla g_t(u,\omega)\|_{L^2}
\le
L_\nabla(\omega):=\sqrt T\,L_f\,e^{\Lambda_{r_\rho}}L_g(\omega).
$$
\item[(iii)] $J$ is Fr\'echet differentiable on $\mathcal O$, with
Lipschitz gradient given by
$
\nabla J=\int_\Omega\nabla g_T(\cdot,\omega)\,d\mu(\omega),
$
the integral being understood in the Bochner sense in $\Vspace$. In particular,
$\Jreg$ is Fr\'echet differentiable on $\mathcal O$ with Lipschitz gradient
\begin{equation}\label{grad.Jreg}
\nabla\Jreg(u)(s)
=
\int_\Omega
F\big(x^u(s,\omega),\omega\big)^{\!\top}p(s,\omega)\,d\mu(\omega)
+\alpha\,u(s),
\  \text{ for a.e. } s\in[0,T],
\end{equation}
where $p$ solves \eqref{adjoint.t} with $t=T$.
\end{itemize}
\end{proposition}

\begin{proof}
(i) follows from Proposition  \ref{Diff.x} and an argument analogous to that of \cite[Corollary 5.6.9]{polak2012},
the Lipschitz constant being controlled by $L_g$ through Lemma  \ref{p.lc.w.t.u}(i).
For (ii), Lemma  \ref{p.lc.w.t.u}(i) gives
$\sup_{s\in[0,t]}\|p(s,\omega)\|\le e^{\Lambda_{r_\rho}}L_g(\omega)$, which,
combined with $\|F(x,\omega)\|\le L_f$ for $x\in \mathcal{X}_0$ (Remark  \ref{bounded.deriv.f}),
yields the stated bound, and $L_\nabla\in L^2(\mu,\Omega)$ since $L_g\in L^2(\mu,\Omega)$.
Finally, (iii) follows from (i)--(ii): by Lemma  \ref{gTlipschitz}(i) the
difference quotients of $g_T(\cdot,\omega)$ are dominated by the integrable
function $L_T$, uniformly in $u\in\mathcal O$, so one may differentiate under
the integral; the Lipschitz continuity of $\nabla J$ follows by integrating the
Lipschitz estimate in (i).
\end{proof}

\begin{lemma}[Time regularity of the gradients]\label{gradient.regularity}
Suppose that Assumptions  \ref{A1}-- \ref{A2} hold. Then, for every $u\in\Uad$ and every
$\omega\in\Omega$, one has $\nabla g_T(u,\omega)\in W^{1,\infty}([0,T];\mathbb R^m)$,
with
$$
\|\nabla g_T(u,\omega)\|_{L^\infty}\le L_f\,e^{\Lambda_{r_\rho}}L_g(\omega),
\qquad
\Bigl\|\tfrac{d}{dt}\nabla g_T(u,\omega)\Bigr\|_{L^\infty}\le C_\Psi\,L_g(\omega),
$$
where
$
C_\Psi
:=
L_f e^{\Lambda_{r_\rho}}(C_{\dot x}+L_f)
$
and $C_{\dot x}:=L(1+R)(1+\sqrt m\,\rho).$
Consequently, $\nabla J(u)\in W^{1,\infty}([0,T];\mathbb R^m)$, with
the same bounds with $L_g(\omega)$ replaced by $\|L_g\|_{L^1(\mu,\Omega)}$.
\end{lemma}

\begin{proof}
Fix $u\in\Uad$ and $\omega\in\Omega$, and write $x(t):=x^u(t,\omega)$ and
$p(t):=p(t,\omega)$, so that $\nabla g_T(u,\omega)(t)=F(x(t),\omega)^{\!\top}p(t)$
by Proposition  \ref{1order}. Since $\|u(t)\|\le\rho$ a.e. and $\|x(t)\|\le R$,
\eqref{growth.f} gives $\|\dot x(t)\|\le C_{\dot x}$ for a.e. $t$. By
Lemma  \ref{p.lc.w.t.u}, $\|p(t)\|\le e^{\Lambda_{r_\rho}}L_g(\omega)$ and
$\|\dot p(t)\|\le L_f e^{\Lambda_{r_\rho}}L_g(\omega)$. Since $\|F(x,\omega)\|\le L_f$ and
$F(\cdot,\omega)$ is $L_f$-Lipschitz on $\mathcal{X}_0$ (Remark  \ref{bounded.deriv.f}), the first bound follows, and
adding and subtracting $F(x(s),\omega)^{\!\top}p(t)$ yields, for all
$s,t\in[0,T]$, 
$
\|\nabla g_T(u,\omega)(t)-\nabla g_T(u,\omega)(s)\|
\le\bigl(L_fC_{\dot x}+L_f^2\bigr)e^{\Lambda_{r_\rho}}L_g(\omega)\,|t-s|
=C_\Psi L_g(\omega)\,|t-s| .
$
The statement on $\nabla J(u)$ follows by integrating with respect to $\mu$.
\end{proof}

\subsection{Second order
Fr\'echet derivatives}

As in the previous results, we first deduce the second-order derivatives  the integrand $g_t.$ We stress that here the control-affine structure
plays an essential role: it is precisely because
$f_{uu}\equiv 0$ and $f_u$, $f_{xu}$ do not depend on $u$  that second-order differentiability holds
directly on the Hilbert space $\Vspace$, and not only with respect to the
$L^\infty$-topology.

We define the {\it random Pontryagin Hamiltonian} $H:\mathbb{R}^n\times\mathbb{R}^m\times\mathbb{R}^n\times \Omega\to \mathbb{R}$
as
$$
H(x,u,p,\omega)
:= \langle p,f(x,u,\omega) \rangle
= \langle p,f_0(x,\omega)\rangle+\langle p,F(x,\omega)u\rangle .
$$
The following result is an adaptation of \cite[Theorem 5.6.10]{polak2012} to our setting. For the sake of completeness, we include its proof in Appendix  \ref{appendixA}.
\begin{proposition}\label{2order}
Suppose that Assumptions  \ref{A1}-- \ref{A2} hold. Then,  for all $\omega  \in \Omega$ and $t \in [0,T]$,
\begin{enumerate}
\item[(i)] the function $g_t(\cdot,\omega)$ is twice
 Fr\'echet differentiable on $\mathcal O.$ Moreover, for every $u\in\mathcal O$ and every $\delta u,\delta v \in \Vspace, $ the second order Fr\'echet differential
is given by
\begin{equation}\label{F2order}
\begin{aligned}
D^2 g_t(u,\omega)&[\delta u,\delta v]
=
\langle
\nabla_{xx} g(x(t),\omega)
\mathcal S(\delta u)(t),\mathcal S(\delta v)(t)\rangle
\\
&\hspace{2cm}  +
\int_0^t
\bigg\langle
\begin{pmatrix}
H_{xx}(s) & H_{xu}(s)\\
H_{ux}(s) & 0
\end{pmatrix}
\begin{pmatrix}
\mathcal S(\delta u)(s)\\
\delta u(s)
\end{pmatrix}, \begin{pmatrix}
\mathcal S(\delta v)(s)\\
\delta v(s)
\end{pmatrix}\bigg\rangle_{\mathbb{R}^{n}\times\mathbb{R}^{m}} ds,
\end{aligned}
\end{equation}
where $H_{ij}(s):=\partial^2_{ij} H(x(s),u(s),p(s),\omega)$ for ${i,j}\in \{x,u\},$ $x(s):=x^u(s,\omega)$ solves \eqref{dynamics} with control $u,$ $p(s):=p(s,\omega)$ solves the adjoint equation \eqref{adjoint.t} with control $u$,
$
\mathcal S : \Vspace
\to
C([0,t];\mathbb R^n)
$ is the linear and bounded operator given by $\mathcal S(\delta u)(s)=\delta x(s)$, where $\delta x$  solves the variational equation \eqref{variational}.
\item[(ii)] The Hessian mapping $u\mapsto D^2 g_t(u,\omega)$ is Lipschitz continuous on $\mathcal O$: there exists $C_{\mathcal H}\in L^2(\mu,\Omega;\mathbb{R}^+)$ depending on  $T,$ $L_f $ and $L_g$ such that
\begin{equation}\label{D2gLIPS}
\|D^2 g_t(u,\omega)
-
D^2 g_t(v,\omega)\|_{\mathcal B(L^2)}\leq C_{\mathcal H}(\omega)\|u-v\|_{L^2} 
\ \text{ 
 for all } \  u,v\in \mathcal O.
 \end{equation}

\item[(iii)] For every $u\in\mathcal O$, the mapping
$
\omega\mapsto D^2g_t(u,\omega)
$
belongs to
$
L^2\left(\mu,\Omega;\mathcal B(L^2)\right).
$
\end{enumerate}
\end{proposition}

\begin{proof}
See Appendix  \ref{proof.2order}.
\end{proof}

\begin{lemma}\label{2orderw}
Suppose that Assumptions  \ref{A1}-- \ref{A2} hold. Then
\begin{enumerate}
\item[(i)] $J$  is twice continuously Fr\'echet differentiable on $\mathcal O$ and for all $u\in\mathcal O$ and all $\delta u,\delta v \in \Vspace,$ the second-order Fr\'echet differential of $J$ is given by
    $
    D^2 J(u)[\delta u, \delta v] = \int_{\Omega} D^2 g_T(u,\omega)[\delta u, \delta v] \, d\mu(\omega).
    $
 \item[(ii)] The Hessian mapping $u\mapsto D^2 J(u)$   is Lipschitz continuous on $\mathcal O$.
 \item[(iii)] The same statements hold for $\Jreg$, with
 $
 D^2\Jreg(u)[\delta u,\delta v]
 =
 D^2J(u)[\delta u,\delta v]+\alpha\langle\delta u,\delta v\rangle_{L^2},
 $
 and $u\mapsto D^2\Jreg(u)$ has the same Lipschitz constant as
 $u\mapsto D^2J(u)$.
    \end{enumerate}
\end{lemma}
\begin{proof}
    See Appendix  \ref{Proof.lemma.2orderw}.
\end{proof}

\subsection{The Hessian as a compact perturbation of the identity}\label{section.hessian_perturbation}

The previous results show that, for every $u\in\mathcal O$, the bilinear form
$D^2 J(u)$ is continuous on $\Vspace$. By the Riesz representation theorem
we may therefore associate with it a bounded operator; the next proposition
shows that this operator is, in fact, of Hilbert--Schmidt type. This structural
information is the key ingredient both for the second-order analysis of
Section  \ref{sec:sosc} and for the asymptotic results of
Section  \ref{sec:asymptotics}.

We denote by $\HS(\Vspace)$ the space of Hilbert--Schmidt operators on
$\Vspace$, that is, the space of bounded operators $\mathcal T$ on $\Vspace$
such that $\sum_{k=1}^{\infty}\|\mathcal Te_k\|_{L^2}^2<\infty$ for
any orthonormal basis $\{e_k\}$ of $\Vspace$. Endowed with the inner product
$\langle \mathcal T_1,\mathcal T_2\rangle_{\HS}=\sum_{k=1}^{\infty}\langle \mathcal
T_1e_k,\mathcal T_2e_k\rangle_{L^2},$ the space $\HS(\Vspace)$ is a separable
Hilbert space \cite[Theorem 4.4.5]{Hsing2013}, one has
$ \|\cdot\|_{\mathcal L(L^2)}\le\|\cdot\|_{\HS}$, and every
Hilbert--Schmidt operator is compact \cite[Theorem 4.4.3]{Hsing2013}. Recall
finally that every kernel $\mathbf K\in L^2([0,T]^2;\mathbb R^{m\times m})$
defines, through
$$
\mathcal T[\delta u](\tau):=\int_0^T\mathbf K(\tau,\sigma)\,\delta u(\sigma)\,d\sigma ,
$$
an operator $\mathcal T\in\HS(\Vspace)$ with
$\|\mathcal T\|_{\HS}^2=\int_0^T\int_0^T\|\mathbf K(\tau,\sigma)\|^2_{\rm
F}\,d\tau\,d\sigma$, where $\|\cdot\|_{\rm F}$ denotes the Frobenius norm.
\begin{proposition}[Kernel representation of the Hessian]\label{hessian.HS}
Suppose that Assumptions  \ref{A1}-- \ref{A2} hold and let $u\in\mathcal O$ and
$\omega\in\Omega$. Then there exists a kernel
$
\mathbf K^{u,\omega}\in L^\infty\big([0,T]^2;\mathbb R^{m\times m}\big)
$
such that the operator
$\mathcal H_g(u,\omega)\in\mathcal L(L^2)$ defined by
\begin{equation}\label{eq:H.definition}
\mathcal H_g(u,\omega)[\delta u](\tau)
:=
\int_0^T\mathbf K^{u,\omega}(\tau,\sigma)\,\delta u(\sigma)\,d\sigma ,
\quad  \tau\in[0,T],
\end{equation}
is the Riesz representative of $D^2g_T(u,\omega)$, that is,
$
D^2g_T(u,\omega)[\delta u,\delta v]
=
\langle \mathcal H_g(u,\omega)\delta u,\delta v\rangle_{L^2}
$ for $\delta u,\delta v\in\Vspace .
$
Moreover, there is a constant $c_{\HS}>0$, depending only on
$T,m,n,L_f$ and $\Lambda_{r_\rho}$, such that
\begin{equation}\label{HS.bound}
\|\mathcal H_g(u,\omega)\|_{\HS}
\le
c_{\HS}\,L_g(\omega)
\  \ \text{ for all } u\in\mathcal O,\ \omega\in\Omega .
\end{equation}
In particular, $\mathcal H_g(u,\cdot)\in L^2\big(\mu,\Omega;\HS(\Vspace)\big)$,
and the operator
$
\mathcal H_J(u):=\int_\Omega \mathcal H_g(u,\omega)\,d\mu(\omega)
$
is a self-adjoint Hilbert--Schmidt operator on $\Vspace$
representing $D^2J(u)$. Consequently, the Riesz representative of
$D^2\Jreg(u)$ is
\begin{equation}\label{H.def}
\mathcal H(u)=\mathcal H_J(u)+\alpha\,\mathrm{Id},
\end{equation}
that is, a compact perturbation of $\alpha\,\mathrm{Id}$.
\end{proposition}

\begin{proof}
See Appendix  \ref{proof.hessian.HS}.
\end{proof}

\begin{remark}\label{extension.hessian}
For dynamics that are not affine in the control, $\Jreg$ is in general twice
Fr\'echet differentiable only in the $L^\infty$-topology, while second-order
conditions are formulated in $L^2$: this is the classical two-norm discrepancy
\cite{Ioffe1979III,casas2012} (see also Remark  \ref{no.L2.diff}). Here
$f_{uu}\equiv0$, hence $H_{uu}\equiv0$, and Lemma  \ref{2orderw} gives
$\Jreg\in C^2(\mathcal O)$ in the $L^2$-topology itself. Thus no discrepancy
arises, and the theory of \cite{casas2012} applies with
$U_\infty=U_2=\Vspace$; the derivative formulas also hold on
$\mathcal O\cap L^\infty([0,T];\mathbb R^m)$, by the continuous embedding
$L^\infty\hookrightarrow L^2$.
\end{remark}

\section{Existence of solutions and first-order optimality conditions}\label{sec:necessary}

We first show that \eqref{Pmu}, and hence each SAA problem, admits a solution,
so that the empirical minimizers $\hat u_M$ are well defined. The argument
relies on the control-affine structure, which makes the control-to-state map
continuous from the weak topology of $\Vspace$ to the strong topology of
$C([0,T];\mathbb R^n)$.

\begin{lemma}[Weak continuity of the control-to-state mapping]\label{weak.continuity}
Suppose that Assumption  \ref{A1} holds and let $\omega\in\Omega$. If
$\{u_k\}\subset\Uad$ and $u_k\rightharpoonup u$ weakly in $\Vspace$, then
$u\in\Uad$ and
$
x^{u_k}(\cdot,\omega)\rightarrow x^{u}(\cdot,\omega)
$ strongly in $ C([0,T];\mathbb R^n).
$
\end{lemma}
\begin{proof}
    See Appendix  \ref{App.weak.continuity}
\end{proof}

\begin{theorem}[Existence of optimal controls]\label{existence}
Suppose that Assumptions  \ref{A1}-- \ref{A2} hold. Then problem \eqref{Pmu} admits at least one
global solution $\hat u\in\Uad$. The same holds for the SAA problems: for every
$M\in\mathbb N$, the random functional $\Jreg_M$ attains its infimum on
$\Uad$ pointwise on $\Omega_{\mathbb P}$.
\end{theorem}

\begin{proof}
The set $\Uad$ is convex, bounded and closed in $\Vspace$, hence weakly compact
by reflexivity. Let $\{u_k\}\subset\Uad$ be a minimizing sequence for $\Jreg$;
up to a subsequence, $u_k\rightharpoonup\hat u\in\Uad$. By
Lemma  \ref{weak.continuity}, $x^{u_k}(T,\omega)\to x^{\hat u}(T,\omega)$ for
every $\omega\in\Omega$, so that $g_T(u_k,\omega)\to g_T(\hat u,\omega)$ by
continuity of $g(\cdot,\omega)$. Since $|g_T(u_k,\omega)|\le a(\omega)$ with
$a\in L^1(\mu,\Omega)$ by Lemma  \ref{gTlipschitz}(ii), dominated convergence
yields $J(u_k)\to J(\hat u)$; that is, $J$ is weakly continuous on $\Uad$. As
the $L^2$-norm is weakly lower semicontinuous,
$
\Jreg(\hat u)\le\liminf_{k\to\infty}\Jreg(u_k)=\inf_{u\in\Uad}\Jreg(u),
$
so $\hat u$ is a global minimizer. The argument applies to
$\Jreg_M$, which is the cost associated with the empirical random measure
$\mu_M$, by Remark  \ref{JM.notation}.
\end{proof}
Since $\Uad$ is convex and $\Jreg$ is Fr\'echet differentiable on the open set
$\mathcal O\supset\Uad$, the first-order necessary conditions take the form of
a variational inequality. Because of the regularization term, this inequality can be
solved pointwise in terms of a projection, which in turn yields regularity of
the stationary controls. This last property will be crucial in
Sections  \ref{sec:quantitative} and  \ref{sec:asymptotics}. Throughout this section we will use the projection operator onto the set of pointwise control constraints
$
\Proj_{U_c}:\mathbb R^m \ni u 
    \mapsto 
    \argmin_{v \in U_c} \frac{1}{2}\|u - v\|^2
    \in U_c $
to formulate optimality conditions. Since $U_c$ is nonempty, closed and convex, $\Proj_{U_c}$ is a maximally monotone
operator \cite[Example 20.32]{BauschkeCombettes2017} and nonexpansive
\cite[Proposition 4.16]{BauschkeCombettes2017}, i.e., 
$
    \|\Proj_{U_c}(u_1)-\Proj_{U_c}(u_2)\|\le\|u_1-u_2\|$
for all $
    u_1,u_2\in\mathbb R^m. 
$

\begin{definition}
We say that $\hat u\in\Uad$ is a \emph{local minimizer of \eqref{Pmu} in the
$L^\infty$-topology} if there exists $\varepsilon>0$ such that
$
\Jreg(\hat u)\le \Jreg(u)$  for every $
u\in\Uad \cap
B_{L^\infty}(\hat u,\varepsilon).
$
If, in addition, the inequality is strict for every
$u\in\Uad\cap B_{L^\infty}(\hat u,\varepsilon)$ with $u\neq\hat u$, then
$\hat u$ is called a \emph{strict local minimizer in the $L^\infty$-topology};
the corresponding notions in the $L^2$-topology are defined analogously.
\end{definition}

Observe that every local minimizer in the $L^2$-topology is also a local
minimizer in the $L^\infty$-topology, since
$B_{L^\infty}(\hat u,\varepsilon/\sqrt T)\subset B_{L^2}(\hat u,\varepsilon)$
for every $\varepsilon>0$.

\begin{definition}
We say that $\hat u\in\Uad$ is a {\it stationary point} of \eqref{Pmu} if
\begin{equation}\label{VI}
\langle\nabla\Jreg(\hat u),u-\hat u\rangle_{L^2}\ge0
\ \ \text{ for all } u\in\Uad .
\end{equation}
For $\hat u\in\Uad$, we write
$$
N_{\Uad}(\hat u):=\bigl\{\xi\in\Vspace:\ \langle\xi,u-\hat u\rangle_{L^2}\le0
\ \text{ for all }u\in\Uad\bigr\}
$$
for the {\em normal cone} of $\Uad$ at $\hat u$. 
\end{definition}
From latter definition, we can deduce that $\hat u$ is a stationary point if,
and only if, $-\nabla\Jreg(\hat u)\in N_{\Uad}(\hat u)$.

\begin{remark}
Under Assumptions  \ref{A1}-- \ref{A2}, if  $\hat u\in\Uad$ is a local minimizer
of \eqref{Pmu} in the $L^2$- or in the $L^\infty$-topology, then it is a
stationary point. Indeed, let $u\in\Uad$. Since $\Uad$ is convex,
$\hat u+\lambda(u-\hat u)\in\Uad$ for every $\lambda\in(0,1]$ and, for
$\lambda$ sufficiently small, $\hat u+\lambda(u-\hat u)$ belongs to the
corresponding neighbourhood of $\hat u$ both in the $L^2$- and in the
$L^\infty$-topology, $u-\hat u$ being essentially bounded. Therefore,
$
0
\le
\lim_{\lambda\downarrow0}
\left[\Jreg\big(\hat u+\lambda(u-\hat u)\big)-\Jreg(\hat u)\right]/{\lambda}
=
\langle\nabla\Jreg(\hat u),u-\hat u\rangle_{L^2},
$
by Proposition  \ref{1order}. Hence $\hat u$ satisfies \eqref{VI}.
\end{remark}

\begin{theorem}[Characterization of stationary points]\label{necessary.conditions}
Suppose that Assumptions  \ref{A1}-- \ref{A2} hold and let $\hat u\in\Uad.$ Then the following statements are
equivalent:
\begin{itemize}
\item[(i)] $\hat u$ satisfies \eqref{VI};
\item[(ii)] for a.e. $t\in[0,T]$,
\begin{equation}\label{projection.formula}
\hat u(t)
=
\Proj_{U_c}
\left(
-\frac1\alpha\,\Psi_\mu[\hat u](t)
\right),
\ \
\text{where} \,\, \Psi_\mu[u](t)
:=
\int_\Omega
F\big(x^{u}(t,\omega),\omega\big)^{\!\top}p(t,\omega)\,d\mu(\omega),
\end{equation}
\item[(iii)] for a.e. $t\in[0,T]$,
$
\hat u(t)
=
\displaystyle \argmin_{v\in U_c}
\Bigg\{
 \int_\Omega H\big(x^{\hat u}(t,\omega),v,p(t,\omega),\omega\big)\,d\mu(\omega)
+\frac\alpha2\|v\|^2
\Bigg\},
$

\noindent the minimizer being unique.
\end{itemize}
Here $p(\cdot,\omega)$ denotes the solution of the adjoint equation
\eqref{adjoint.t} associated with $\hat u$ and $t=T$.
\end{theorem}
\begin{proof}
 (i)$\Leftrightarrow$(ii).
By \eqref{grad.Jreg},
$
\nabla\Jreg(\hat u)
=
\Psi_\mu[\hat u]+\alpha\hat u,
$
so that \eqref{VI} is equivalent to
$$
\int_0^T
\big\langle
\Psi_\mu[\hat u](t)+\alpha\hat u(t),
u(t)-\hat u(t)
\big\rangle\,dt
\ge0
 \ 
\text{ for all }u\in\Uad.
$$
We first prove that (i) implies (ii). Fix $v\in U_c$ and define
$
h_v(t):=
\bigl\langle
\alpha\hat u(t)+\Psi_\mu[\hat u](t),
v-\hat u(t)
\bigr\rangle.
$
Since $\Psi_\mu[\hat u]$ is continuous on $[0,T]$ and
$\|\hat u\|_{L^\infty}\le\rho$, we have $h_v\in L^1([0,T])$. Let $A\subset[0,T]$ be an arbitrary measurable set and define
$
u(t)
:=
v\,\mathbf 1_A(t)
+
\hat u(t)\,\mathbf 1_{[0,T]\setminus A}(t).
$
Then $u\in\Uad$ and, since $u-\hat u$ vanishes outside $A$,
\eqref{VI} gives
$
\int_A h_v(t)\,dt\ge0
$  for every measurable $A\subset[0,T].
$
Therefore, by  \cite[Theorem 3.2.5(a)]{florescu2021lebesgue},
$
h_v(t)\ge0
$ for a.e. $t\in[0,T].$ The null set obtained in this way depends on $v$. Since $U_c$ is compact,
it is separable, and we may choose a countable dense subset
$
D=\{v_k\}_{k\in\mathbb N}\subset U_c.
$
For each $k$, let $N_{v_k}$ be a null set such that $h_{v_k}(t)\ge0$ for every $t\notin N_{v_k}$, and set
$
N:=\bigcup_{k\in\mathbb N}N_{v_k}.
$
Then $N$ is a null set and, for every $t\notin N$, $h_{v_k}(t)\ge0$ for all $k$. Moreover,
for fixed $t\notin N$, the map
$
v\mapsto
h_v(t)
$
is continuous. Hence, by the density of $D$ in $U_c$, $h_v(t)\ge0$ for all $v\in U_c$ and all $t\notin N$. Equivalently,
$$
\left\langle
-\frac1\alpha\Psi_\mu[\hat u](t)-\hat u(t),
v-\hat u(t)
\right\rangle
\le0
\
\text{ for all }v\in U_c,
\ \text{ for a.e. } t\in[0,T].
$$
Thus, by the variational characterization of the orthogonal projection onto
a closed convex set \cite[Theorem 5.2]{brezis2010functional},
(ii) holds. Conversely, suppose that (ii) holds. By the same variational characterization,
$
\bigl\langle
\Psi_\mu[\hat u](t)+\alpha\hat u(t),
v-\hat u(t)
\bigr\rangle
\ge0
$ for all $v\in U_c$ and a.e. $t\in[0,T]$. Given $u\in\Uad$, we may take
$v=u(t)$ in the previous inequality and integrate over $[0,T]$, which gives \eqref{VI}.

\medskip

\noindent
(ii)$\Leftrightarrow$(iii).
For fixed $t\in[0,T]$, define
$
c_0(t)
:=
\int_\Omega
\left\langle
p(t,\omega),
f_0\bigl(x^{\hat u}(t,\omega),\omega\bigr)
\right\rangle
\,d\mu(\omega).
$
Then, for every $v\in U_c$, completing the square,
$$
\begin{aligned}
\int_\Omega
H\bigl(x^{\hat u}(t,\omega),v,p(t,\omega),\omega\bigr)
\,d\mu(\omega)
+\frac{\alpha}{2}\|v\|^2
&=
c_0(t)
+
\left\langle
\Psi_\mu[\hat u](t),v
\right\rangle
+
\frac{\alpha}{2}\|v\|^2
=
\frac{\alpha}{2}
\left\|
v+\frac1\alpha\Psi_\mu[\hat u](t)
\right\|^2
+c(t),
\end{aligned}
$$
where
$
c(t)
:=
c_0(t)
-
\frac{1}{2\alpha}
\left\|\Psi_\mu[\hat u](t)\right\|^2
$
is independent of $v$. Since $U_c$ is closed and convex, by \cite[Theorem 5.2]{brezis2010functional}, the
minimizer of the right-hand side over $U_c$ is unique and equal to
$
\Proj_{U_c}
\bigl(
-\frac1\alpha\Psi_\mu[\hat u](t)
\bigr).
$
Hence (ii) and (iii) are equivalent.
\end{proof}
\begin{remark}[The case $\alpha=0$: Pontryagin's principle]\label{remark.PMP}
Statement (iii) is the corresponding Hamiltonian form. For $\alpha=0$, since
the Hamiltonian is affine in $v$ and $U_c$ is convex, \eqref{VI} is equivalent
to Pontryagin's principle for the averaged Hamiltonian,
$$
\hat u(t)\in\argmin_{v\in U_c}
\int_\Omega H\bigl(x^{\hat u}(t,\omega),v,p(t,\omega),\omega\bigr)\,d\mu(\omega)
=\argmin_{v\in U_c}\bigl\langle\Psi_\mu[\hat u](t),v\bigr\rangle
\quad\text{for a.e. } t\in[0,T],
$$
which in general is stronger than \eqref{VI}; see \cite{BettiolKhalil} and, for
control-affine dynamics with uncertain initial data,
\cite{AronnaMonteiroSierra2025}, where it is stated as a maximum principle with
the opposite sign convention for the adjoint state. This is an inclusion: for
box constraints it gives the bang--bang law $\hat u_i(t)\in\{a_i,b_i\}$
according to the sign of the switching function $\Psi_{\mu,i}[\hat u](t)$, and no
information on singular arcs, where $\Psi_{\mu,i}[\hat u]\equiv0$, which are
characterized via Lie brackets in \cite{AronnaMonteiroSierra2025}; stationary
controls may thus be discontinuous. For $\alpha>0$, instead, the regularized
Hamiltonian is strongly convex in $v$ and selects the unique point
\eqref{projection.formula}, which makes Corollary  \ref{regularity.stationary} and
the second-order analysis of Section  \ref{sec:sosc} possible.
\end{remark}

The projection formula \eqref{projection.formula} has two useful consequences.
The first one is a gain of regularity: the right-hand side involves $\hat u$
only through the state and the adjoint state, which are one derivative smoother
than the control, and $\Proj_{U_c}$ is nonexpansive, so $\hat u$ inherits their
Lipschitz regularity. 
\begin{corollary}[Regularity of stationary controls]
\label{regularity.stationary}
Suppose that Assumptions  \ref{A1}-- \ref{A2} hold and let $\hat u\in\Uad$ be a stationary
point of \eqref{Pmu}. Then $\hat u$ admits a Lipschitz continuous
representative on $[0,T]$, still denoted by $\hat u$, and
$
\hat u\in W^{1,\infty}([0,T];\mathbb R^m)
$ with $
\|\dot{\hat u}\|_{L^\infty}
\le
\frac{C_\Psi}{\alpha}\,
\|L_g\|_{L^1(\mu,\Omega)},
$
where $C_\Psi$ is given in Lemma  \ref{gradient.regularity}.
\end{corollary}

\begin{proof}
By Lemma  \ref{gradient.regularity}, $\Psi_\mu[\hat u]=\nabla J(\hat u)$ is
Lipschitz continuous on $[0,T]$ with Lipschitz constant
$C_\Psi\|L_g\|_{L^1(\mu,\Omega)}$. Since $\Proj_{U_c}$ is nonexpansive
\cite[Proposition 5.3]{brezis2010functional}, the function
$t\mapsto\Proj_{U_c}\bigl(-\alpha^{-1}\Psi_\mu[\hat u](t)\bigr)$ is Lipschitz
continuous with constant $\alpha^{-1}C_\Psi\|L_g\|_{L^1(\mu,\Omega)}$, and by
Theorem  \ref{necessary.conditions} it coincides with $\hat u$ a.e. Hence $\hat u$
has a Lipschitz representative, which is absolutely continuous with derivative
bounded by the same constant.
\end{proof}

\begin{remark}\label{remark.uniqueness.small.alpha}
The second consequence of \eqref{projection.formula} is that the stationarity
system is a fixed-point equation for the map
$u\mapsto\Proj_{U_c}(-\alpha^{-1}\Psi_\mu[u])$ on $\Uad$. Since, by
Proposition  \ref{traj.lipschitz} and Lemma  \ref{p.lc.w.t.u}, the map
$u\mapsto\Psi_\mu[u]$ is Lipschitz continuous from $(\mathcal O,\|\cdot\|_{L^2})$
into $C([0,T];\mathbb R^m)$ with some constant $L_\Psi>0$, the map above is a
contraction of the closed set $\Uad\subset\Vspace$ whenever $\alpha>\sqrt T\,L_\Psi$. In this case
\eqref{Pmu} has a unique stationary point, which is therefore its unique global
minimizer.
\end{remark}

\section{Second-order optimality conditions}\label{sec:sosc}

By Lemma  \ref{2orderw},  the
functional $\Jreg$ is of class $C^2$ on the open set $\mathcal O$
with respect to the $L^2$-topology, with Lipschitz continuous gradient and
Hessian, and by Proposition  \ref{hessian.HS} its Hessian is a compact
perturbation of $\alpha\,\mathrm{Id}$. As explained in
Remark  \ref{extension.hessian}, the two-norm discrepancy typical of nonlinear
control problems is therefore absent. This allows us to apply the abstract
second-order theory of \cite{casas2012,CasTro2012} collected in
Appendix  \ref{appendixB}.

For $\hat u\in\Uad$ we consider the cones \eqref{CT.cones} introduced in Appendix  \ref{appendixB} with
$\mathcal U=\Uad$ and $\mathfrak J=\Jreg$.
Since $\Uad$ is bounded in
$L^\infty([0,T];\mathbb R^m)$, the radial cone $\mathcal R(\hat u)$ is the same
for both choices $U_\infty=\Vspace$ and $U_\infty=L^\infty([0,T];\mathbb R^m)$,
and so are $\mathcal C(\hat u)$ and $\mathcal D(\hat u)$. We also use the
tangent cone
$$
\mathfrak T(\hat u)
:=
\bigl\{\delta u\in\Vspace:\
\exists\,\varepsilon_k\downarrow0,\
\delta u_k\to\delta u \text{ in }\Vspace,\
\hat u+\varepsilon_k\delta u_k\in\Uad\bigr\}.
$$
Since $\Uad$ is convex, $\overline{\mathcal R(\hat u)}^{\,L^2}=\mathfrak T(\hat u)$
(see \cite{BonnansShapiro2000}); in particular, $\mathcal C(\hat u)$ is a closed
convex cone, hence weakly closed in $\Vspace$. We stress that the critical cone
is defined through the derivative of the cost, that is, through
$\nabla\Jreg(\hat u)=\nabla J(\hat u)+\alpha\hat u$.

\begin{remark}[Box constraints]\label{remark.crit_cone_box_constraints}
Let $U_c=\prod_{i=1}^m [a_i,b_i]$. Then the tangent cone admits the pointwise
characterization
$$
\begin{aligned}
\mathfrak T(\hat u)=\Bigg\{\delta u\in \Vspace\ \Bigg | \
&\begin{cases}
\delta u_i(t)\ge0,
& \text{if }\hat u_i(t)=a_i, \\
\delta u_i(t)\le0,
& \text{if }\hat u_i(t)=b_i,\\
\delta u_i(t)\in\mathbb R,
& \text{if }a_i<\hat u_i(t)<b_i,
\end{cases}
\text{ a.e. } t\in[0,T], \ i=1,...,m \Bigg\}.
\end{aligned}
$$
If, moreover, $\hat u$ is a stationary point and
$d:=\nabla\Jreg(\hat u)=\Psi_\mu[\hat u]+\alpha\hat u$, then \eqref{VI} gives
$d_i\ge0$ where $\hat u_i=a_i$, $d_i\le0$ where $\hat u_i=b_i$, and $d_i=0$
where $a_i<\hat u_i<b_i$. Hence $d_i\,\delta u_i\ge0$ a.e. for every
$\delta u\in\mathfrak T(\hat u)$, and the condition
$D\Jreg(\hat u)[\delta u]=\int_0^T\langle d,\delta u\rangle\,dt=0$ forces
$d_i\,\delta u_i=0$ a.e. Therefore
$$
\mathcal C(\hat u)
=
\Bigl\{
\delta u\in\mathfrak T(\hat u):\ \delta u_i(t)=0 \text{ if } d_i(t)\neq0,\
\text{ for a.e. } t\in[0,T],\ i=1,\dots,m
\Bigr\}.
$$
By Corollary  \ref{regularity.stationary}, $\hat u$ and $d$ admit continuous
representatives; with this choice, the active sets
$\{t:\hat u_i(t)=a_i\}$ and $\{t:\hat u_i(t)=b_i\}$ are closed, while
$\{t:a_i<\hat u_i(t)<b_i\}$ and $\{t:d_i(t)\neq0\}$ are open.
\end{remark}

Throughout this section, for $\hat u\in\mathcal O$ we write
$$
Q(\delta u):=D^2\Jreg(\hat u)[\delta u,\delta u]
=
\alpha\|\delta u\|_{L^2}^2
+
\langle\mathcal H_J(\hat u)\delta u,\delta u\rangle_{L^2},
\ \  \delta u\in\Vspace .
$$
The following elementary consequence of the decomposition \eqref{H.def} is
the only structural property of the Hessian used in the sequel.

\begin{lemma}\label{weak.cont.lemma}
Suppose that Assumptions  \ref{A1}-- \ref{A2} hold, let $\hat u\in\mathcal O$ and let
$\delta u_k\rightharpoonup\delta u$ weakly in $\Vspace$. Then
\begin{equation}\label{weak.cont.compact}
\langle\mathcal H_J(\hat u)\delta u_k,\delta u_k\rangle_{L^2}
\rightarrow
\langle\mathcal H_J(\hat u)\delta u,\delta u\rangle_{L^2}.
\end{equation}
In particular, $Q$ is a Legendre form, that is, it is weakly lower
semicontinuous on $\Vspace$ and, if in addition $Q(\delta u_k)\to Q(\delta u)$,
then $\delta u_k\to\delta u$ strongly in $\Vspace$.
\end{lemma}

\begin{proof}
By Proposition  \ref{hessian.HS}, $\mathcal H_J(\hat u)$ is compact, so
$\mathcal H_J(\hat u)\delta u_k\to\mathcal H_J(\hat u)\delta u$ strongly in
$\Vspace$; since $\{\delta u_k\}$ is bounded and converges weakly, the pairing
of the two sequences converges, which is \eqref{weak.cont.compact}. The weak
lower semicontinuity of $Q$ follows from \eqref{weak.cont.compact} and that of
$\delta u\mapsto\alpha\|\delta u\|_{L^2}^2$. If moreover
$Q(\delta u_k)\to Q(\delta u)$, then by \eqref{weak.cont.compact}
$
\alpha\|\delta u_k\|_{L^2}^2
=Q(\delta u_k)-\langle\mathcal H_J(\hat u)\delta u_k,\delta u_k\rangle_{L^2}
\rightarrow\alpha\|\delta u\|_{L^2}^2 ,
$
so that $\|\delta u_k\|_{L^2}\to\|\delta u\|_{L^2}$, which together with
$\delta u_k\rightharpoonup\delta u$ gives strong convergence in the Hilbert
space $\Vspace$.
\end{proof}

As a consequence, on the critical cone positivity and coercivity of $Q$ are
equivalent.

\begin{lemma}[Positivity implies coercivity on the critical cone]
\label{coercivity.cone}
Suppose that Assumptions  \ref{A1}-- \ref{A2} hold and let $\hat u\in\Uad$. If
$Q(\delta u)>0$ for all $\delta u\in\mathcal C(\hat u)\setminus\{0\}$, then
there exists $\beta>0$ such that
$
Q(\delta u)\ge\beta\|\delta u\|_{L^2}^2
$ for all $\delta u\in\mathcal C(\hat u).
$
\end{lemma}
\begin{proof}
By Lemma  \ref{weak.cont.lemma}, $Q$ is a Legendre form on the Hilbert space
$\Vspace$, and $\mathcal C(\hat u)$ is a closed convex cone. The claim is
therefore \cite[Lemma 3.75]{BonnansShapiro2000}, applied with the constant
family of cones $C_\eta:=\mathcal C(\hat u)$, for which the weak-closedness
condition of that lemma holds trivially.
\end{proof}
We now state the second-order necessary conditions, which are the particular
case $U_\infty=L^\infty([0,T];\mathbb R^m)$, $U_2=\Vspace$ of
Theorem  \ref{CT.necessary}.

\begin{theorem}[Second-order necessary conditions]\label{son.conditions}
Suppose that Assumptions  \ref{A1}-- \ref{A2} hold and that the compact convex set $U_c$ is
polyhedral, that is,
$U_c=\{v\in\mathbb R^m:\ \langle a_j,v\rangle\le b_j,\ j=1,\dots,n_c\}$ for some
$n_c\in\mathbb N$, $a_j\in\mathbb R^m$ and $b_j\in\mathbb R$. If $\hat u\in\Uad$
is a local minimizer of \eqref{Pmu} in the $L^\infty$-topology, then
$
Q(\delta u)\ge0
$ for all $\delta u\in\mathcal C(\hat u).
$
\end{theorem}
\begin{proof}
We apply Theorem  \ref{CT.necessary}. Assumption {\bf (A1)} of
Appendix  \ref{appendixB} holds
by Lemma  \ref{2orderw} and Remark  \ref{extension.hessian}. Since $\hat u$ is a
local minimizer, it is a stationary point, so that $-\nabla\Jreg(\hat u)\in N_{\Uad}(\hat u)$. The regularity
condition $\mathcal C(\hat u)=\overline{\mathcal D(\hat u)}^{\,L^2}$ is then the
polyhedricity of $\Uad$ at $\hat u$ in the sense of
\cite[Definition 3.51]{BonnansShapiro2000}, and sets of the form $\Uad$ defined
by finitely many pointwise linear constraints are polyhedric in $\Vspace$ by
\cite[Proposition 6.35]{BonnansShapiro2000}, applied with $s=2.$ Finally, since $\Uad$ is bounded
in $L^\infty([0,T];\mathbb R^m)$, its radial directions are essentially
bounded, so the cones $\mathcal R(\hat u)$ and $\mathcal D(\hat u)$ are the same
whether computed in $L^\infty([0,T];\mathbb R^m)$ or in $\Vspace$, and the
claim follows.
\end{proof}

The sufficient condition is the particular case $U_\infty=U_2=\Vspace$ of
Theorem  \ref{CT.growth}.

\begin{theorem}[Quadratic growth]\label{sos.conditions}
Suppose that Assumptions  \ref{A1}-- \ref{A2} are satisfied. Let $\hat u\in\Uad$ be a
stationary point, and assume that
$
Q(\delta u)>0$
 for all $ \delta u\in
\mathcal C(\hat u)\setminus\{0\}.
$
Then there exist $\varepsilon>0$ and $\delta>0$ such that
\begin{equation}\label{quadratic.growth}
\Jreg(\hat u)
+
\frac{\delta}{2}
\|u-\hat u\|_{L^2}^2
\le
\Jreg(u)
\ \
\text {for all }
u\in\Uad \cap
B_{L^2}(\hat u,\varepsilon).
\end{equation}
In particular, $\hat u$ is a strict local minimizer of \eqref{Pmu} in the
$L^2$-topology.
\end{theorem}
\begin{proof}
We apply Theorem  \ref{CT.growth}. The first-order condition is \eqref{VI}.
Assumption {\bf (A1)} is Lemma  \ref{2orderw}. For {\bf (A2)}, let
$\{u_k\}\subset\Uad$ with $u_k\to\hat u$ in $\Vspace$ and
$\delta u_k\rightharpoonup\delta u$. Since $\nabla\Jreg$ and $\mathcal H_J$ are
Lipschitz continuous (Proposition  \ref{1order}(iii) and
Lemma  \ref{2orderw}(ii)) and $\{\delta u_k\}$ is bounded, we get
$D\Jreg(u_k)[\delta u_k]\to D\Jreg(\hat u)[\delta u]$, which is
\eqref{CT.2.3}, and, by \eqref{weak.cont.compact},
$$
D^2\Jreg(u_k)[\delta u_k,\delta u_k]
=\alpha\|\delta u_k\|_{L^2}^2
+\langle\mathcal H_J(\hat u)\delta u,\delta u\rangle_{L^2}+o(1).
$$
The weak lower semicontinuity of the norm yields \eqref{CT.2.4}, and, if
$\delta u=0$, the same identity yields \eqref{CT.2.5} with $\vartheta=\alpha$.
\end{proof}

\begin{remark}\label{remark.no.gap}
When $U_c$ is a polytope, Theorems  \ref{son.conditions} and
 \ref{sos.conditions} form a pair of second-order conditions \emph{with no gap},
as in finite dimension: the necessary condition is $Q\ge0$ on
$\mathcal C(\hat u)$, while the sufficient one is $Q>0$ on
$\mathcal C(\hat u)\setminus\{0\}$, which by Lemma  \ref{coercivity.cone} is
equivalent to $Q\ge\beta\|\cdot\|_{L^2}^2$ on $\mathcal C(\hat u)$. 
When $U_c$ is not a polytope, for instance a Euclidean ball in $\mathbb R^m$ with $m\ge2$, the set $\Uad$ is in general not polyhedric, and the necessary condition must be corrected by a curvature term of the constraint set, the
so-called sigma term; see \cite[Theorem 3.45 and \S3.2.3]{BonnansShapiro2000}.
\end{remark}

\begin{remark}\label{remark.alpha.large}
 If $\alpha>\|\mathcal H_J(\hat u)\|_{\mathcal L(L^2)}$, then
$
Q(\delta u)\ge\bigl(\alpha-\|\mathcal H_J(\hat u)\|_{\mathcal L(L^2)}\bigr)\|\delta u\|^2_{L^2}$  for all $\delta u\in\Vspace,
$
so that the second-order sufficient condition holds on the whole space. By
Proposition  \ref{hessian.HS} and \eqref{HS.bound},
$\|\mathcal H_J(\hat u)\|_{\mathcal L(L^2)}\le c_{\HS}\|L_g\|_{L^1(\mu,\Omega)}$,
so this is guaranteed as soon as $\alpha>c_{\HS}\|L_g\|_{L^1(\mu,\Omega)}$, a
condition involving only the data of the problem.
\end{remark}

\section{Convergence of the SAA scheme}\label{sec:quantitative}
\subsection{Consistency}\label{subsec:consistency} In this subsection 
we establish the almost sure consistency of the SAA scheme, both for the
optimal values and, in the uniform norm, for the empirical minimizers; these
qualitative results are the only ones needed in Section  \ref{sec:asymptotics}.

Our discussion of convergence starts at order zero, that is, with the
qualitative notion of convergence of variational problems known as
\textit{epi-convergence}, also called
$\Gamma$-convergence.

\begin{definition}[Epi-convergence, see \cite{AttouchWets1983}]
Let $(\mathcal{U}, d)$ be a Polish space, and consider a sequence of lower
semicontinuous functions $F_M: \mathcal{U}\rightarrow (-\infty, \infty]$. The
sequence $\{F_M\}$ is said to epi-converge to $F$, denoted
$F_M\xrightarrow[M \to \infty]{\mathrm{epi}}F$, if the following two
conditions are satisfied:
\begin{enumerate}
    \item[(1)] whenever $\displaystyle\lim_{M \to \infty}u_M = u$, we have
    $\displaystyle \liminf_{M \to \infty} F_M(u_M) \geq F(u)$;
    \item[(2)] for every $u \in \mathcal{U}$, there exists a sequence
    $u_M \rightarrow u$ such that
    $\displaystyle \lim_{M \to \infty} F_M(u_M) = F(u)$.
\end{enumerate}
\end{definition}

\begin{definition}[see {\cite{AttouchWets1983}}]
Let $(\mathcal{U},d)$ be a separable complete metric space with Borel
$\sigma$-field $\Bor(\mathcal{U})$, and let $(\Omega,\mathcal{A},\mu)$
be a probability space with $\mathcal{A}$ $\mu$-complete. A function
$g : \mathcal{U} \times \Omega \to \mathbb{R}$ is a {\it random lower
semicontinuous} function (or {\it normal integrand}) if $u \mapsto g(u,\omega)$
is lower semicontinuous for all $\omega \in \Omega$ and $g$ is
$\Bor(\mathcal{U})\otimes\mathcal{A}$-measurable.
\end{definition}

\begin{remark}\label{gt.rlsc}
Lemma  \ref{gTlipschitz} shows that the function $g_T$,
restricted to $\Uad$, is a random lower semicontinuous function that satisfies
the hypotheses of \cite[Theorem 2.3]{ArtsteinWets1995}; it is
continuous in $u$, jointly measurable, and bounded from below by $-a$ with
$a\in L^1(\mu,\Omega)$.
\end{remark}

\begin{theorem}[Epi-convergence of the SAA functionals]\label{JM.epi.J}
Under Assumptions  \ref{A1}-- \ref{A2}, $J_M \xrightarrow{\text{epi}} J$ and
$\Jreg_M \xrightarrow{\text{epi}} \Jreg$ almost surely on $\Uad$ as
$M \to \infty$.
\end{theorem}

\begin{proof}
The set $\Uad$ is a closed subset of the separable Hilbert space $\Vspace$,
hence a Polish space. Replacing $\Bor(\Omega)$ by its $\mu$-completion if
necessary, which does not affect almost sure statements, the epi-convergence of
$J_M$ to $J$ follows from Remark  \ref{gt.rlsc} and
 \cite[Theorem 2.3]{ArtsteinWets1995}, together with an argument
analogous to that used in \cite[Theorem 3.7]{Phelps2013}. The statement for
$\Jreg_M$ follows because the term $\tfrac\alpha2\|\cdot\|_{L^2}^2$ is
deterministic, continuous and independent of $M$, so that both the $\liminf$
inequality and the existence of a recovery sequence are preserved.
\end{proof}

The major interest of epi-convergence is that it provides consistency
properties for sequences of minimizers: if $\hat u$ is a cluster point of a
sequence of minimizers of $\Jreg_M$, then it minimizes $\Jreg$.

\begin{proposition}[Consistency of empirical minimizers]\label{consistency.uM}
Suppose that Assumptions  \ref{A1}-- \ref{A2} hold. Almost surely, if, along a
subsequence $\{M_k\}$, $\hat u_{M_k}\to \hat u$ strongly in $\Vspace$, then
$$
\hat u\in\argmin_{u\in \Uad} \Jreg(u)
\qquad\text{and}\qquad
\lim_{k\to\infty}V_{M_k}=V.
$$
\end{proposition}

\begin{proof}
On the full-measure set where the epi-convergence of
Theorem  \ref{JM.epi.J} holds, each $\Jreg_M$ is continuous on $\Uad$, so
\cite[Theorem 2.5]{AttouchWets1983} applies.
\end{proof}

The existence of a convergent subsequence is guaranteed by the regularity of
stationary controls given by Corollary  \ref{regularity.stationary}. For $\Lambda>0$ we consider the set of admissible controls
with Lipschitz constant at most $\Lambda$,
\begin{equation}\label{K.Lambda}
\mathcal K_\Lambda
:=
\bigl\{u\in\Uad:\ \|u(t)-u(s)\|\le\Lambda|t-s|\ \text{ for all } s,t\in[0,T]\bigr\},
\end{equation}
where, as usual, $u$ is identified with its Lipschitz representative. By the
Arzel\`a--Ascoli theorem, $\mathcal K_\Lambda$ is compact in
$C([0,T];\mathbb R^m)$, hence in $L^\infty([0,T];\mathbb R^m)$ and in
$\Vspace$; indeed, uniform limits of $\Lambda$-Lipschitz functions with values
in the closed set $U_c$ are again $\Lambda$-Lipschitz with values in $U_c$.

The sets $\mathcal K_\Lambda$ also provide the measurability of empirical
minimizers.

\begin{lemma}[Measurable selection of empirical minimizers]\label{measurable.selection}
Suppose that Assumptions  \ref{A1}-- \ref{A2} hold and let $M\in\mathbb N$. Then
there exists an $\mathcal F$-measurable map
$\hat u_M:\Omega_{\mathbb P}\to\Vspace$ such that, pointwise,
$\hat u_M\in\argmin_{u\in\Uad}\Jreg_M(u)$.
\end{lemma}

\begin{proof}
Set
$\Lambda_M:=\alpha^{-1}C_\Psi\frac1M\sum_{i=1}^ML_g(X_i)
=\alpha^{-1}C_\Psi\|L_g\|_{L^1(\mu_M,\Omega)}$, which is a
random variable. By Theorem  \ref{existence}, the set
$\argmin_{u\in\Uad}\Jreg_M(u)$ is nonempty, and by
Corollary  \ref{regularity.stationary}, applied with the measure
$\mu_M$, it is contained in $\mathcal K_{\Lambda_M}$.
Hence
$$
\argmin_{u\in\Uad}\Jreg_M(u)
=\argmin_{u\in\mathcal K_{\Lambda_M}}\Jreg_M(u)
\ \text{ pointwise on }\Omega_{\mathbb P} .
$$
The random multifunction $\Gamma:=\mathcal K_{\Lambda_M}$ has
nonempty compact values in the separable metric space $\Vspace$ (it contains
the constant controls with values in $U_c$). It is weakly measurable: given an
open set $G\subset\Vspace$, since $\Lambda\mapsto\mathcal K_\Lambda$ is
nondecreasing, the set $I_G:=\{\Lambda\ge0:\mathcal K_\Lambda\cap G\neq\emptyset\}$
is either empty or an interval unbounded from above, so that
$\{\Gamma\cap G\ne\emptyset\}=\{\Lambda_M\in I_G\}$ is measurable. Moreover,
$\Jreg_M$ is a Carath\'eodory function: it is continuous in $u$ by
Lemma  \ref{gTlipschitz}(i) and measurable in its random argument by the measurability part of
Lemma  \ref{gTlipschitz}(i). The Measurable Maximum Theorem
\cite[Theorem 18.19]{aliprantis2006infinite}, applied to $-\Jreg_M$, shows that
$\argmin_{u\in\Gamma}\Jreg_M(u)$ admits a measurable selector.
\end{proof}

\begin{proposition}[Uniform consistency of empirical minimizers]
\label{uniform.consistency.Linfty}
Suppose that Assumptions  \ref{A1}-- \ref{A2} hold and that $\Jreg$ admits a unique global
minimizer $\hat u\in\Uad$. Then
$
\|\hat u_M-\hat u\|_{L^\infty}\xrightarrow[]{}0$ as $M \to \infty$ a.s.
\end{proposition}
\begin{proof}
Set $\Lambda:=2\alpha^{-1}C_\Psi\|L_g\|_{L^1(\mu,\Omega)}$ and consider the event
\begin{equation}\label{BM.def}
B_M:=\bigl\{\|L_g\|_{L^1(\mu_M,\Omega)}>2\|L_g\|_{L^1(\mu,\Omega)}\bigr\}.
\end{equation}
Outside $B_M$, every empirical minimizer is a stationary point of the problem
associated with $\mu_M$, so that Corollary  \ref{regularity.stationary},
applied with the measure $\mu_M$, yields $\argmin_{\Uad}\Jreg_M\subset\mathcal K_\Lambda$.
By the strong law of large numbers,
$\|L_g\|_{L^1(\mu_M,\Omega)}=\frac1M\sum_{i=1}^ML_g(X_i)\to\|L_g\|_{L^1(\mu,\Omega)}$
almost surely. Hence, almost surely, $B_M$ occurs only for finitely many $M$.
Fix a realization in the full-measure set where this holds, and the
epi-convergence of Theorem  \ref{JM.epi.J} holds, and let $M_0$ be such that
$B_M$ does not occur for $M\ge M_0$. Then $\hat u_M\in\mathcal K_\Lambda$ for
$M\ge M_0$. Let $\{\hat u_{M_k}\}$ be any subsequence. By compactness of
$\mathcal K_\Lambda$ in $C([0,T];\mathbb R^m)$, a further subsequence converges
uniformly, and hence in $\Vspace$, to some $\bar u\in\mathcal K_\Lambda\subset\Uad$,
which minimizes $\Jreg$ by Proposition  \ref{consistency.uM}. By uniqueness,
$\bar u=\hat u$. Since every subsequence has a further subsequence converging
uniformly to $\hat u$, the whole sequence converges uniformly to $\hat u$.
\end{proof}
We finally record a global version of the quadratic growth condition on
compact subsets of $\Uad$, which will be used below to convert estimates on
the optimal values into estimates on the minimizers.

\begin{lemma}[Quadratic growth on compact sets]\label{global.growth}
Suppose that Assumptions  \ref{A1}-- \ref{A2} hold, that $\Jreg$ admits a unique global
minimizer $\hat u\in\Uad$, and that $\hat u$ satisfies the assumptions of
Theorem  \ref{sos.conditions}, with constants $\varepsilon,\delta>0$ in
\eqref{quadratic.growth}. Then, for every compact  $\mathcal K\subset\Uad$, there exists $\kappa_{\mathcal K}>0$ such that
$
\kappa_{\mathcal K}\|u-\hat u\|_{L^2}^2\le\Jreg(u)-\Jreg(\hat u)$
for all $u\in\mathcal K .
$
\end{lemma}

\begin{proof}
The set $\mathcal K_\varepsilon:=\{u\in\mathcal K:\|u-\hat u\|_{L^2}\ge\varepsilon\}$
is compact and $\Jreg$ is continuous on $\Uad$ by Lemma  \ref{gTlipschitz}(iii);
if $\mathcal K_\varepsilon\neq\emptyset$, the minimum value $\eta$ of
$\Jreg-\Jreg(\hat u)$ over $\mathcal K_\varepsilon$ is attained, and $\eta>0$ by
uniqueness of the global minimizer. Since $\|u-\hat u\|_{L^2}\le2\sqrt T\rho$ for
all $u\in\Uad$, the claim holds with
$\kappa_{\mathcal K}:=\min\{\delta/2,\eta/(4T\rho^2)\}$, by
\eqref{quadratic.growth} on $B_{L^2}(\hat u,\varepsilon)$ and by the definition
of $\eta$ on $\mathcal K_\varepsilon$.
\end{proof}

\subsection{Convergence rates}\label{subsec:rates}

As discussed in the introduction, epi-convergence does not provide rates. We
obtain them by a metric entropy argument on the compact set
$\mathcal K_\Lambda$ (Proposition  \ref{prop.conv_rates_energy}), in the spirit
of \cite{melnikov2026convergence}.

The following estimate on the metric entropy of the sets $\mathcal K_\Lambda$
is used in Proposition  \ref{prop.conv_rates_energy} and in Section  \ref{sec:asymptotics}. We denote by
$N(\varepsilon,\mathcal K,\|\cdot\|)$ the minimal number of balls of radius
$\varepsilon$ in the norm $\|\cdot\|$ needed to cover $\mathcal K$, the centers
being allowed outside $\mathcal K$. This result quantifies the classical entropy estimate for Lipschitz
functions \cite{KolmogorovTikhomirov1959} (see also
\cite[Theorem 2.7.1]{vanderVaartWellner1996}, without the logarithmic factor
but with a non-explicit constant); its short proof gives explicit constants and
centers in $\mathcal O$. 
\begin{lemma}[Metric entropy of $\mathcal K_\Lambda$]\label{entropy.K}
For every $\Lambda>0$ and every $\varepsilon>0$,
$$
\log N\bigl(\varepsilon,\mathcal K_\Lambda,\|\cdot\|_{L^\infty}\bigr)
\le
m\Bigl(1+\frac{2\Lambda T}{\varepsilon}\Bigr)\log\Bigl(1+\frac{4\rho}{\varepsilon}\Bigr),
$$
and the centers can be chosen in $\mathcal O$. Consequently,
$N(\varepsilon,\mathcal K_\Lambda,\|\cdot\|_{L^2})
\le N(\varepsilon/\sqrt T,\mathcal K_\Lambda,\|\cdot\|_{L^\infty})$.
\end{lemma}

\begin{proof}
Fix $\varepsilon>0$. We first construct a finite discretization of the
time interval $[0,T]$. Set
$
h:=\frac{\varepsilon}{2\Lambda},
$
$n:=\left\lceil\frac{T}{h}\right\rceil,$
$t_j:=\min\{jh,T\},$ $ j=0,\ldots,n.
$
Then $t_0=0$, $t_n=T$, and
$
t_{j+1}-t_j\leq h,
$ $ j=0,\ldots,n-1.
$
Furthermore,
$
n=\left\lceil\frac{2\Lambda T}{\varepsilon}\right\rceil
\leq 1+\frac{2\Lambda T}{\varepsilon}.
$
Next, let $Y\subset\overline{B(0,\rho)}$ be an $\varepsilon/2$-net of this
ball in $\mathbb{R}^m$. By a standard volumetric
argument, $Y$ can be chosen so that
$
\operatorname{card}(Y)\leq \left(1+\frac{4\rho}{\varepsilon}\right)^m.
$
For each $u\in\mathcal{K}_{\Lambda}$ and each $j=0,\ldots,n-1$,
choose $y_j\in Y$ such that
$
\|u(t_j)-y_j\|\leq \frac{\varepsilon}{2}.
$
Define the step function $w:[0,T]\to\mathbb{R}^m$ by
$
w(t):=y_j,
$ $t\in[t_j,t_{j+1}),$ $ j=0,\ldots,n-1.
$
The value at $t=T$ can be assigned arbitrarily in $Y$, since it does
not affect the $L^\infty$ or $L^2$ equivalence class. If
$t\in[t_j,t_{j+1})$, the Lipschitz continuity of $u$ yields
$$
\begin{aligned}
\|u(t)-w(t)\|
&\leq \|u(t)-u(t_j)\|+\|u(t_j)-y_j\|
\leq \Lambda |t-t_j|+\frac{\varepsilon}{2}\leq \Lambda h+\frac{\varepsilon}{2}
=\varepsilon.
\end{aligned}
$$
Hence,
$
\|u-w\|_{L^\infty([0,T];\mathbb R^m)}\leq\varepsilon.
$
Thus, the collection of all step functions constructed in this way
forms an $\varepsilon$-net of $\mathcal{K}_{\Lambda}$ in
$L^\infty([0,T];\mathbb R^m)$. There are at most $\operatorname{card}(Y)$ possible choices for each of the $n$ values
$y_0,\ldots,y_{n-1}$. Therefore, the number of such step functions is
bounded above by
$
\operatorname{card}(Y)^n
\leq
\left(1+\frac{4\rho}{\varepsilon}\right)^{mn}.
$
It follows that
$$
\begin{aligned}
\log N\bigl(\varepsilon,\mathcal{K}_{\Lambda},\|\cdot\|_{L^\infty}\bigr)
&\leq mn\log\left(1+\frac{4\rho}{\varepsilon}\right)\leq
m\left(1+\frac{2\Lambda T}{\varepsilon}\right)
\log\left(1+\frac{4\rho}{\varepsilon}\right),
\end{aligned}
$$
which proves the entropy estimate. We now verify that the covering centers belong to $\mathcal{O}$.
Since $y_j\in\overline{B(0,\rho)}$ for every $j$, each constructed
step function satisfies
$
\|w\|_{L^2}^2
=\sum_{j=0}^{n-1}(t_{j+1}-t_j)\|y_j\|^2
\leq T\rho^2.
$
Consequently,
$
\|w\|_{L^2}
\leq \sqrt{T}\rho<r_\rho.
$ Thus, if $\mathcal{O}$ contains the $L^2$-ball of radius $r_\rho$,
all the covering centers belong to $\mathcal{O}$. Finally, since $
\|\cdot\|_{L^2}
\leq \sqrt{T}\,
\|\cdot\|_{L^\infty},$
 any $L^\infty$-cover of $\mathcal{K}_{\Lambda}$ with
radius $\varepsilon/\sqrt{T}$ is also an $L^2$-cover with radius
$\varepsilon$. This gives
$N(\varepsilon,\mathcal K_\Lambda,\|\cdot\|_{L^2})
\le N(\varepsilon/\sqrt T,\mathcal K_\Lambda,\|\cdot\|_{L^\infty})$ 
and completes the proof.
\end{proof}
The next result is a variant of \cite[Theorem 6.1]{melnikov2026convergence},
which covers general strongly convex regularizers but needs Lipschitz constants
uniform in the uncertain parameter for its sub-Gaussian maximal inequality. We
use instead a bracketing inequality. The bracketing number
$N_{[\,]}(\varepsilon,\mathscr G,L^2(\mu))$ of a class
$\mathscr G\subset L^2(\mu,\Omega)$ is the minimal number of brackets
$[\ell,r]:=\{f:\ \ell\le f\le r\}$, with $\ell,r\in L^2(\mu,\Omega)$ and
$\|r-\ell\|_{L^2(\mu,\Omega)}\le\varepsilon$, needed to cover $\mathscr G$
\cite[Definition 2.1.6]{vanderVaartWellner1996}. Since only the width $r-\ell$
must be small in $L^2(\mu)$, it suffices that $L_g$ be square integrable, as in
Assumption  \ref{A2}.

Since $\Jreg_M$ is continuous in $u$ and measurable in $\omega_{\mathbb P}$ (Lemma  \ref{gTlipschitz}(i)), $V_M$ is an infimum over a countable dense subset of $\Uad$, hence measurable, and the expectation below is well defined.

\begin{proposition}\label{prop.conv_rates_energy}
Suppose that Assumptions  \ref{A1}-- \ref{A2} hold. Then there exists a constant $C>0$,
independent of $M$ and of $\alpha$, such that
\begin{equation}\label{rate.energy}
\mathbb E\bigl[|V_M-V|\bigr]\le\frac{C}{\sqrt M}\Bigl(1+\frac{1}{\sqrt\alpha}\Bigr)
\ \text{ for all } M\in\mathbb N .
\end{equation}
\end{proposition}
\begin{proof}
Set $\Lambda:=2\alpha^{-1}C_\Psi\|L_g\|_{L^1(\mu,\Omega)}$, $\mathcal K:=\mathcal K_\Lambda$,
$D_\rho:=2\sqrt T\rho$, and let $a=|g(0,\cdot)|+L_gR\in L^2(\mu,\Omega)$ be as in
Lemma  \ref{gTlipschitz}(ii). We write $\|\cdot\|_{2}:=\|\cdot\|_{L^2(\mu,\Omega)}$.

We first establish a maximal inequality over $\mathcal K$. Fix $u_0\in\mathcal K$ and consider
the class $\mathscr G_{\mathcal K}:=\{g_T(u,\cdot)-g_T(u_0,\cdot):u\in\mathcal K\}$, which has
envelope $D_\rho L_T$ by Lemma  \ref{gTlipschitz}(i). Let $\{u_k\}\subset\mathcal O$ be an
$\varepsilon$-net of $\mathcal K$ (Lemma  \ref{entropy.K}). By Lemma  \ref{gTlipschitz}, the brackets
$[g_T(u_k,\cdot)-g_T(u_0,\cdot)\mp\varepsilon L_T]$ cover $\mathscr G_{\mathcal K}$ and have
$L^2(\mu)$-width $2\varepsilon\|L_T\|_2$, so that
$N_{[\,]}(2\varepsilon\|L_T\|_2,\mathscr G_{\mathcal K},L^2(\mu))\le N(\varepsilon,\mathcal K,\|\cdot\|_{L^2})$
(cf.\ \cite[Theorem 2.7.11]{vanderVaartWellner1996}), and this number equals one for $\varepsilon\ge D_\rho$.
Since $J_M-J$ is continuous on the separable set $\mathcal K$, its supremum equals the supremum over a
countable dense subset and is therefore measurable, so the bracketing maximal inequality
\cite[Corollary 19.35]{vanderVaart1998} applies with the usual expectation. After the change of
variable $\varepsilon\mapsto2\varepsilon\|L_T\|_2$ it yields
$$
\mathbb E\Bigl[\sup_{u\in\mathcal K}\bigl|(J_M-J)(u)-(J_M-J)(u_0)\bigr|\Bigr]
\le\frac{c_0\|L_T\|_2}{\sqrt M}\int_0^{D_\rho}\sqrt{\log N(\varepsilon,\mathcal K,\|\cdot\|_{L^2})}\,d\varepsilon .
$$
for some $c_0>0.$ By Lemma  \ref{entropy.K}, $\log N(\varepsilon,\mathcal K,\|\cdot\|_{L^2})\le
m\bigl(1+\tfrac{2\Lambda T^{3/2}}{\varepsilon}\bigr)\log\bigl(1+\tfrac{4\sqrt T\rho}{\varepsilon}\bigr)$,
and since $\sqrt{x+y}\le\sqrt x+\sqrt y$,
$$
\int_0^{D_\rho}\!\sqrt{\log N(\varepsilon,\mathcal K,\|\cdot\|_{L^2})}\,d\varepsilon\le\sqrt m\int_0^{D_\rho}\Bigl(1+\sqrt{\tfrac{2\Lambda T^{3/2}}{\varepsilon}}\Bigr)\sqrt{\log\Bigl(1+\tfrac{4\sqrt T\rho}{\varepsilon}\Bigr)}\,d\varepsilon=:c_1+c_2\sqrt\Lambda,
$$
where $c_1,c_2<\infty$ depend only on $m$, $T$ and $\rho$, the singularity at $\varepsilon=0$ being
integrable. Using also $\mathbb E|(J_M-J)(u_0)|\le\|g_T(u_0,\cdot)\|_2/\sqrt M\le\|a\|_2/\sqrt M$, we obtain
\begin{equation}\label{maximal.K}
\mathbb E\Bigl[\sup_{u\in\mathcal K}|J_M(u)-J(u)|\Bigr]
\le\frac{\|a\|_2+c_0\|L_T\|_2(c_1+c_2\sqrt\Lambda)}{\sqrt M}.
\end{equation}

It remains to control the event $B_M$ on which the empirical minimizers may leave $\mathcal K$. By
\eqref{BM.def}, $B_M\subset\{|\,\|L_g\|_{L^1(\mu_M,\Omega)}-\|L_g\|_{L^1(\mu,\Omega)}|\ge\|L_g\|_{L^1(\mu,\Omega)}\}$.
Since $\|L_g\|_{L^1(\mu_M,\Omega)}$ is the mean of $M$ i.i.d.\ variables with expectation
$\|L_g\|_{L^1(\mu,\Omega)}$ and variance $\operatorname{Var}_\mu(L_g)/M$, Chebyshev's inequality together
with $\operatorname{Var}_\mu(L_g)\le\|L_g\|_2^2$ and $\|L_g\|_{L^1(\mu,\Omega)}\ge1$ (because $L_g\ge1$) gives
\begin{equation}\label{PBM}
\mathbb P(B_M)\le\frac{\operatorname{Var}_\mu(L_g)}{M\|L_g\|^2_{L^1(\mu,\Omega)}}\le\frac{\|L_g\|_2^2}{M}.
\end{equation}

Now, by Corollary  \ref{regularity.stationary}, applied with $\mu$ and with $\mu_M$, the minimizers of
$\Jreg$ lie in $\mathcal K$, and so do those of $\Jreg_M$ outside $B_M$. Since
$\Jreg_M-\Jreg=J_M-J$, we get $|V_M-V|\le\sup_{\mathcal K}|J_M-J|$ outside $B_M$, whereas on $B_M$ we
use $|V_M-V|\le b_M:=\|a\|_{L^1(\mu,\Omega)}+\frac1M\sum_{i=1}^Ma(X_i)$, which satisfies
$\mathbb E[b_M^2]\le4\|a\|_2^2$. Splitting the expectation according to $B_M$ and its complement, and
applying the Cauchy--Schwarz inequality and \eqref{PBM} on $B_M$,
$$
\mathbb E|V_M-V|\le\mathbb E\Bigl[\sup_{\mathcal K}|J_M-J|\Bigr]+\bigl(\mathbb E[b_M^2]\bigr)^{1/2}\mathbb P(B_M)^{1/2}
\le\mathbb E\Bigl[\sup_{\mathcal K}|J_M-J|\Bigr]+\frac{2\|a\|_2\|L_g\|_2}{\sqrt M}.
$$
Combining this with \eqref{maximal.K},
$$
\mathbb E|V_M-V|\le\frac{\|a\|_2\bigl(1+2\|L_g\|_2\bigr)+c_0\|L_T\|_2(c_1+c_2\sqrt\Lambda)}{\sqrt M}.
$$
Finally, $\sqrt\Lambda=\sqrt{2C_\Psi\|L_g\|_{L^1(\mu,\Omega)}}\,\alpha^{-1/2}$ and no other constant
depends on $\alpha$, which proves the claim.
\end{proof}
\begin{remark}[Rates for the empirical minimizers]\label{remark.rates.controls}
Under the assumptions of Lemma  \ref{global.growth}, outside $B_M$ both $\hat u$
and $\hat u_M$ lie in $\mathcal K=\mathcal K_\Lambda$, so that
$\kappa_{\mathcal K}\|\hat u_M-\hat u\|_{L^2}^2\le\Jreg(\hat u_M)-\Jreg(\hat u)
\le2\sup_{\mathcal K}|J_M-J|$, while $\|\hat u_M-\hat u\|_{L^2}^2\le4T\rho^2$ on
$B_M$. By \eqref{maximal.K} and \eqref{PBM}, this gives
$\mathbb E\|\hat u_M-\hat u\|_{L^2}^2=O(M^{-1/2})$. This non-asymptotic rate is
not sharp: by the proof of Theorem  \ref{asymp.active},
$\sqrt M\|\hat u_M-\hat u\|_{L^2}$ is bounded in probability.
\end{remark}

\section{Asymptotics of the SAA scheme}\label{sec:asymptotics}

In this section we describe the fluctuations of the SAA estimators around their
limits, first for the optimal values and then for the empirical minimizers. The
consistency required in both cases is provided by
Subsection  \ref{subsec:consistency}.

\subsection{Optimal values}
As announced in the introduction, Theorem  \ref{asymp.values} follows the proof of
\cite[Theorem 5.7]{shapiro2021} on the compact set $\mathcal K_\Lambda$. Both
statements were obtained in \cite[Theorem B.2 and Corollary 3.1]{JaveedMilz2026}
without regularization, but for an $\omega$-independent terminal cost and a
compact parameter set, hence with Lipschitz constants uniform in $\omega$.

We first state the following result, which is a central limit theorem for
i.i.d.\ random elements of the space $C(\mathcal K)$ of continuous functions
on a compact metric space, due to Jain and Marcus \cite{JainMarcus1975}; see
also \cite[Corollary 7.17]{AraujoGine1980}.
\begin{theorem}[{\cite[Theorem 1]{JainMarcus1975}}]\label{CLT.JainMarcus}
Let $(\mathcal K,d)$ be a compact metric space, and let
$\upchi_1,\upchi_2,\dots$ be i.i.d.\ random elements of $C(\mathcal K)$ with
$\mathbb E[\upchi_1(u)]=0$ for every $u\in\mathcal K$ and
$\sup_{u\in\mathcal K}\mathbb E[\upchi_1(u)^2]<\infty$. Assume that there exists
a nonnegative random variable $\Upsilon$ with $\mathbb E[\Upsilon^2]<\infty$ such
that
$$
|\upchi_1(u)-\upchi_1(v)|\le\Upsilon\,d(u,v)
\ \text{ for all }u,v\in\mathcal K,
$$
and that
$
\int_0^1\sqrt{\log N(\varepsilon,\mathcal K,d)}\,d\varepsilon<\infty .
$
Then
$
\frac1{\sqrt M}\sum_{i=1}^M\upchi_i\xrightarrow[]{\mathfrak D}Y \ 
\text{ in }C(\mathcal K),
$
where $Y$ is a centered Gaussian process on $\mathcal K$ with continuous paths
and covariance
$\mathbb E[Y(u)Y(v)]=\mathbb E[\upchi_1(u)\upchi_1(v)]$.
\end{theorem}

\begin{theorem}[Asymptotics of the optimal values]\label{asymp.values}
Suppose that Assumptions  \ref{A1}-- \ref{A2} hold,
let $\Lambda:=\frac{2 C_\Psi}{\alpha}\| L_g\|_{L^1(\mu,\Omega)}$, and set
$$
\hat S:=\argmin_{u\in\Uad}\Jreg(u),
\qquad
\sigma^2(u):=\operatorname{Var}_\mu\big[g_T(u,\cdot)\big],\quad u\in\Uad .
$$
Then the following statements hold.
\begin{itemize}
\item[(i)] (Functional central limit theorem.) In the space
$C(\mathcal K_\Lambda)$,
$
\sqrt M\,(J_M-J)\xrightarrow[]{\mathfrak D}Y,
$
where $Y$ is a centered Gaussian process on $\mathcal K_\Lambda$ with
continuous paths and covariance
$
\mathbb E[Y(u)Y(v)]
=
\operatorname{Cov}_\mu\big(g_T(u,\cdot),g_T(v,\cdot)\big)
$.
\item[(ii)] (Optimal values.) One has
\begin{equation}\label{first.order.value}
V_M=\min_{u\in\hat S}\Jreg_M(u)+o_p(M^{-1/2}),
\qquad
\sqrt M\,(V_M-V)\xrightarrow[]{\mathfrak D}\min_{u\in\hat S}Y(u).
\end{equation}
If, moreover, $\hat S=\{\hat u\}$ is a singleton, then
$
\sqrt M\,(V_M-V)\xrightarrow[]{\mathfrak D}\mathcal N\big(0,\sigma^2(\hat u)\big).
$
\end{itemize}
\end{theorem}
\begin{proof}
We follow the proof of \cite[Theorem 5.7]{shapiro2021}, where compactness of
the feasible set is used twice. Here it is replaced by the reduction to
$\mathcal K_\Lambda$ below and by Lemma  \ref{entropy.K}. We first reduce the optimization problem to the compact set
$\mathcal K_\Lambda$. Let $B_M$ be the event of
\eqref{BM.def}, so that $\mathbb P(B_M)\le\|L_g\|^2_{L^2(\mu,\Omega)}/M$ by
\eqref{PBM}. By Corollary  \ref{regularity.stationary}, applied with $\mu$ and
with $\mu_M$, we have $\hat S\subset\mathcal K_\Lambda$ and, outside $B_M$,
$\argmin_{\Uad}\Jreg_M\subset\mathcal K_\Lambda$. Let
$\mathcal V(\psi):=\min_{\mathcal K_\Lambda}\psi$ for
$\psi\in C(\mathcal K_\Lambda)$. Then $V=\mathcal V(\Jreg)$,
$\argmin_{\mathcal K_\Lambda}\Jreg=\hat S$, and $V_M=\mathcal V(\Jreg_M)$ outside
$B_M$. Since $V_M$ is measurable, we obtain $\mathbb P\bigl(V_M\neq\mathcal V(\Jreg_M)\bigr)\le\mathbb P(B_M)\to0,$ so that
\begin{equation}\label{reduction.KLambda}
V_M=\mathcal V(\Jreg_M)+o_p(M^{-1/2}).
\end{equation}

We next establish the functional central limit theorem on
$C(\mathcal K_\Lambda)$. Since the Tikhonov term
cancels, $\sqrt M(\Jreg_M-\Jreg)=\sqrt M(J_M-J)=M^{-1/2}\sum_{i=1}^M\upchi_i$ on
$\mathcal K_\Lambda$, where $\upchi_i:=g_T(\cdot,X_i)-J$ are i.i.d.\ centered
random elements of the separable space $C(\mathcal K_\Lambda)$. By
Lemma  \ref{gTlipschitz},
$\sup_{u\in\mathcal K_\Lambda}\mathbb E[\upchi_1(u)^2]\le\|a\|^2_{L^2(\mu,\Omega)}$
and $|\upchi_1(u)-\upchi_1(v)|\le\Upsilon\|u-v\|_{L^2}$ with
$\Upsilon:=L_T(X_1)+\|L_T\|_{L^1(\mu,\Omega)}\in L^2$, while
Lemma  \ref{entropy.K} gives
$\int_0^1\sqrt{\log N(\varepsilon,\mathcal K_\Lambda,\|\cdot\|_{L^2})}\,d\varepsilon<\infty$. Thus 
Theorem  \ref{CLT.JainMarcus}  yields (i).

It remains to identify the asymptotic behavior of the optimal values. The function $\mathcal V$ is $1$-Lipschitz on
$C(\mathcal K_\Lambda)$ and, by Danskin's theorem
\cite[Proposition 4.12 and the subsequent remark]{BonnansShapiro2000},
Hadamard directionally differentiable at $\Jreg$ with
$\mathcal V'_{\Jreg}(\eta)=\min_{u\in\hat S}\eta(u)$. Since $\Jreg=V$ on
$\hat S$, the Delta theorem \cite[Theorem 9.74]{shapiro2021} and (i) give
$\mathcal V(\Jreg_M)=\min_{u\in\hat S}\Jreg_M(u)+o_p(M^{-1/2})$ and
$\sqrt M(\mathcal V(\Jreg_M)-V)\xrightarrow[]{\mathfrak D}\min_{u\in\hat S}Y(u)$.
By \eqref{reduction.KLambda} and \cite[Theorem 3.1]{bill99}, both statements
hold with $V_M$ in place of $\mathcal V(\Jreg_M)$, which proves (ii); if
$\hat S=\{\hat u\}$, the limit is $Y(\hat u)\sim\mathcal N(0,\sigma^2(\hat u))$.
\end{proof}

\subsection{Asymptotic distribution of empirical minimizers}\label{subsec:minimizers}

Throughout this subsection, $\hat u_M$ is a measurable selection of
$\argmin_{\Uad}\Jreg_M$ (Lemma  \ref{measurable.selection}), hence a random element
of $\Vspace$. Its limit is Gaussian when the constraints are inactive and, in
general, is the solution of a quadratic problem over the critical cone with a
Gaussian linear term, as in finite-dimensional stochastic programming
\cite{Shapiro1993,KingRockafellar1993}. We treat both cases in a single framework
in the infinite-dimensional space $L^2$, where the limiting problem is governed
by the Hessian of the regularized objective.

The following result provides both a strong law of large numbers and a central
limit theorem for i.i.d. random elements $\upchi_1,\upchi_2,\dots$ in
$\mathbb{H}$, with common mean $\mathfrak m$ and sample means
$\mathfrak m_M:=\frac{1}{M}\sum_{i=1}^{M}\upchi_i$.

\begin{theorem}[{\cite[Theorem 8.1.1]{Hsing2013}}]\label{CLT.Hilbert}
If $\mathbb{E}\left[\|\upchi_1\|_{\mathbb H}\right]<\infty$, then
$\|\mathfrak{m}_M-\mathfrak m\|_{\mathbb H}\rightarrow 0$ almost surely. If, in addition,
$\mathbb{E}\left[\|\upchi_1\|_{\mathbb H}^{2}\right]<\infty$, then
$
\sqrt{M}\,(\mathfrak m_M-\mathfrak m)
\xrightarrow[]{\mathfrak{D}}
\mathcal N(0,\mathscr K)
$ in $\mathbb H,
$ where $
\mathscr K
:=
\mathbb E\big[(\upchi_1-\mathfrak m)\otimes(\upchi_1-\mathfrak m)\big].
$
\end{theorem}

We apply Theorem  \ref{CLT.Hilbert} in $\mathbb H=\Vspace$, to the empirical
gradients, and in $\mathbb H=\HS(\Vspace)$, to the empirical Hessians. Recall
from Proposition  \ref{hessian.HS} that $\mathcal H(u)=\mathcal H_J(u)+\alpha\,\mathrm{Id}$
is the Riesz representative of $D^2\Jreg(u)$; accordingly we set
\begin{equation}\label{HM.def}
\mathcal H_M(u)
:=
\frac1M\sum_{i=1}^M\mathcal H_g(u,X_i)
+
\alpha\,\mathrm{Id},
\end{equation}
which is the Riesz representative of $D^2\Jreg_M(u)$. The regularization term
cancels in the differences:
\begin{equation}\label{cancellation}
\nabla\Jreg_M(u)-\nabla\Jreg(u)=\nabla J_M(u)-\nabla J(u),
\qquad
\mathcal H_M(u)-\mathcal H(u)
=
\frac1M\sum_{i=1}^M\mathcal H_g(u,X_i)-\mathcal H_J(u).
\end{equation}

\begin{lemma}[SLLN and CLT for the empirical gradients and Hessians]
\label{CLT.gradient}
Suppose that Assumptions  \ref{A1}-- \ref{A2} hold. Then, for every $u\in\mathcal O$,
$
\|\nabla \Jreg_M(u)-\nabla \Jreg(u)\|_{L^2}\rightarrow 0
\ \text{ a.s.},
\
\|\mathcal H_M(u)-\mathcal H(u)\|_{\mathcal L(L^2)}\rightarrow 0
\ \text{a.s.},$ and
$$
\sqrt M
\bigl(
\nabla \Jreg_M(u)-\nabla \Jreg(u)
\bigr)
\overset{\mathfrak D}{\longrightarrow}
\mathcal N(0,\mathscr K(u))
\qquad\text{in }\Vspace,
$$
where
$
\mathscr K(u)
:=
\mathbb E_\mu
\bigl[
(\nabla g_T(u,\cdot)-\nabla J(u))
\otimes
(\nabla g_T(u,\cdot)-\nabla J(u))
\bigr]
$
is the covariance operator of  \ $\nabla g_T(u,\cdot)$. Moreover, setting
$\bar C_M:=\frac1M\sum_{i=1}^MC_{\mathcal H}(X_i)$, which converges almost
surely to $\mathbb E_\mu[C_{\mathcal H}]<\infty$, one has, for all
$u,v\in\mathcal O$,
\begin{equation}\label{HM.lip}
\|\mathcal H_M(v)-\mathcal H(u)\|_{\mathcal L(L^2)}
\le\bar C_M\,\|v-u\|_{L^2}+\|\mathcal H_M(u)-\mathcal H(u)\|_{\mathcal L(L^2)}.
\end{equation}
In particular, if $u_M\in\mathcal O$ are random elements with $u_M\to u$ in
$\Vspace$ almost surely, then $\mathcal H_M(u_M)\to\mathcal H(u)$ in operator norm
almost surely.
\end{lemma}

\begin{proof}
Fix $u\in\mathcal O$. By  Proposition  \ref{1order},
$\nabla J_M(u)=\frac1M\sum_{i=1}^M\nabla g_T(u,X_i)$, and
$\upchi:=\nabla g_T(u,\cdot)$ is an
$\Vspace$-valued random element with $\|\upchi\|_{L^2}\le L_\nabla\in L^2(\mu,\Omega)$
and mean $\nabla J(u)$. Theorem  \ref{CLT.Hilbert} in $\mathbb H=\Vspace$ gives
the almost sure convergence and the central limit theorem for
$\nabla J_M(u)-\nabla J(u)$, and hence, by \eqref{cancellation}, for
$\nabla\Jreg_M(u)-\nabla\Jreg(u)$. For the Hessians, by
Proposition  \ref{hessian.HS}, $\omega\mapsto\mathcal H_g(u,\omega)$ is a random
element of the separable Hilbert space $\HS(\Vspace)$ with
$\|\mathcal H_g(u,\omega)\|_{\HS}\le c_{\HS}L_g(\omega)$ and
$\mathbb E_\mu[\mathcal H_g(u,\cdot)]=\mathcal H_J(u)$. Theorem  \ref{CLT.Hilbert} in
$\mathbb H=\HS(\Vspace)$ and
$\|\cdot\|_{\mathcal L(L^2)}\le\|\cdot\|_{\HS}$ give the second convergence.
Finally, by Lemma  \ref{2orderw}(ii)--(iii) and \eqref{D2gLIPS}, the Riesz
isometry gives
$\|\mathcal H_M(v)-\mathcal H_M(u)\|_{\mathcal L(L^2)}\le\bar C_M\|v-u\|_{L^2}$,
and \eqref{HM.lip} follows from the triangle inequality. Since
$C_{\mathcal H}\in L^2(\mu,\Omega)\subset L^1(\mu,\Omega)$, the strong law of
large numbers yields $\bar C_M\to\mathbb E_\mu[C_{\mathcal H}]$ almost surely, which
gives the last assertion.
\end{proof}

Throughout the rest of this subsection, $\hat u\in\Uad$ denotes a stationary
point of \eqref{Pmu}, $d:=\nabla\Jreg(\hat u)$, so that $-d\in N_{\Uad}(\hat u)$,
and $\mathcal H:=\mathcal H(\hat u)$. Besides the critical cone $\mathcal C(\hat u)$ and
the set $\mathcal D(\hat u)$ of critical radial directions of \eqref{CT.cones}, we
consider the closed linear subspace generated by the critical cone,
$
\mathcal E(\hat u):=\overline{\mathcal C(\hat u)-\mathcal C(\hat u)}^{\,L^2},
$
and the following strong second-order sufficient condition:
\begin{equation}\label{SSOSC}
Q(\delta u)=D^2\Jreg(\hat u)[\delta u,\delta u]>0
\qquad\text{for all }\delta u\in\mathcal E(\hat u)\setminus\{0\}.
\end{equation}
Since $\mathcal C(\hat u)\subset\mathcal E(\hat u)$, \eqref{SSOSC} implies the
assumption of Theorem  \ref{sos.conditions}. Moreover, the proof of
Lemma  \ref{coercivity.cone} only uses that $\mathcal C(\hat u)$ is a weakly
closed cone, so it applies to the closed subspace $\mathcal E(\hat u)$:
under \eqref{SSOSC} there exists $\beta>0$ such that
\begin{equation}\label{SSOSC.beta}
Q(\delta u)\ge\beta\|\delta u\|_{L^2}^2
\qquad\text{for all }\delta u\in\mathcal E(\hat u).
\end{equation}
For box constraints,
$\mathcal E(\hat u)=\{\delta u:\ \delta u_i=0\ \text{a.e. where }d_i\neq0\}$,
so that \eqref{SSOSC} only involves directions supported where the constraints
are inactive or weakly active.
We also make precise the nondegeneracy condition on the control constraints
under which the limit is Gaussian. Recall from
Corollary  \ref{regularity.stationary} that a stationary control is, after
modification on a null set, continuous on $[0,T]$.

\begin{definition}[Uniformly inactive constraints]\label{def.inactive}
We say that the control constraints are {\it uniformly inactive} at
$\hat u\in\Uad$ if there exists $\varepsilon_0>0$ such that
$\overline{B(\hat u(t),\varepsilon_0)}\subset U_c$ for a.e. $t\in[0,T]$, or,
equivalently, $\hat u\in\intinf(\Uad)$, the interior of $\Uad$ relative to the
$L^\infty$-topology.
\end{definition}

The only property of the constraint set that the asymptotic analysis requires is
the density of the critical radial directions in the critical cone, that is,
polyhedricity of $\Uad$ at $\hat u$. It holds in the two situations of interest.

\begin{lemma}[Density of the critical radial directions]\label{polyhedricity}
Suppose that Assumptions  \ref{A1}-- \ref{A2} hold and let $\hat u\in\Uad$ be a
stationary point. Then
\begin{equation}\label{poly.density}
\mathcal C(\hat u)=\overline{\mathcal D(\hat u)}^{\,L^2}
\end{equation}
in each of the following cases:
\begin{itemize}
\item[(i)] $U_c$ is a polytope;
\item[(ii)] the constraints are uniformly inactive at $\hat u$, for an
arbitrary compact convex set $U_c$. In this case, moreover, $d=0$ and
$\mathcal C(\hat u)=\mathcal E(\hat u)=\Vspace$.
\end{itemize}
\end{lemma}

\begin{proof}
(i) is the polyhedricity of $\Uad$ at $\hat u$, established in the proof of
Theorem  \ref{son.conditions} through \cite[Proposition 6.35]{BonnansShapiro2000}.
For (ii), let $\varepsilon_0$ be as in Definition  \ref{def.inactive}. Then
$\hat u\pm w\in\Uad$ for every $w\in L^\infty([0,T];\mathbb R^m)$ with
$\|w\|_{L^\infty}\le\varepsilon_0$, and \eqref{VI} gives
$\langle d,w\rangle_{L^2}=0$ for all such $w$. These $w$ span
$L^\infty([0,T];\mathbb R^m)$, which is dense in $\Vspace$, so $d=0$. The same
$w$ belong to $\mathcal R(\hat u)$, which is a cone, so
$L^\infty([0,T];\mathbb R^m)\subset\mathcal R(\hat u)$. Since $d=0$, we get
$\mathcal D(\hat u)=\mathcal R(\hat u)\supset L^\infty([0,T];\mathbb R^m)$ and
$\mathcal C(\hat u)=\overline{\mathcal R(\hat u)}^{\,L^2}=\Vspace$. Hence
\eqref{poly.density} holds, and $\mathcal E(\hat u)=\Vspace$.
\end{proof}
\begin{remark}\label{remark.heuristic}
We explain heuristically why the limit of the rescaled minimizers is
governed by problem \eqref{limit.problem} below. Set
$
t_M:=M^{-1/2},$ $ e_M:=\hat u_M-\hat u,$ $ h_M:=\sqrt M\,e_M,$ $
G_M:=\sqrt M\bigl(\nabla J_M(\hat u)-\nabla J(\hat u)\bigr),
$
so that $\hat u_M=\hat u+t_Mh_M$ and, by \eqref{cancellation},
$\nabla\Jreg_M(\hat u)=d+t_MG_M$. By Lemma  \ref{CLT.gradient},
$G_M\overset{\mathfrak D}{\longrightarrow}\mathcal G\sim\mathcal N(0,\mathscr K(\hat u))$
and $\mathcal H_M(\hat u)\to\mathcal H$ almost surely. A second-order Taylor
expansion of $\Jreg_M$ around $\hat u$ then gives, for fixed $h\in\Vspace$,
$$
M\bigl(\Jreg_M(\hat u+t_Mh)-\Jreg_M(\hat u)\bigr)
=\sqrt M\,\langle d,h\rangle_{L^2}+\langle G_M,h\rangle_{L^2}
+\tfrac12\langle\mathcal H_M(\hat u)h,h\rangle_{L^2}+o(1).
$$
where $o(1)\to0$ as $M\to\infty$, uniformly for $h$ in bounded sets. Since $\hat u_M=\hat u+t_Mh_M$ minimizes $\Jreg_M$ over $\Uad$, $h_M$ minimizes the left-hand
side over $(\Uad-\hat u)/t_M$, a set that increases to $\mathcal R(\hat u)$. By
\eqref{VI}, $\langle d,h\rangle_{L^2}\ge0$ for $h\in\mathfrak T(\hat u)$, and
directions with $\langle d,h\rangle_{L^2}>0$ are penalized by the diverging
term $\sqrt M\langle d,h\rangle_{L^2}$. The relevant directions are therefore
those of the critical cone
$\mathcal C(\hat u)=\mathfrak T(\hat u)\cap\{h\in\Vspace:\ \langle d,h\rangle_{L^2}=0\},$
on which the right-hand side converges to
$\frac12Q(h)+\langle\mathcal G,h\rangle_{L^2}$. One thus expects
$h_M\overset{\mathfrak D}{\longrightarrow}\Theta(\mathcal G)$, the minimizer of this
functional over $\mathcal C(\hat u)$; when the constraints are uniformly
inactive, $\mathcal C(\hat u)=\Vspace$ and $\Theta(\mathcal G)=-\mathcal H^{-1}\mathcal G$.
In infinite dimension, however, this passage to the limit cannot rely on
compactness. The proof below therefore works with the first-order condition
satisfied by $h_M$, namely the perturbed variational inequality \eqref{VI.t} below,
which is tested with critical radial directions using the density
\eqref{poly.density} and passed to the limit through the Legendre property of
$Q$ (Lemma  \ref{key.lemma}).
\end{remark}
\begin{lemma}[The limiting quadratic problem]\label{stampacchia}
Suppose that Assumptions  \ref{A1}-- \ref{A2} and \eqref{SSOSC} hold. Then, for every
$G\in\Vspace$, the problem
\begin{equation}\label{limit.problem}
\min_{h\in\mathcal C(\hat u)}\bigg\{\tfrac12\,Q(h)+\langle G,h\rangle_{L^2}\bigg\}
\end{equation}
admits a unique solution $\bar \Theta(G)\in\mathcal C(\hat u)$, which is also the unique $ h \in \mathcal C(\hat u)$
satisfying
\begin{equation}\label{limit.VI}
\langle\mathcal H  h+G,\,w- h\rangle_{L^2}\ge0
\ \ \text{for all }w\in\mathcal C(\hat u).
\end{equation}
Moreover,
$\|\Theta(G_1)-\Theta(G_2)\|_{L^2}\le\beta^{-1}\|G_1-G_2\|_{L^2}$ for all
$G_1,G_2\in\Vspace$, and $\Theta(\lambda G)=\lambda\Theta(G)$ for $\lambda\ge0$. Here $\beta>0$ is the constant in \eqref{SSOSC.beta}.
\end{lemma}
\begin{proof}
The closed subspace $\mathcal E(\hat u)$ is a Hilbert space on which, by
\eqref{SSOSC.beta}, the symmetric bilinear form
$(h,w)\mapsto\langle\mathcal Hh,w\rangle_{L^2}$ is coercive, and
$\mathcal C(\hat u)\subset\mathcal E(\hat u)$ is a nonempty closed convex set.
Since $w-h\in\mathcal E(\hat u)$ for $h,w\in\mathcal C(\hat u)$, we may replace
$G$ in \eqref{limit.problem} and \eqref{limit.VI} by its orthogonal projection
onto $\mathcal E(\hat u)$. Existence and uniqueness of $\Theta(G)$, and its
characterization by \eqref{limit.VI}, then follow from Stampacchia's theorem
\cite[Theorem 5.6]{brezis2010functional}
. 
If $h_i=\Theta(G_i)$, $i=1,2$, testing \eqref{limit.VI}
for $h_1$ with $w=h_2$ and for $h_2$ with $w=h_1$ and adding gives
$Q(h_1-h_2)\le\langle G_2-G_1,h_1-h_2\rangle_{L^2}$, and \eqref{SSOSC.beta}
yields the Lipschitz estimate. Positive homogeneity follows from uniqueness,
since $\mathcal C(\hat u)$ is a cone.
\end{proof}

To make Remark  \ref{remark.heuristic} rigorous, we pass to the limit in the
variational inequality satisfied by $h_M=\sqrt M(\hat u_M-\hat u)$ (see
\eqref{VI.t} and \eqref{VI.perturbed.proof}) through the following
deterministic stability result, in which polyhedricity enters through the
density of critical radial directions, and the Legendre property of $Q$
(Lemma  \ref{weak.cont.lemma}) replaces compactness.

\begin{lemma}\label{key.lemma}
Suppose that Assumptions  \ref{A1}-- \ref{A2} hold, that \eqref{SSOSC} holds at the
stationary point $\hat u$, and that \eqref{poly.density} holds
. Then, for every compact
set $\mathbb K\subset\Vspace$, every $\bar R>0$ and every $\varepsilon>0$, there exist
$t_0>0$ and $\eta>0$ with the following property: if $t\in(0,t_0)$,
$G\in\mathbb K$, $A\in\mathcal L(L^2)$ with $\|A-\mathcal H\|_{\mathcal L(L^2)}\le\eta$,
and $h\in\Vspace$ with $\|h\|_{L^2}\le\bar R$ satisfy $\hat u+th\in\Uad$ and
\begin{equation}\label{VI.t}
\bigl\langle d+t\,(G+Ah),\,v-\hat u-th\bigr\rangle_{L^2}\ge0
\quad\text{for all }v\in\Uad,
\end{equation}
then $\|h-\Theta(G)\|_{L^2}\le\varepsilon$.
\end{lemma}
\begin{proof}
    See Appendix  \ref{App.key.lemma}
\end{proof}
\begin{theorem}[Asymptotic distribution of empirical minimizers]
\label{asymp.active}
Suppose that Assumptions  \ref{A1}-- \ref{A2} hold and let
$\hat u\in\Uad$ be the unique global minimizer of \eqref{Pmu}. Assume that the
strong second-order sufficient condition \eqref{SSOSC} holds at $\hat u$, and
that either $U_c$ is a polytope or the control constraints are uniformly
inactive at $\hat u$. Then
\begin{equation}\label{conv.emp.u}
\sqrt M(\hat u_M-\hat u)
\overset{\mathfrak D}{\longrightarrow}
\Theta(\mathcal G)
\ \ \text{in }\Vspace,
\qquad
\mathcal G\sim\mathcal N\bigl(0,\mathscr K(\hat u)\bigr),
\end{equation}
where $\Theta$ is the solution of the quadratic problem \eqref{limit.problem}
and $\mathscr K(\hat u)$ is as in Lemma  \ref{CLT.gradient}.
\end{theorem}
\begin{proof} We use throughout the notation introduced in Remark  \ref{remark.heuristic}. In
particular,
$
t_M=M^{-1/2},
$ $ e_M=\hat u_M-\hat u,$ $
h_M=\sqrt M\,e_M,
$
and
$
G_M
=
\sqrt M\bigl(
\nabla\Jreg_M(\hat u)-\nabla\Jreg(\hat u)
\bigr)
=
\sqrt M\bigl(
\nabla J_M(\hat u)-\nabla J(\hat u)
\bigr).
$
By Lemma  \ref{CLT.gradient},
$
G_M\overset{\mathfrak D}{\longrightarrow}
\mathcal G\sim\mathcal N\bigl(0,\mathscr K(\hat u)\bigr)$ in $ \Vspace. $
Since $\hat u_M$ minimizes $\Jreg_M$ over the convex set $\Uad$, it satisfies
$$
\left\langle
\nabla\Jreg_M(\hat u_M),v-\hat u_M
\right\rangle_{L^2}
\geq0,
\ \text{ for all }
v\in\Uad.
$$
The segment joining $\hat u$ and $\hat u_M$ is contained in $\Uad\subset\mathcal O$. Hence, by the fundamental theorem of calculus,
$$
\nabla\Jreg_M(\hat u_M)
=
\nabla\Jreg_M(\hat u)
+
\overline{\mathcal H}_M e_M,
\quad \text{ 
where } \quad
\overline{\mathcal H}_M
:=
\int_0^1
\mathcal H_M\bigl(\hat u+\tau e_M\bigr)\,d\tau.
$$
We therefore have
$
\nabla\Jreg_M(\hat u_M)
=
d+t_M
\bigl(
G_M+\overline{\mathcal H}_M h_M
\bigr).
$
Consequently,
\begin{equation}\label{VI.perturbed.proof}
\left\langle
d+t_M
\bigl(
G_M+\overline{\mathcal H}_M h_M
\bigr),
v-\hat u-t_Mh_M
\right\rangle_{L^2}
\geq0,
\ \text{ for all } 
v\in\Uad.
\end{equation}
 By the strong second-order sufficient condition and Theorem  \ref{sos.conditions}, there exist $\varepsilon_1,\delta>0$ such that
$
\Jreg(\hat u+\mathfrak{e})-\Jreg(\hat u)
\geq
\frac{\delta}{2}\|\mathfrak{e}\|_{L^2}^2
$
whenever $\hat u+\mathfrak{e}\in\Uad$ and $\|\mathfrak{e}\|_{L^2}<\varepsilon_1$. Since
$
\Jreg_M(\hat u_M)\leq\Jreg_M(\hat u)
$
and $\Jreg_M-\Jreg=J_M-J$,
using the fundamental theorem of calculus,
$$
\Jreg(\hat u_M)-\Jreg(\hat u)
\leq(J-J_M)(\hat u_M)-(J-J_M)(\hat u)
=
-\int_0^1
\left\langle
\nabla(J_M-J)(\hat u+\tau e_M),
e_M
\right\rangle_{L^2}\,d\tau.
$$
Moreover,
$$
\nabla(J_M-J)(\hat u+\tau e_M)
=
t_MG_M
+
\int_0^\tau
\bigl(
\mathcal H_M-\mathcal H
\bigr)
(\hat u+s\,e_M)e_M\,ds.
$$
Thus, whenever $\|e_M\|_{L^2}<\varepsilon_1$,
$
\frac{\delta}{2}\|e_M\|_{L^2}^2
\leq
t_M\|G_M\|_{L^2}\|e_M\|_{L^2}
+
\varepsilon_M\|e_M\|_{L^2}^2,
$
where
$$
\varepsilon_M
:=
\sup_{s\in[0,1]}
\left\|
\mathcal H_M(\hat u+s\,e_M)
-
\mathcal H(\hat u+s\,e_M)
\right\|_{\mathcal L(L^2)}.
$$
By Lemma  \ref{2orderw}(ii) and \eqref{D2gLIPS},
$$
\varepsilon_M
\leq
\left\|
\mathcal H_M(\hat u)-\mathcal H(\hat u)
\right\|_{\mathcal L(L^2)}
+
\left(
\bar C_M+\mathbb E_\mu[C_{\mathcal H}]
\right)
\|e_M\|_{L^2},
$$
where $\bar C_M=M^{-1}\sum_{i=1}^M C_{\mathcal H}(X_i)$. By Lemma  \ref{CLT.gradient}, the first term in the r.h.s. of the inequality above converges to zero a.s., while the strong law of large numbers gives $\bar C_M\to\mathbb E_\mu[C_{\mathcal H}]$ a.s. Together with the consistency $\hat u_M\to\hat u$, this yields $
\varepsilon_M\longrightarrow0$ a.s.
Consequently, a.s. for all sufficiently large $M$,
$$
\|h_M\|_{L^2}
\leq
4\delta^{-1}\|G_M\|_{L^2}.
$$
 Since $(G_M)$ converges in distribution in the Polish space $\Vspace$, its laws
are uniformly tight \cite[Theorem 5.2]{bill99}: given $\zeta>0$, there is a
compact set $\mathbb K\subset\Vspace$ such that
$\sup_M\mathbb P(G_M\notin\mathbb K)\le\zeta$. Fix also $\varepsilon>0$.
Since $\mathbb K$ is bounded, there is a 
$\bar R>0$ such that
$\limsup_{M\to\infty}\mathbb P(\|h_M\|_{L^2}>\bar R)\le\zeta$. Let $t_0$ and $\eta$ be given by
Lemma  \ref{key.lemma}, whose assumption \eqref{poly.density} holds by
Lemma  \ref{polyhedricity}, and recall that $\overline{\mathcal H}_M\to\mathcal H$ a.s. in $\mathcal L(L^2)$. Therefore, for $M>t_0^{-2}$, on the event
$
\left\{
G_M\in\mathbb K,\quad
\|h_M\|_{L^2}\le\bar R,\quad
\|\overline{\mathcal H}_M-\mathcal H\|_{\mathcal L(L^2)}
\le\eta
\right\},
$
the variational inequality \eqref{VI.perturbed.proof} satisfies all the hypotheses of
Lemma  \ref{key.lemma}, so that $\|h_M-\Theta(G_M)\|_{L^2}\le\varepsilon$. Hence
$$
\limsup_{M\to\infty}
\mathbb P\left(
\|h_M-\Theta(G_M)\|_{L^2}>\varepsilon
\right)
\le 2\zeta,
$$
and, since $\zeta>0$ is arbitrary,
$
h_M-\Theta(G_M)
\longrightarrow 0$ in probability.

Finally, by Lemma  \ref{stampacchia},  $\Theta$ is continuous on $\Vspace$. Therefore, by Lemma  \ref{CLT.gradient} and the continuous mapping theorem for random elements in metric spaces
\cite[Theorem 18.1]{vanderVaart1998},
$
\Theta(G_M)
\overset{
}{\longrightarrow}
\Theta(\mathcal G)
$ in distribution. 
Combining this convergence with
$
h_M-\Theta(G_M)\rightarrow0
$ in probability, and recalling that $h_M=\sqrt M(\hat u_M-\hat u)$, we obtain \eqref{conv.emp.u}. 
\end{proof}
The limit $\Theta(\mathcal G)$ is Gaussian whenever the critical cone is a
linear subspace, in which case \eqref{SSOSC} reduces to the second-order
sufficient condition of Theorem  \ref{sos.conditions}.

\begin{corollary}[Gaussian limits]\label{asymp.gaussian}
Under the assumptions of Theorem  \ref{asymp.active}, suppose that
$\mathcal C(\hat u)$ is a closed linear subspace of $\Vspace$, and let $P$ be the
orthogonal projection onto it. Then $\mathcal H_{\mathcal C}:=P\mathcal H|_{\mathcal C(\hat u)}$
is an isomorphism of $\mathcal C(\hat u)$, the operator
$\mathscr B:=\mathcal H_{\mathcal C}^{-1}P\in\mathcal L(L^2)$ is self-adjoint, and
$$
\sqrt M(\hat u_M-\hat u)
\overset{\mathfrak D}{\longrightarrow}
\mathcal N\bigl(0,\Sigma(\hat u)\bigr)
\ \text{ in }\Vspace,
\qquad
\Sigma(\hat u):=\mathscr B\,\mathscr K(\hat u)\,\mathscr B .
$$
This applies in particular in the following two cases.
\begin{itemize}
\item[(a)] \emph{Uniformly inactive constraints.} If the constraints are uniformly
inactive at $\hat u$, then $\mathcal C(\hat u)=\mathcal E(\hat u)=\Vspace$,
$\mathscr B=\mathcal H^{-1}$ and
$
\Sigma(\hat u)=\mathcal H(\hat u)^{-1}\mathscr K(\hat u)\mathcal H(\hat u)^{-1}.
$

Here $U_c$ is an arbitrary compact convex set, and \eqref{SSOSC} reads
$Q(\delta u)>0$ for all $\delta u\in\Vspace\setminus\{0\}$. By
Remark  \ref{remark.alpha.large}, this holds in particular whenever
$\alpha>\|\mathcal H_J(\hat u)\|_{\mathcal L(L^2)}$.
\item[(b)] \emph{Box constraints with strict complementarity.} If
$U_c=\prod_{i=1}^m[a_i,b_i]$ and strict complementarity holds,
that is, $d_i(t)\neq0$ for a.e. $t$ in the active set
$\mathcal A_i:=\{t:\hat u_i(t)\in\{a_i,b_i\}\}$, $i=1,\dots,m$, then
\begin{equation}\label{strict.cone}
\mathcal C(\hat u)=\bigl\{\delta u\in\Vspace:\ \delta u_i=0\ \text{a.e. on }\mathcal A_i,\ i=1,\dots,m\bigr\}.
\end{equation}
In particular, the components of the limit vanish on the active sets.
\end{itemize}
\end{corollary}
\begin{proof}
Suppose first that $\mathcal C(\hat u)$ is a closed subspace of $\Vspace$, and
let $P$ be the orthogonal projection onto it. Then
$\mathcal E(\hat u)=\mathcal C(\hat u)$, and the variational inequality
\eqref{limit.VI} is equivalent to
$
h\in\mathcal C(\hat u),$ $ P(\mathcal Hh+G)=0,
$
that is, to $\mathcal H_{\mathcal C}h=-PG$, where
$\mathcal H_{\mathcal C}:=P\mathcal H|_{\mathcal C(\hat u)}$. Since $Ph=h$ for
$h\in\mathcal C(\hat u)$, \eqref{SSOSC.beta} gives
$\langle\mathcal H_{\mathcal C}h,h\rangle_{L^2}=\langle\mathcal Hh,h\rangle_{L^2}
\ge\beta\|h\|_{L^2}^2$ for all $h\in\mathcal C(\hat u)$. Hence, by the
Lax--Milgram theorem \cite[Corollary 5.8]{brezis2010functional},
$\mathcal H_{\mathcal C}$ is an isomorphism of $\mathcal C(\hat u)$ onto itself
with $\|\mathcal H_{\mathcal C}^{-1}\|_{\mathcal L(\mathcal C(\hat u))}\le\beta^{-1}$,
and $\Theta(G)=-\mathcal H_{\mathcal C}^{-1}PG$. Since $\mathcal H$ and $P$ are
self-adjoint, so is $\mathcal H_{\mathcal C}$, and the operator
$\mathscr B:=\mathcal H_{\mathcal C}^{-1}P\in\mathcal L(\Vspace)$, which vanishes on
$\mathcal C(\hat u)^\perp$, is self-adjoint and satisfies
$$
\Theta(G)=-\mathscr BG,
\qquad
\|\mathscr B\|_{\mathcal L(L^2)}\le\|\mathcal H_{\mathcal C}^{-1}\|\,\|P\|\le\beta^{-1}.
$$
Therefore, since $\mathcal G$ is centered Gaussian in $\Vspace$,
$\Theta(\mathcal G)=-\mathscr B\mathcal G$ is also centered Gaussian, with
covariance operator $\mathscr B\,\mathscr K(\hat u)\,\mathscr B$.

It remains to verify the two cases in the statement. In case (a),
Lemma  \ref{polyhedricity}(ii) gives $\mathcal C(\hat u)=\mathcal E(\hat u)=\Vspace$,
so that $P=\mathrm{Id}$ and $\mathscr B=\mathcal H^{-1}$; moreover, the density
property \eqref{poly.density} holds by the same lemma, so all the assumptions of
Theorem  \ref{asymp.active} are satisfied. In case (b), by Remark  \ref{remark.crit_cone_box_constraints},
$\delta u\in\mathcal C(\hat u)$ if and only if $\delta u$ satisfies the sign
conditions defining $\mathfrak T(\hat u)$ on the active sets and
$\delta u_i=0$ a.e.\ on $\{d_i\neq0\}$. Since $d_i=0$ on the inactive set and,
by strict complementarity, $d_i\neq0$ a.e.\ on $\mathcal A_i$, the latter
condition reads $\delta u_i=0$ a.e.\ on $\mathcal A_i$. This makes the sign
conditions redundant and we obatain the characterization \eqref{strict.cone}.
\end{proof}

\begin{remark}[Dependence on $\alpha$]\label{remark.trace}
In the setting of Corollary  \ref{asymp.gaussian}, the covariance operator
$\mathscr K(\hat u)$ is of trace class, being the covariance of a random element
of $\Vspace$ with finite second moment; together with
$\|\mathscr B\|_{\mathcal L(L^2)}\le\beta^{-1}$ this gives
$
\operatorname{tr}\mathscr K(\hat u)
\le
\|L_\nabla\|^2_{L^2(\mu,\Omega)}$ and $\operatorname{tr}\Sigma(\hat u)
\le
\|L_\nabla\|^2_{L^2(\mu,\Omega)}/\beta^{2},
$
with $L_\nabla$ as in Proposition  \ref{1order}(ii) and $\beta$ the coercivity
constant of \eqref{SSOSC.beta}. When
$\alpha>\|\mathcal H_J(\hat u)\|_{\mathcal L(L^2)}$ one may take
$\beta=\alpha-\|\mathcal H_J(\hat u)\|_{\mathcal L(L^2)}$, so that the
asymptotic variance decreases as $\alpha$ increases: a larger $\alpha$
stabilizes the SAA estimator, at the price of a larger bias with respect to the
unregularized problem.
\end{remark}
\begin{corollary}[Second-order asymptotics of the optimal values]
\label{cor.second.order.values}
Under the assumptions of Theorem  \ref{asymp.active}, let $V_M=\Jreg_M(\hat u_M)$
be the SAA optimal value  and $V=\Jreg(\hat u)$ that of \eqref{Pmu}, let
$G_M$ be as in Remark  \ref{remark.heuristic}, and for $G\in\Vspace,$ denote by
$ \mathfrak M(G)
$ 
the optimal value of \eqref{limit.problem}. Then
\begin{equation}\label{second.order.value}
V_M=\Jreg_M(\hat u)+\frac1M\,\mathfrak M(G_M)+o_p\bigl(M^{-1}\bigr),
\end{equation}
and
\begin{equation}\label{second.order.value.limit}
M\bigl(V_M-\Jreg_M(\hat u)\bigr)\overset{\mathfrak D}{\longrightarrow}
\mathfrak M(\mathcal G)=-\tfrac12\,Q\bigl(\Theta(\mathcal G)\bigr),
\qquad
\mathcal G\sim\mathcal N\bigl(0,\mathscr K(\hat u)\bigr).
\end{equation}
If, in addition, $\mathcal C(\hat u)$ is a closed linear subspace, then
$\mathfrak M(\mathcal G)=-\frac12\langle\mathscr B\mathcal G,\mathcal G\rangle_{L^2}$, with
$\mathscr B$ as in Corollary  \ref{asymp.gaussian}; in particular, if the
constraints are uniformly inactive at $\hat u$, then
$\mathfrak M(\mathcal G)=-\frac12\langle\mathcal H^{-1}\mathcal G,\mathcal G\rangle_{L^2}$.
\end{corollary}
\begin{proof}
We use the notation of Remark  \ref{remark.heuristic} and of the proof of
Theorem  \ref{asymp.active}, and write $F_G(h):=\frac12Q(h)+\langle G,h\rangle_{L^2}$.
Testing \eqref{limit.VI} with $w=0$ and $w=2\Theta(G)$ gives
\begin{equation}\label{Theta.orth}
\langle\mathcal H\Theta(G)+G,\Theta(G)\rangle_{L^2}=0,
\qquad\text{hence}\qquad
\mathfrak M(G)=F_G(\Theta(G))=-\tfrac12Q(\Theta(G)),
\end{equation}
and $\mathfrak M$ is continuous by Lemma  \ref{stampacchia}. A second-order Taylor expansion
of $\Jreg_M$ along $[\hat u,\hat u_M]\subset\Uad$, with
$\nabla\Jreg_M(\hat u)=d+t_MG_M$, gives
$$
M\bigl(V_M-\Jreg_M(\hat u)\bigr)
=\sqrt M\langle d,h_M\rangle_{L^2}+F_{G_M}(h_M)+o_p(1)
=\sqrt M\langle d,h_M\rangle_{L^2}+\mathfrak M(G_M)+o_p(1).
$$
Indeed, the remainder is $\frac12\langle(\widetilde{\mathcal H}_M-\mathcal H)h_M,h_M\rangle_{L^2}$,
where $\widetilde{\mathcal H}_M$ is an average of $\mathcal H_M$ over
$[\hat u,\hat u_M]$ and tends to $\mathcal H$ almost surely by \eqref{HM.lip};
moreover, $(h_M)$ and $(G_M)$ are bounded in probability,
$h_M=\Theta(G_M)+o_p(1)$ by the proof of Theorem  \ref{asymp.active}, and
$(G,h)\mapsto F_G(h)$ is Lipschitz on bounded sets. Finally, testing \eqref{VI}
with $u=\hat u_M$ and \eqref{VI.perturbed.proof} with $v=\hat u$, as in the proof
of Lemma  \ref{key.lemma}, and using \eqref{Theta.orth},
$$
0\le\sqrt M\langle d,h_M\rangle_{L^2}
\le-\bigl\langle G_M+\overline{\mathcal H}_Mh_M,h_M\bigr\rangle_{L^2}
=-\bigl\langle G_M+\mathcal H\Theta(G_M),\Theta(G_M)\bigr\rangle_{L^2}+o_p(1)
=o_p(1).
$$
This proves \eqref{second.order.value}, and \eqref{second.order.value.limit}
follows from Lemma  \ref{CLT.gradient} and the continuous mapping theorem
\cite[Theorem 18.1]{vanderVaart1998}, \cite[Theorem 3.1]{bill99}. If
$\mathcal C(\hat u)$ is a closed subspace, then $\Theta(G)=-\mathscr BG$ by the
proof of Corollary  \ref{asymp.gaussian}, and \eqref{Theta.orth} gives
$Q(\Theta(G))=\langle G,\mathscr BG\rangle_{L^2}$.
\end{proof}
Since $\mathbb E[\Jreg_M(\hat u)]=V$ and $\mathfrak M\le0$,
Corollary  \ref{cor.second.order.values} shows that the $O(M^{-1})$ term of
$V_M-V$ is nonpositive, in agreement with the downward bias of SAA optimal
values; under uniform integrability of $\{M(V_M-\Jreg_M(\hat u))\}_M$, one
obtains
$\mathbb E[V_M]-V=-\frac1{2M}\,\mathbb E\bigl[Q(\Theta(\mathcal G))\bigr]+o(M^{-1})$.
\subsection{Asymptotic confidence intervals}\label{subsec:inference}

The covariance operators in Corollary  \ref{asymp.gaussian} can be estimated
consistently from the sample, which yields asymptotic confidence intervals for
linear functionals of the optimal control.
\begin{corollary}[Plug-in estimation of the covariance]\label{inference}
Suppose that the assumptions of Corollary  \ref{asymp.gaussian} hold, let
$\xi_i:=\nabla g_T(\hat u_M,X_i)$ and
$$
\mathscr K_M:=\frac1M\sum_{i=1}^M\bigl(\xi_i-\nabla J_M(\hat u_M)\bigr)\otimes\bigl(\xi_i-\nabla J_M(\hat u_M)\bigr),
$$
and set $\Sigma_M:=\mathscr B_M\mathscr K_M\mathscr B_M$, where
\begin{itemize}
\item[(a)] if the constraints are uniformly inactive,
$\mathscr B_M:=\mathcal H_M(\hat u_M)^{-1}$;
\item[(b)] for box constraints with strict complementarity,
$\mathscr B_M:=(P_M\mathcal H_M(\hat u_M)P_M)^{-1}P_M$, with the inverse taken on
the range of $P_M$, where $P_M$ is the orthogonal projection onto
$\{\delta u:\ \delta u_i=0\ \text{a.e. where }\hat u_{M,i}\in\{a_i,b_i\},\ i=1,\dots,m\}$.
\end{itemize}
Then, almost surely, $\Sigma_M$ is defined for $M$ large and
$\|\Sigma_M-\Sigma(\hat u)\|_{\mathcal L(L^2)}\to0$. Consequently, for every
$h\in\Vspace$ with $\sigma_h^2:=\langle\Sigma(\hat u)h,h\rangle_{L^2}>0$,
$$
\frac{\sqrt M\,\langle\hat u_M-\hat u,h\rangle_{L^2}}{\langle\Sigma_Mh,h\rangle_{L^2}^{1/2}}
\xrightarrow{\mathfrak D}\mathcal N(0,1),
$$
so that, for $\gamma\in(0,1)$ and $z_{1-\gamma/2}$ the $(1-\gamma/2)$-quantile of the
standard normal distribution,
$\langle\hat u_M,h\rangle_{L^2}\pm z_{1-\gamma/2}\bigl(\langle\Sigma_Mh,h\rangle_{L^2}/M\bigr)^{1/2}$
is an asymptotic confidence interval of level $1-\gamma$ for
$\langle\hat u,h\rangle_{L^2}$.
\end{corollary}

\begin{proof}
(a) By Proposition  \ref{uniform.consistency.Linfty} and \eqref{HM.lip},
$\mathcal H_M(\hat u_M)\to\mathcal H(\hat u)$ in operator norm almost surely; hence
$\mathcal H_M(\hat u_M)$ is invertible for $M$ large and
$\mathcal H_M(\hat u_M)^{-1}\to\mathcal H(\hat u)^{-1}$
\cite[Theorem 10.7]{rudin1973functional}. Next,
$\mathscr K_M=\frac1M\sum_{i=1}^M\xi_i\otimes\xi_i-\nabla J_M(\hat u_M)\otimes\nabla J_M(\hat u_M)$.
Since $\|a\otimes a-b\otimes b\|_{\HS}\le(\|a\|_{L^2}+\|b\|_{L^2})\|a-b\|_{L^2}$,
Proposition  \ref{1order} gives
$\|\xi_i\otimes\xi_i-\nabla g_T(\hat u,X_i)\otimes\nabla g_T(\hat u,X_i)\|_{\HS}
\le2L_\nabla(X_i)\ell(X_i)\|\hat u_M-\hat u\|_{L^2}$, with
$L_\nabla\ell\in L^1(\mu,\Omega)$; the strong laws of large numbers in $\HS(\Vspace)$
and in $\mathbb R$ then give $\mathscr K_M\to\mathscr K(\hat u)$ in $\HS(\Vspace)$,
and similarly $\nabla J_M(\hat u_M)\to\nabla J(\hat u)$ in $\Vspace$, almost
surely. Hence $\Sigma_M\to\Sigma(\hat u)$ in $\mathcal L(L^2)$.

(b) Let $\hat{\mathcal A}_{M,i}:=\{t:\hat u_{M,i}(t)\in\{a_i,b_i\}\}$. By
Lemma  \ref{gradient.regularity}, the convergence
$\nabla J_M(\hat u_M)\to\nabla J(\hat u)$ holds with uniformly bounded Lipschitz
constants in time, hence uniformly on $[0,T]$. Together with
\eqref{projection.formula}, Proposition  \ref{uniform.consistency.Linfty} and strict
complementarity, this shows that
$\|\mathbf 1_{\hat{\mathcal A}_{M,i}}-\mathbf 1_{\mathcal A_i}\|_{L^1(0,T)}\to0$, so
that $P_M\to P$ strongly, $P$ being the orthogonal projection onto
$\mathcal C(\hat u)$. The Legendre property (Lemma  \ref{weak.cont.lemma}) and
\eqref{SSOSC.beta} then give, as in Lemmas  \ref{coercivity.cone}
and  \ref{key.lemma}, a uniform coercivity of $\mathcal H_M(\hat u_M)$ on
$\operatorname{ran}P_M$ and the strong convergence $\mathscr B_M\to\mathscr B$.
Since $\mathscr K_M\to\mathscr K(\hat u)$ in $\HS(\Vspace)$, as in (a), and
$\mathscr K(\hat u)$ is compact, $\Sigma_M\to\Sigma(\hat u)$ in $\mathcal L(L^2)$.
Since these steps only combine the arguments of (a) with those of
Lemmas  \ref{coercivity.cone} and  \ref{key.lemma}, now applied to the converging
subspaces $\operatorname{ran}P_M$, we omit the details.

In both cases, $\langle\Sigma_Mh,h\rangle_{L^2}\to\sigma_h^2$ almost surely, and the
asymptotic normality of the studentized statistic follows from
Corollary  \ref{asymp.gaussian} and Slutsky's lemma \cite[Theorem 3.1]{bill99}; the
confidence interval is obtained by inverting it.
\end{proof}

\section{Numerical illustration}\label{sec:numerics}

We illustrate Sections  \ref{sec:quantitative} and  \ref{sec:asymptotics} on two
examples: a linear oscillator with uncertain frequency, which is convex and for
which all assumptions can be checked analytically, and a nonconvex model of the
sterile insect technique with uncertain biological rates, for which uniqueness
and the second-order conditions are verified numerically.

\subsection{Methodology}
Controls are piecewise constant on a uniform grid of $[0,T]$ with $N$ intervals,
and the analysis applies to the resulting finite-dimensional problem; in the second example the state equation is
integrated with a semi-implicit Euler scheme (four substeps per control interval), and the gradients used are the exact
gradients of the discrete problem. The ``true'' problem is discretized in $\omega$ by tensor Gauss--Legendre quadrature
(200 nodes in the first example and $20\times20=400$ in the second, where for the truncated normal law of $\delta_S$ the
nodes are taken in the quantile variable), so that $\hat u$, $V$,
$\sigma^2(\hat u)=\operatorname{Var}_\mu[g_T(\hat u,\cdot)]$ and
$\Sigma(\hat u)=\mathscr B\mathscr K(\hat u)\mathscr B$
(Theorem  \ref{asymp.values}, Corollary  \ref{asymp.gaussian}) are computed to quadrature accuracy;
SAA solutions are compared with $\hat u$ on the same grid and with the same scheme, so that
the fluctuations are not contaminated by discretization errors.
For each $M$ we generate $N_{\rm rep}$ independent samples and report
$\mathbb E|V_M-V|$, $\mathbb E\|\hat u_M-\hat u\|_{L^2}^2$, the variances of
$Z_h:=\sqrt M\langle\hat u_M-\hat u,h\rangle_{L^2}$ for $h_1\equiv1$ and
$h_2(t):=\cos(\pi t/T)$, the pointwise standard deviation of
$\sqrt M(\hat u_M-\hat u)(t)$, and the empirical coverage of the 95\% intervals
of Corollary  \ref{inference}. In both examples the minimizer $\hat u$ is
unique (proved for the first example; for the second, which is not convex, it is observed numerically from ten different
initial controls), so that the limit in Theorem  \ref{asymp.values} is Gaussian and the
asymptotic predictions are
$$
\mathbb E|V_M-V|\approx\sigma(\hat u)\sqrt{2/\pi}\,M^{-1/2},\qquad
\mathbb E\|\hat u_M-\hat u\|_{L^2}^2\approx\operatorname{tr}\Sigma(\hat u)\,M^{-1},\qquad
\operatorname{Var}Z_h\approx\langle\Sigma(\hat u)h,h\rangle_{L^2}.
$$
In both examples the parameters range over compact sets and the data satisfy the local version of
Assumptions  \ref{A1}-- \ref{A2} of Remark  \ref{local.assumptions}; this is
immediate for the linear oscillator and is verified in
Subsection  \ref{subsec:sit} for the second example.
\subsection{An oscillator with uncertain frequency}
\begin{figure}[ht]
\centering
\includegraphics[width=\textwidth]{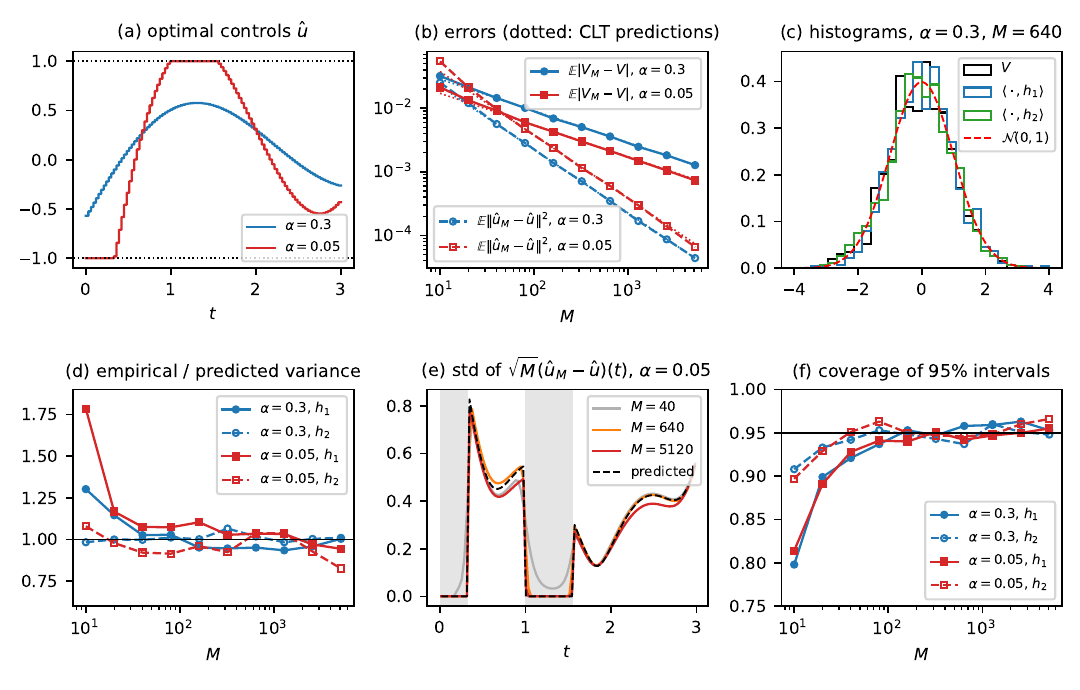}
\caption{Oscillator with uncertain frequency. (a) Optimal controls. (b) Errors
versus $M$, with asymptotic predictions (dotted). (c) Histograms of
$\sqrt M(V_M-V)/\sigma(\hat u)$ and $Z_{h_i}/\langle\Sigma h_i,h_i\rangle^{1/2}$,
$\alpha=0.3$, $M=640$. (d) Empirical over predicted variances of $Z_{h_i}$.
(e) Pointwise standard deviation of $\sqrt M(\hat u_M-\hat u)$, $\alpha=0.05$,
active arcs shaded. (f) Empirical coverage of the 95\% confidence intervals.}
\label{fig:oscillator}
\end{figure}

We consider
$$
\dot x_1=x_2,\qquad\dot x_2=-\omega^2x_1+u,\qquad x(0)=(1,0),\qquad
g(x)=\tfrac12|x|^2,\qquad |u(t)|\le1,
$$
with $T=3$, $\omega\sim\mathcal U[0.5,1.5]$ and $N=100$: a single control must
steer an ensemble of oscillators with different frequencies to the origin. Since
$x^u(T,\omega)$ is affine in $u$ and $g$ is convex, $\mathcal H_J\ge0$; hence the
problem is strictly convex, its minimizer is unique and \eqref{SSOSC} holds for
every $\alpha>0$. The SAA problems are convex quadratic programs with box
constraints, solved exactly by bounded-variable least squares. We take
$\alpha=0.3$, for which the constraints are inactive ($\max|\hat u|\approx0.58$),
and $\alpha=0.05$, for which $\hat u$ has two active arcs covering $30$ of the
$100$ intervals, with strict complementarity (Figure  \ref{fig:oscillator}(a));
$N_{\rm rep}=1000$ and $M\in\{10,20,\dots,5120\}$.

The results agree closely with the theory (Figure  \ref{fig:oscillator}). For
$\alpha=0.3$, $\sqrt M\,\mathbb E|V_M-V|$ stays in $[0.087,0.10]$ (prediction
$0.091$) and $M\,\mathbb E\|\hat u_M-\hat u\|^2_{L^2}$ in $[0.21,0.25]$
($\operatorname{tr}\Sigma(\hat u)=0.223$); the squared error thus decays like
$M^{-1}$, faster than the bound $M^{-1/2}$ of Remark  \ref{remark.rates.controls},
which does not use the central limit theorem. The histograms are close to the
standard normal density already for moderate $M$, and the variance ratios lie in
$[0.93,1.07]$ for $M\ge40$, within the Monte Carlo standard error of about
$4.5\%$. For $\alpha=0.05$, the fluctuations vanish on the active arcs, as
predicted by Corollary  \ref{asymp.gaussian}(b), and the pointwise standard
deviation converges to the predicted profile, which peaks at the junctions.
Convergence is slower there: grid points that are inactive at $\hat u$ but close
to the bound are clipped for finite $M$, which reduces the variance until the
fluctuations become smaller than the distance to the bound; this explains the
ratio $0.83$ for $h_2$ at $M=5120$. The empirical coverage lies in $[0.94,0.97]$
for $M\ge80$, in agreement with Corollary  \ref{inference}(a) for $\alpha=0.3$ and
Corollary  \ref{inference}(b) for $\alpha=0.05$. Finally, $\operatorname{tr}\Sigma(\hat u)$ decreases
monotonically with $\alpha$, from $1.13$ for $\alpha=0.01$ to $0.095$ for
$\alpha=1$ and $0.0033$ for $\alpha=10$, which quantifies the stabilizing effect
of the regularization (Remark  \ref{remark.trace}).

\subsection{Sterile insect technique with uncertain biological rates}\label{subsec:sit}
The sterile insect technique (SIT) reduces a wild population by releasing sterilized males.
We use the model of \cite{AMS.SIT}, a modification of those in \cite{AronnaDumont,BidiCoronHayat},
normalized by the carrying capacity $K$, with aquatic stage $A$, females $F$, wild males $M$,
sterile males $S$, and control $u\in[0,1]$ (release rate):
\begin{equation}\label{sit.model}
\begin{aligned}
\dot A&=b\,\frac{MF}{M+qS}\,(1-A)-(\delta_A+\lambda)A, &\qquad
\dot F&=r\lambda A-\delta_F F,\\
\dot M&=(1-r)\lambda A-\delta_M M, &\qquad
\dot S&=\mathfrak{n}\,u-\delta_S S,
\end{aligned}
\end{equation}
where $q$ is the mating competitiveness of sterile males and $\delta_S$ is their daily mortality rate, so that a released male lives on average $1/\delta_S$ days. Under a constant release $u$, the sterile population settles at $S^*=\mathfrak{n} u/\delta_S$ and its weight in the mating term is $qS^*=\mathfrak{n} u\,q/\delta_S$, so that $q/\delta_S$ measures the sterilizing efficiency per unit of release. These two parameters, which describe the quality of the released males, are the uncertain ones: $\omega=(q,\delta_S)$ with independent components $q\sim\mathcal U[0.91,1]$ and $\delta_S$ normal with mean $0.12$ (the value in \cite{BidiCoronHayat}, i.e., a mean lifespan $1/\delta_S\approx8$ days, shorter than the $1/\delta_M=10$ days of wild males) and standard deviation $0.01$, truncated to $[0.10,0.14]$, so that sterile males do not outlive wild ones; the last two values are our choice.

The remaining parameters are deterministic and equal to the values of Table 1 of \cite{BidiCoronHayat}:
$b=8$, $\lambda=0.25$, $\delta_A=0.03$, $\delta_F=0.04$, $\delta_M=0.10$, $r=0.49$, $K=5\cdot10^4$ and
$\mathfrak{n}=u_{\max}/K=6$ with $u_{\max}=3\cdot10^5$. Since the mating factor appears here in $\dot A$ and not in $\dot F$ as in
\cite{BidiCoronHayat}, these are only reference values. The horizon is $T=90$ days and the initial
state is the wild equilibrium without sterile males. The cost is $\mathbb E_\mu[F(T)+M(T)]+\frac\alpha2\|u\|_{L^2}^2$; we omit the
linear release cost, which kept the bound $u=0$ active near $T$ for all the values of $\alpha$ we tested. The problem is not convex. The field is singular at $M+qS=0$, which
admissible controls do not reach but which may occur in the neighborhood $\mathcal O$ of $\Uad$ where the assumptions are
required. Replacing $S$ in the denominator by a $C^2$ function equal to $S$ for $S\ge0$ and bounded below by
$-\varepsilon$ does not change the problem on $\Uad$, and for $\varepsilon$ small enough it keeps the denominator
away from zero on $\mathcal O$, so Assumptions  \ref{A1}-- \ref{A2} hold locally (Remark  \ref{local.assumptions}).

We take $N=90$ intervals of one day, $N_{\rm rep}=500$, $M\in\{10,20,\dots,2560\}$ and two regimes: $\alpha=0.1$
(inactive constraints) and $\alpha=0.03$ (upper bound active). The coverage of the 95\% intervals is estimated with $200$
replications at $M\in\{20,40,\dots,1280\}$.

\begin{figure}[ht]
\centering
\includegraphics[width=\textwidth]{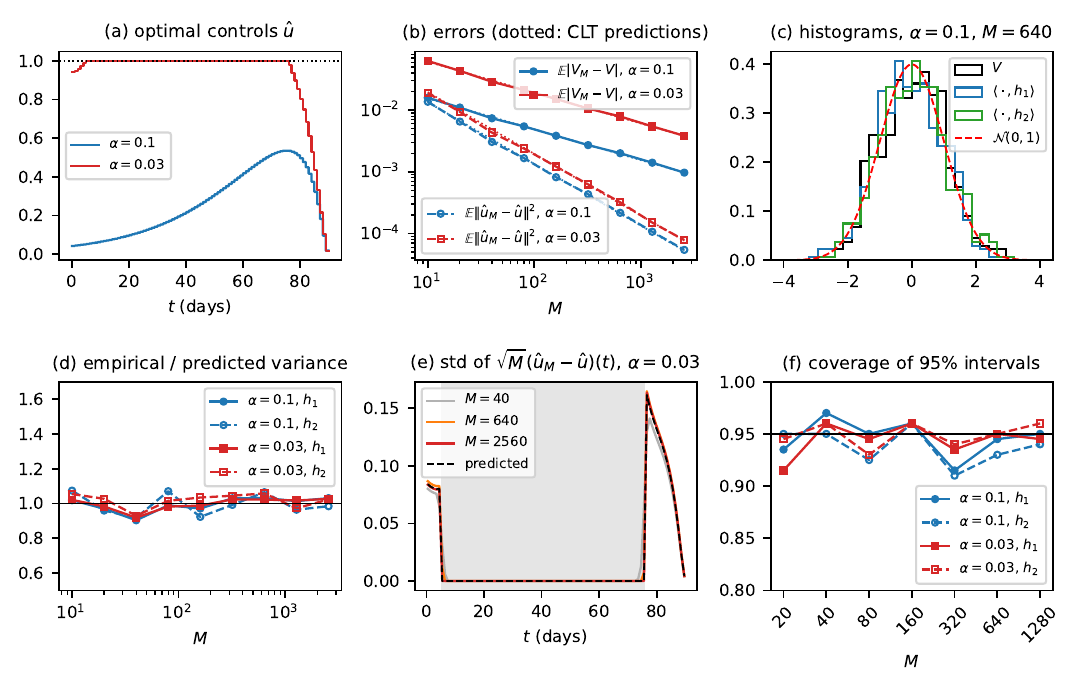}
\caption{Sterile insect technique. (a) Optimal controls. (b) Errors versus $M$ with the asymptotic predictions (dotted).
(c) Normalized histograms for $\alpha=0.1$, $M=640$. (d) Empirical/predicted variance of $Z_{h_i}$.
(e) Pointwise standard deviation of $\sqrt M(\hat u_M-\hat u)$ for $\alpha=0.03$; active arc shaded.
(f) Empirical coverage of the 95\% intervals.}
\label{fig:sit}
\end{figure}

For $\alpha=0.1$ the optimal control satisfies $0.019\le\hat u\le0.54$, so no bound is active (only in the discretized problem, since $\hat u(T)$ is expected to vanish as the grid is refined). For $\alpha=0.03$ the bound $\hat u=1$ is active on $71$ consecutive intervals, from day $5$ to day $76$. Ten optimizations from different initial controls end at the same point (distance $\le8\cdot10^{-8}$), so the minimizer is unique in practice. In both regimes the Hessian is positive definite (smallest eigenvalues $0.070$ and $0.020$), so \eqref{SSOSC} holds, and strict complementarity holds for $\alpha=0.03$, although the multiplier is small near the junctions ($\ge1.6\cdot10^{-4}$). The numerical results agree with the asymptotic theory (Figure~\ref{fig:sit}).
Panel (b): by Theorem~\ref{asymp.values}, $\sqrt M\,(V_M-V)$ is asymptotically
$\mathcal N(0,\sigma^2(\hat u))$, so $\sqrt M\,\mathbb E|V_M-V|\to\sigma(\hat u)\sqrt{2/\pi}$; by
Corollary~\ref{asymp.gaussian}, $M\,\mathbb E\|\hat u_M-\hat u\|_{L^2}^2\to\operatorname{tr}\Sigma(\hat u)$.
Numerically, both quantities stay within $5\%$ and $11\%$ of these limits ($0.049$ and $0.134$
for $\alpha=0.1$; $0.193$ and $0.196$ for $\alpha=0.03$).
Panel (c): each projection $Z_h$ is asymptotically $\mathcal N(0,\langle\Sigma h,h\rangle)$;
the normalized histograms of $V_M$ and of $Z_{h_1},Z_{h_2}$ are close to the standard normal density.
Panel (d): the variance ratios of $Z_{h_1},Z_{h_2}$ lie in $[0.90,1.07]$ for $M\ge40$.
Panel (e): for $\alpha=0.03$, the pointwise standard deviation of $\sqrt M(\hat u_M-\hat u)$
matches the predicted profile and vanishes on the active arc, as predicted by
Corollary~\ref{asymp.gaussian}.
Panel (f): the coverage of the $95\%$ intervals lies in $[0.91,0.97]$ (mean $0.944$),
consistent with the nominal level given the Monte Carlo error ($1.5\%$ with $200$ replications),
in agreement with Corollary~\ref{inference}(a) for $\alpha=0.1$ and
Corollary~\ref{inference}(b) for $\alpha=0.03$.
\section{Concluding remarks}\label{sec:conclusions} We have developed a first- and second-order analysis of the sample average approximation of optimal control problems with uncertain control-affine dynamics and Tikhonov regularization, which culminates in the asymptotic distribution of the empirical optimal controls in $L^2$. Two structural features made this possible. The affine dependence on the control removes the two-norm discrepancy and makes the Hessian of the uncertain cost a random Hilbert--Schmidt operator. The Tikhonov term turns the first-order condition into a projection formula, which provides compactness, and it is the only source of coercivity of the Hessian, a compact perturbation of $\alpha\,\mathrm{Id}$. All constants deteriorate as $\alpha\downarrow0$, while the bias with respect to the unregularized problem grows with $\alpha$.

Several questions remain open: dynamics that are not affine in the control (two-norm discrepancy); the case $\alpha=0$, with bang--bang and singular arcs; non-polyhedral constraint sets, where a curvature term is expected in the limiting problem of Theorem  \ref{asymp.active}.
\appendix






\section{Proofs of technical results}\label{appendixA}
 The main tool to prove the well-posedness is the following generalization of the Picard Lemma
\cite[Lemma 5.6.3]{polak2012}, allowing the Lipschitz coefficient to depend on
time through an integrable function $K\in L^1([0,T])$ rather than being a
constant. 
\begin{lemma}
\label{lem:picard}
Let $u\in\Vspace$ and consider the differential
equation
$\dot x(t)=f(x(t),u(t))$ for a.e. $t\in[0,T],$ and
$x(0)=\xi.$
Assume that there exists a function $K\in L^1([0,T])$ such that
$$
\|f(x,u(t))-f(y,u(t))\|
\le
K(t)\|x-y\|
$$
for all $x,y\in\mathbb{R}^n$ and for almost every
$t\in[0,T]$. Then, given any absolutely continuous function
$x_0:[0,T]\to\mathbb{R}^n$ such that $f(x_0(\cdot),u(\cdot))\in L^1([0,T];\mathbb R^n)$,
there exists a solution $x(\cdot)$ of the
differential equation such that, for all $t\in[0,T]$,
{\small $$
\|x(t)-x_0(t)\|
\le
e^{\int_0^t K(s)ds}
E(x_0,u),
\ \
E(x_0,u)
:=
\|\xi-x_0(0)\|
+
\int_0^T
\|\dot x_0(t)-f(x_0(t),u(t))\|\,dt.
$$}
Moreover, the solution is unique.
\end{lemma}
\begin{proof}
Let $x_0:[0,T]\to\mathbb{R}^n$ be absolutely continuous with
$f(x_0(\cdot),u(\cdot))\in L^1([0,T];\mathbb R^n)$. Define the
Picard iterates by
$
x_{i+1}(t)
=
\xi+\int_0^t f(x_i(s),u(s))\,ds,
$ $ i=0,1,2,\ldots
$
By induction and the Lipschitz condition, all $x_i(\cdot)$ are absolutely continuous. Hence,
for all $t\in[0,T]$,
$$
\|x_1(t)-x_0(t)\|
\le
\|\xi-x_0(0)\|
+
\int_0^t
\|\dot x_0(s)-f(x_0(s),u(s))\|\,ds
\le E(x_0,u).
$$
For $i\ge1$, using the Lipschitz condition,
$$
\|x_{i+1}(t)-x_i(t)\|
\le
\int_0^t
K(s)\|x_i(s)-x_{i-1}(s)\|\,ds.
$$
Proceeding by induction, we obtain
$$
\|x_{j+1}(t)-x_j(t)\|
\le
\frac{E(x_0,u)}{j!}
\left(\int_0^t K(s)\,ds\right)^j,
\qquad j\ge0,
$$
the induction step following from
$
\int_0^t
K(s)
\left(\int_0^s K(r)\,dr\right)^j ds
=
\frac{1}{j+1}
\left(\int_0^t K(s)\,ds\right)^{j+1}.
$
Consequently, for any $i,k\in\mathbb{N}$,
$$
\|x_{i+k}(t)-x_i(t)\|
\le
E(x_0,u)
\sum_{j=i}^{i+k-1}
\frac{1}{j!}
\left(\int_0^t K(s)\,ds\right)^j.
$$
Since $K\in L^1([0,T])$, the series
$
\sum_{j=0}^{\infty}
\frac{1}{j!}
\left(\int_0^t K(s)\,ds\right)^j
=e^{\int_0^tK(s)ds}
$
is uniformly bounded on $[0,T]$. Hence $\{x_i\}$ is a Cauchy
sequence in $C([0,T];\mathbb{R}^n)$ and converges uniformly to some
$x(\cdot)$. Passing to the limit in the Picard iteration gives
$
x(t)
=
\xi+\int_0^t f(x(s),u(s))\,ds,
$
so $x(\cdot)$ is a solution. Taking $i=0$ and letting $k\to\infty$ yields
$$
\|x(t)-x_0(t)\|
\le
e^{\int_0^tK(s)ds}
E(x_0,u).
$$
Uniqueness follows from the same Lipschitz estimate and
Gr\"onwall's inequality.
\end{proof}
\subsection{Proof of Proposition \ref{wellposedness}}\label{App.wellposedness}
\begin{proof}
(i) By \eqref{Lipschitz.bounds}, the map $x\mapsto f(x,u(t),\omega)$ is Lipschitz
continuous with constant $K_u(t):=L\big(1+\sqrt m\,\|u(t)\|\big)$, and
$K_u\in L^1([0,T])$ since $u\in\Vspace$. Hence the existence and uniqueness of an
absolutely continuous solution follow from Lemma~\ref{lem:picard}, applied with
the constant reference trajectory $x_0\equiv\varphi(\omega)$, for which
$f(x_0,u(\cdot),\omega)\in L^1([0,T];\mathbb R^n)$ by \eqref{growth.f}.

(ii) For all $t\in[0,T]$ and $\omega\in\Omega,$
$
x^u(t,\omega)
=
\varphi(\omega)
+
\int_0^t
f(x^u(s,\omega),u(s),\omega)\,ds,
$
so that, by \eqref{growth.f} and the Cauchy--Schwarz inequality,
$$
\begin{aligned}
\|x^u(t,\omega)\|
\le&
\|\varphi(\omega)\|
+
L\int_0^t
\big(1+\sqrt m\,\|u(s)\|\big)\,ds
+
L\int_0^t
\big(1+\sqrt m\,\|u(s)\|\big)\,\|x^u(s,\omega)\|
\,ds \\
\le& \|\varphi(\omega)\|
+
\Lambda_u
+
L\int_0^t
\big(1+\sqrt m\,\|u(s)\|\big)\,\|x^u(s,\omega)\|
\,ds.
\end{aligned}
$$
Applying Gr\"onwall's inequality yields
$
\|x^u(t,\omega)\|
\le
\big(\|\varphi(\omega)\|+\Lambda_u\big)\,e^{\Lambda_u}$ for all $ t\in[0,T].$ 

(iii) We define
$
\Lambda_{r_\rho}:=L\big(T+\sqrt{mT}\,r_\rho\big).$
For  $u\in\mathcal O$ we have $\|u\|_{L^2} <r_\rho$, hence
$
\Lambda_u<\Lambda_{r_\rho},
$
and since $\varphi(\omega)\in K_0$ for all $\omega \in \Omega,$ with $K_0$ compact, the constant
$
R:=\Big(\max_{\xi\in K_0}\|\xi\|+\Lambda_{r_\rho}\Big)e^{\Lambda_{r_\rho}}
$
satisfies $\|x^u(t,\omega)\|\le R$ for all $(t,\omega)\in[0,T]\times\Omega$ and every $u\in\mathcal O$. Therefore,
$
x^u(t,\omega)
\in
\mathcal{X}_0
:=
\overline{B(0,R)}
$
for all $ (t,\omega)\in[0,T]\times\Omega$ and $u\in\mathcal O.$
Moreover, from the bound \eqref{growth.f}, for a.e. $t\in[0,T]$ and every $\omega\in\Omega$,
$
\left\|\dot x^u(t,\omega)\right\|
=
\left\|f(x^u(t,\omega),u(t),\omega)\right\|
\le
L(1+R)\big(1+\sqrt m\,\|u(t)\|\big).
$
Since $u\in\Vspace$, it follows that $\dot x^u(\cdot,\omega)\in L^2([0,T];\mathbb R^n)$, and hence
$
x^u(\cdot,\omega)\in W^{1,2}([0,T];\mathbb R^n).
$
\end{proof}

\subsection{Proof of Proposition \ref{traj.lipschitz}}\label{App.traj.lipschitz}
\begin{proof}
Fix $u,v\in\mathcal O$ and $\omega\in\Omega$, and let $x^u(\cdot,\omega)$ and $x^v(\cdot,\omega)$ be the unique solutions of \eqref{dynamics} with controls $u$ and $v$, respectively. By \eqref{Lipschitz.bounds}, for a.e. $t\in[0,T]$ and $x,y\in \mathbb{R}^n$,
$$
\|
f(x,u(t),\omega)
-
f(y,u(t),\omega)
\|
\le
K_u(t)\|x-y\|,
$$
with
$
K_u(t)=L\bigl(1+\sqrt m\,\|u(t)\|\bigr)
$
and
$
\int_0^tK_u(s)\,ds
\le
\Lambda_{r_\rho}
$ for $ t\in[0,T].$
We apply Lemma~\ref{lem:picard} with the reference trajectory
$
x_0(t):=x^v(t,\omega).
$
Since both trajectories have the same initial condition, the residual appearing in
Lemma~\ref{lem:picard} is
$$
E(x_0,u)
=
\int_0^T
\bigl\|
f(x^v(t,\omega),v(t),\omega)
-
f(x^v(t,\omega),u(t),\omega)
\bigr\|\,dt.
$$
Using the control-Lipschitz estimate in \eqref{Lipschitz.bounds}, the
uniform bound
$
\|x^v(t,\omega)\|\le R$
from Proposition~\ref{wellposedness}, and the Cauchy--Schwarz inequality, we obtain
$$
E(x_0,u)
\le
L\sqrt m\,(1+R)
\int_0^T\|u(t)-v(t)\|\,dt \le L\sqrt{mT}\,(1+R)\,
\|u-v\|_{L^2}.
$$
Lemma~\ref{lem:picard} now yields
\begin{equation*}
\|x^u(t,\omega)-x^v(t,\omega)\|
\le
e^{\int_0^tK_u(s)ds}
E(x_0,u)\le L_{x}
\|u-v\|_{L^2},
\end{equation*}
where $L_x:=L\sqrt{mT}\,(1+R)e^{\Lambda_{r_\rho}},$ which concludes the proof.
\end{proof}

\subsection{Proof of Lemma \ref{p.lc.w.t.u}}\label{App.p.lc.w.t.u}
\begin{proof}
(i) The bound follows from  \eqref{f.EDO.bound}
and Remark~\ref{bounded.deriv.g}. The Lipschitz estimate follows the same
strategy as the proof of Theorem~5.6.9(b) in \cite{polak2012}, adapted to the
presence of the random parameter $\omega$ and using Lemma~\ref{lem:picard} to
estimate the differences between the corresponding trajectories and adjoint
solutions; the resulting constant is of the form $L_p=c\,L_g$ with $c$
depending only on $T,L,m,K_0$ and $r_\rho$.

(ii) If $u\in\Uad$, then $u(s)\in U_c$ for a.e. $s$ and
$x^u(s,\omega)\in \mathcal{X}_0$, so that $\|f_x(x^u(s,\omega),u(s),\omega)\|\le L_f$ by
Remark~\ref{bounded.deriv.f}. By \eqref{adjoint.t} and (i),
$\|\dot p^u(s,\omega)\|\le L_fe^{\Lambda_{r_\rho}}L_g(\omega)$ for a.e. $s\in[0,t]$,
which, together with $p^u(\cdot,\omega)\in W^{1,1}([0,t];\mathbb R^n)$, gives the claim.
\end{proof}

\subsection{Proof of Proposition \ref{Diff.x}}\label{App.variational}
\begin{proof}
(i) Fix $u\in\mathcal O$, $\omega\in\Omega$, and $\delta u\in\Vspace$.
By the variation-of-constants formula, the solution of
\eqref{variational} is given by
$
\delta x(t)
=
\int_0^t
\Phi^u(t,s,\omega)\,
f_u\bigl(x^u(s,\omega),\omega\bigr)\,
\delta u(s)\,ds.
$
Since
$f_u(x^u(\cdot,\omega),\omega)\delta u(\cdot)
\in L^1([0,T];\mathbb R^n)$, this representation yields
$
\delta x\in C([0,T];\mathbb R^n).
$
Linearity with respect to $\delta u$ follows directly from
the variational equation, whose coefficients are fixed for given
$u$ and $\omega$, and uniqueness follows from the uniqueness of the
solution of the corresponding homogeneous linear equation. By \eqref{growth.f},
$
\bigl\|f_u(x^u(t,\omega),\omega)\bigr\|
\leq L\sqrt m\,(1+R),
$
hence, from \eqref{f.EDO.bound}, for $t\in[0,T]$,
$$
\|\delta x(t)\|
\leq
e^{\Lambda_{r_\rho}}
L\sqrt m\,(1+R)
\int_0^t\|\delta u(s)\|\,ds
\leq
e^{\Lambda_{r_\rho}}
L\sqrt{mT}\,(1+R)
\|\delta u\|_{L^2}.
$$
Therefore,
\begin{equation}\label{dx<du}
\|\delta x\|_{C([0,T];\mathbb R^n)}
\leq
C_1\|\delta u\|_{L^2},
\ \text{where } \
C_1:=
e^{\Lambda_{r_\rho}}
L\sqrt{mT}\,(1+R).
\end{equation}

(ii) For all $u,v\in \mathcal{O},$ $\omega\in \Omega$ and $\delta u \in \Vspace,$ we have
$$\delta x^u(t)[\delta u]=\int_0^t\Big( f_x(x^u(s,\omega),u(s),\omega)\,\delta x^u(s)[\delta u] +f_u(x^u(s,\omega),\omega)\,\delta u(s)\Big)\,ds,$$
and the same identity with $v$ in place of $u$. Subtracting both identities and
using the Lipschitz continuity of the state mapping
$
u\mapsto x^u(\cdot,\omega)
$
(Proposition \ref{traj.lipschitz}), the estimates \eqref{Lipschitz.bounds} and
\eqref{dx<du}, Gr\"onwall's inequality yields a constant $L_{\delta x}>0$ depending only on $T,L,m,K_0$ and $r_\rho$ such that
$
\|\delta x^{u}(t)[\delta u]-\delta x^{v}(t)[\delta u]\|
\le
L_{\delta x}\|\delta u\|_{L^2}
\|u-v\|_{L^2}
$
for every $t\in[0,T]$.
\end{proof}

\subsection{Proof of Theorem \ref{Diff.x2}}\label{App.Diff.x2}
\begin{proof}
(i) Let $\omega\in\Omega$ fixed, $u\in \mathcal O,$ and   $\delta u_1,\delta u_2\in\Vspace$. For all $t\in[0,T],$ consider $\delta x_i(t)$ the solution of the first variational equation \eqref{variational} for $i=1,2,$ and
$\delta^2x_{12}(t)$  the solution of the second variational equation \eqref{variational2}.

First, we establish a priori bounds on $\delta^2x_{12}$.  From \eqref{variational2} and Remark \ref{A1.implications},
\begin{align*}
\|\delta^2x_{12}(t)\|
&\leq
\int_0^t
L\bigl(1+\sqrt m\,\|u(s)\|\bigr)
\|\delta^2x_{12}(s)\|\,ds
+
\int_0^t
L\bigl(1+\sqrt m\,\|u(s)\|\bigr)
\|\delta x_1(s)\|\,\|\delta x_2(s)\|\,ds
\\
&\quad+
L\sqrt m
\int_0^t
\|\delta x_1(s)\|\,\|\delta u_2(s)\|\,ds
+
L\sqrt m
\int_0^t
\|\delta u_1(s)\|\,\|\delta x_2(s)\|\,ds.
\end{align*}
Using \eqref{dx<du}, and recalling that
$
\int_0^t
L\bigl(1+\sqrt m\,\|u(s)\|\bigr)\,ds
\leq
\Lambda_{r_\rho},
$
we obtain
$$
\|\delta^2x_{12}(t)\|
\leq
\int_0^t
L\bigl(1+\sqrt m\,\|u(s)\|\bigr)
\|\delta^2x_{12}(s)\|\,ds
+
\bigl(\Lambda_{r_\rho}C_1^2
+
2L\sqrt{mT}\,C_1\bigr)
\|\delta u_1\|_{L^2}
\|\delta u_2\|_{L^2}.
$$
Therefore, by Gr\"onwall's inequality,
\begin{equation}
\label{d2x<du}
\|\delta^2x_{12}\|_{L^\infty}
\leq
C_2
\|\delta u_1\|_{L^2}
\|\delta u_2\|_{L^2},
\ \text{ 
where } \
C_2
:=
e^{\Lambda_{r_\rho}}
\left(
\Lambda_{r_\rho}C_1^2
+
2L\sqrt{mT}\,C_1
\right).
\end{equation}
By the variation of constants formula,
$$
\begin{aligned}
\delta^2x_{12}(t)
=
\int_0^t
\Phi^u(t,s,\omega)
\Big(
f_{xx}(s)[\delta x_1(s),\delta x_2(s)]
+
f_{xu}(s)[\delta x_1(s),\delta u_2(s)]
 +
f_{ux}(s)[\delta u_1(s),\delta x_2(s)]
\Big)\,ds.
\end{aligned}
$$
From the above, and since $\delta x_i$ is linear in $\delta u_i,$ $i=1,2,$ the mapping
$
(\delta u_1,\delta u_2)
\mapsto
\delta^2x_{12}(t)
$
is a well-defined bounded bilinear operator from
$\Vspace\times\Vspace$ into $\mathbb R^n$, which is symmetric since $f(\cdot,\cdot,\omega)\in C^2(\mathbb{R}^n\times\mathbb{R}^m).$

We now show that the control-to-state mapping is twice Fr\'echet differentiable in $\mathcal O$. Let $\delta u\in\Vspace$ be such that $v:=u+\delta u\in\mathcal O$, and denote $x^{v}(t):=x^v(t,\omega)$. Let $\delta x$ and $\delta^2x$ denote, respectively, the solutions of the first and second variational equations associated with $\delta u$. Define the remainder
$
r(t,\omega)
:=
x^{v}(t)-x^u(t)-\delta x(t)-\frac{1}{2}\delta^2x(t)
$
and set
$
h(t,\omega):=\delta x(t)+\frac12\delta^2x(t).
$ Applying Lemma~\ref{lem:picard} to the state equation for $x^v$, with reference trajectory
$
x_0(t,\omega):=
x^u(t,\omega)+h(t,\omega),
$
and using the Lipschitz coefficient
$
K_v(t):=L\bigl(1+\sqrt m\,\|v(t)\|\bigr),
$ with $
\int_0^tK_v(s)\,ds
\le
\Lambda_{r_\rho}$ for
$ t\in[0,T],$
we obtain
\begin{equation}\label{rest.bound0}
\|r(t,\omega)\|
\le
e^{\Lambda_{r_\rho}}
\int_0^T
\left\|
\dot{x}^u(s)
+\dot{\delta x}(s)
+\frac12\dot{\delta^2x}(s)
-
f\bigl(x^u(s)+h(s),v(s),\omega\bigr)
\right\|ds.
\end{equation}
Since $f$ is affine with respect to the control variable, we have, with $v-u=\delta u$,
$
f(x^u,v,\omega)
=
f(x^u,u,\omega)
+
f_u(x^u,\omega)\delta u
$
and
$
f_x(x^u,v,\omega)
=
f_x(x^u,u,\omega)
+
f_{xu}(x^u,\omega)\delta u.
$
We next use the integral form of Taylor's theorem with respect to the state variable \cite[Theorem 5.1.29]{polak2012}: for each $s\in[0,t]$,
$$
\begin{aligned}
f\bigl(x^u(s)+h(s),v(s),\omega\bigr)
&=
f\bigl(x^u(s),v(s),\omega\bigr)
+
f_x\bigl(x^u(s),v(s),\omega\bigr)h(s)
\\
&\quad+
\int_0^1(1-\theta)
f_{xx}\bigl(x^u(s)+\theta h(s),v(s),\omega\bigr)
[h(s),h(s)]\,d\theta.
\end{aligned}
$$
Combining this identity with the definitions of $\dot x^u$, $\dot{\delta x}$, and $\dot{\delta^2x}$, all zeroth- and first-order terms cancel, and we obtain, after suppressing the time variable $s$,
$$
\begin{aligned}
\dot{x}^u
+\dot{\delta x}
+\frac12\dot{\delta^2x}
-
f(x^u+h,v,\omega)
&=
\frac12
f_{xx}(x^u,u,\omega)[\delta x,\delta x]
-
\int_0^1(1-\theta)
f_{xx}(x^u+\theta h,v,\omega)[h,h]\,d\theta
\\
&\quad \hspace{6cm}
-\frac12
f_{xu}(x^u,\omega)[\delta^2x,\delta u].
\end{aligned}
$$
We next decompose the first two terms as
$$
\begin{aligned}
\frac12
f_{xx}(x^u,u,\omega)[\delta x,\delta x]
&-
\int_0^1(1-\theta)
f_{xx}(x^u+\theta h,v,\omega)[h,h]\,d\theta
\\
&=
\int_0^1(1-\theta)
\Bigl(
f_{xx}(x^u,u,\omega)
-
f_{xx}(x^u+\theta h,v,\omega)
\Bigr)
[\delta x,\delta x]\,d\theta
\\
&\quad
-
\int_0^1(1-\theta)
f_{xx}(x^u+\theta h,v,\omega)
[\delta x,\delta^2x]\,d\theta
\\
&\quad-\frac14
\int_0^1(1-\theta)
f_{xx}(x^u+\theta h,v,\omega)
[\delta^2x,\delta^2x]\,d\theta.
\end{aligned}
$$
From \eqref{Lipschitz.bounds}, the bounds \eqref{dx<du}--\eqref{d2x<du}, and since $\int_0^1(1-\theta)\theta\,d\theta=\tfrac16$ and $\int_0^1(1-\theta)\,d\theta=\tfrac12$, we have
$$
\begin{aligned}
 \bigg\|\int_0^1(1-\theta)
\Bigl(
f_{xx}&(x^u,u,\omega)
-
f_{xx}(x^u+\theta h,v,\omega)
\Bigr)
[\delta x,\delta x]\,d\theta\bigg\|\\&\le\left[\frac{L}{6}\big(1+\sqrt m\|u(s)\|\big)\Big(\|\delta x\|+\frac{1}{2}\|\delta^2 x\|\Big)+\frac{L}{2}\sqrt m\,\|\delta u(s)\|\right]\,\|\delta x\|^2\\
&\le \frac{L}{6}\big(1+\sqrt m\|u(s)\|\big)\Big(C_1^3\|\delta u\|_{L^2}^3+\frac{1}{2}C_1^2C_2\|\delta u\|_{L^2}^4\Big)\\
&\quad+\frac{L}{2}\sqrt m\, C_1^2\,\|\delta u(s)\|\|\delta u\|_{L^2 }^2,
\end{aligned}
$$
$$
\begin{aligned}
\bigg\|\int_0^1(1-\theta)
f_{xx}(x^u+\theta h,v,\omega)
[\delta x,\delta^2x]\,d\theta\bigg\|&\le \frac{L}{2}\big(1+\sqrt m\,\|v(s)\|\big)\,\|\delta x\|\,\|\delta^2x\|
\\&\le \frac{L}{2}\big(1+\sqrt m\,\|v(s)\|\big)\,C_1 C_2\|\delta u\|_{L^2}^3,
\end{aligned}
$$
and
$$
\begin{aligned}
\bigg\|\frac14
\int_0^1(1-\theta)
f_{xx}(x^u+\theta h,v,\omega)
[\delta^2x,\delta^2x]\,d\theta\bigg\| &\le
\frac{L}{8}\big(1+\sqrt m\,\|v(s)\|\big)\,\|\delta^2x\|^2
\\&\le \frac{L}{8}C_2^2\big(1+\sqrt m\,\|v(s)\|\big)\,\|\delta u\|^4_{L^2}.
\end{aligned}
$$
For the last term, from \eqref{growth.f},
$$
\left\|
\frac{1}{2}f_{xu}(x^u,\omega)
[\delta^2x,\delta u]
\right\|
\le
\frac{L}{2}\sqrt m\,
\|\delta^2x\|\|\delta u(s)\|\le\frac{L}{2}\sqrt m\,
C_2\|\delta u\|_{L^2}^2\|\delta u(s)\|.
$$
Integrating these estimates over $[0,T]$, using $$\int_0^T\|\delta u(s)\|\,ds\le\sqrt T\|\delta u\|_{L^2} \ \text{ and } \ \int_0^TL(1+\sqrt m\|w(s)\|)\,ds\le\Lambda_{r_\rho}$$ for $w\in\{u,v\}$, and combining with \eqref{rest.bound0}, we conclude that there exist constants $C_3,C_4>0$, independent of $\delta u$ and $\omega$, such that $\|r(\cdot,\omega)\|_{L^{\infty}}
\le
C_3\|\delta u\|_{L^2}^3
+
C_4\|\delta u\|_{L^2}^4.$ Hence
$
\|r(\cdot,\omega)\|_{L^{\infty}}=o(\|\delta u\|_{L^2}^2),
$
and the control-to-state mapping is twice Fr\'echet differentiable in $\mathcal O$, with
$$
D^2x^u(t,\omega)[\delta u,\delta u]
=
\delta^2x(t)
$$ and, by Proposition \ref{Diff.x},
$
Dx^u(t,\omega)[\delta u]
=
\delta x(t).
$ By the bilinearity and symmetry of the second variational equation,
$
D^2x^u(t,\omega)[\delta u_1,\delta u_2]
=
\delta^2x_{12}(t,\omega).
$

\noindent(ii) We have
$$
\begin{aligned}
D^2x^{u}(t,\omega)[\delta u_1,\delta u_2]
=
\int_0^t
\Big(&f_x(s)\,D^2x^{u}(s,\omega)[\delta u_1,\delta u_2]+
f_{xx}(s)[\delta x_1(s),\delta x_2(s)]
\\
&+
f_{xu}(s)[\delta x_1(s),\delta u_2(s)]
+
f_{ux}(s)[\delta u_1(s),\delta x_2(s)]
\Big)\,ds.
\end{aligned}
$$
From \eqref{Lipschitz.bounds}, \eqref{dx<du}, \eqref{d2x<du}, Proposition~\ref{traj.lipschitz} and Proposition~\ref{Diff.x}(ii), Gr\"onwall's inequality yields a constant $L_{\delta^2 x}>0$ depending only on $T,L,m,K_0$ and $r_\rho$ such that
$$
\|D^2x^{u}(t,\omega)[\delta u_1,\delta u_2]
-
D^2x^{v}(t,\omega)[\delta u_1,\delta u_2]\|
\le
L_{\delta^2 x}
\|\delta u_1\|_{L^2}
\,
\|\delta u_2\|_{L^2}\|u-v\|_{L^2}
$$
for every
$u,v\in\mathcal O$ and $\delta u_i\in\Vspace,$ $i=1,2,$
which completes the proof.
\end{proof}

\subsection{Proof of Corollary \ref{Diff.x2Linfty}}\label{App.Diff.x2Linfty}
\begin{proof}
Since $
\|\delta u\|_{L^2}
\le \sqrt{T}\,\|\delta u\|_{L^\infty},
$
the remainder term satisfies
$$
\|r(\cdot,\omega)\|_{L^\infty}
\le
C_3T^{3/2}\|\delta u\|_{L^\infty}^3
+
C_4T^2\|\delta u\|_{L^\infty}^4
=
o\big(\|\delta u\|_{L^\infty}^2\big).
$$
The continuity of the differentials follows from items~(ii) of Proposition~\ref{Diff.x} and Theorem~\ref{Diff.x2}, together with the continuous embedding $L^\infty([0,T];\mathbb R^m)\hookrightarrow \Vspace$.
\end{proof}

\subsection{Proof of Proposition \ref{2order}}\label{proof.2order}

\begin{proof} (i)
    For fixed $\omega \in \Omega,$ let $u\in\mathcal O$ and $\delta u\in\Vspace$ be such that $v:=u+\delta u\in\mathcal O,$ and denote by $x^v(t)$ and $x(t)$ the trajectories solving \eqref{dynamics} with controls $v$ and $u$, respectively. From  Proposition  \ref{Diff.x}  and Theorem~\ref{Diff.x2}, we have
\begin{equation}\label{Tay.x2}
x^v(s) = x(s) +  \delta x(s) + \frac{1}{2} \delta^2 x(s) + r(s,\omega),
\end{equation}
where the remainder $r(\cdot,\omega)$  satisfies $\|r(\cdot,\omega)\|_{L^\infty}\le C_3\|\delta u\|_{L^2}^3+C_4\|\delta u\|_{L^2}^4$ with $C_3,C_4$ independent of $\omega$, $\delta x$ satisfies the variational equation \eqref{variational}, and, since $f$ is affine with respect to the control variable, $\delta^2 x$ satisfies
\begin{equation*}
\begin{cases}
\dot{\delta^2 x}(s)
= f_x\delta^2 x(s)
+ f_{xx}[\delta x(s),\delta x(s)]
 + 2 f_{xu}[\delta x(s),\delta u(s)], &s\in[0,t],
\\
\delta^2 x(0)=0,
\end{cases}
\end{equation*}
where we have omitted the dependence of the derivatives of $f$ on $(x(s),u(s),\omega).$

Let
$
h(t,\omega)=x^v(t)-x(t).
$ Since $g(\cdot,\omega)$ is twice differentiable and $ g_t(u,\omega) = g(x(t), \omega)$,
the second-order Taylor formula \cite[Theorem 5.1.29]{polak2012} gives
\begin{equation}\label{gt.taylor0}
g_t(v,\omega) - g_t(u,\omega)
=
\left\langle
\nabla_x g(x(t),\omega),h(t,\omega)
\right\rangle
+\frac12
\left\langle
\nabla_{xx}g(x(t),\omega)h(t,\omega),
h(t,\omega)
\right\rangle
+
r_g(t,\omega),
\end{equation}
where
$$
r_g(t,\omega)
=
\int_0^1
(1-\theta)
\Big\langle
\Big(
\nabla_{xx}g(x(t)+\theta h(t,\omega),\omega)
-
\nabla_{xx}g(x(t),\omega)
\Big)
h(t,\omega),
h(t,\omega)
\Big\rangle
\,d\theta.
$$
Since $\nabla_{xx}g(\cdot,\omega)$ is Lipschitz continuous with constant $L_g(\omega)$,
$$
|r_g(t,\omega)|
\le
\int_0^1
(1-\theta)
L_g(\omega)\theta
\|h(t,\omega)\|^3
\,d\theta
=
\frac{L_g(\omega)}{6}
\|h(t,\omega)\|^3.
$$
Using
$
h(t,\omega)
=
\delta x(t)
+\frac12\delta^2x(t)
+r(t,\omega),
$
the inequalities \eqref{dx<du}, \eqref{d2x<du} and the bound on $r$, we get $
\|h(\cdot,\omega)\|_{L^\infty}
\le C_h\|\delta u\|_{L^2}
$ for $\|\delta u\|_{L^2}\le1$, with $C_h$ independent of $\omega$, and consequently
$
|r_g(t,\omega)|
\le
\frac{C_h^3}{6}\,
L_g(\omega)\,
\|\delta u\|_{L^2}^3,
$
uniformly in $t$.

On the other hand, from \eqref{gt.taylor0},
\begin{equation}\label{gt.taylor1}
\begin{aligned}
g_t(v,\omega) - g_t(u,\omega)
&= \left\langle
\nabla_x g(x(t),\omega),\delta x(t)
\right\rangle
\\&\quad+ \frac{1}{2} \Big(
\langle\nabla_{x} g(x(t),\omega), \delta^2 x(t) \rangle
+ \langle \nabla_{xx} g(x(t),\omega)\delta x(t), \delta x(t) \rangle
\Big)
+ R_g(t,\omega),
\end{aligned}
\end{equation}
with
$$
R_g(t,\omega) = \langle \nabla_x g(x(t),\omega), r(t,\omega) \rangle
+ \sum_{I,J \in \mathcal{I}} c_{IJ} \left\langle \nabla_{xx} g(x(t),\omega) \, I, J \right\rangle
+ r_g(t,\omega),
$$
where
 $\mathcal{I} = \{\delta x(t), \delta^2 x(t), r(t,\omega)\}$,  $c_{IJ}\ge0$ are numerical constants, and $c_{IJ}=0$ for $I=J=\delta x(t).$
Using \eqref{dx<du}, \eqref{d2x<du}, the bound on $r$ and Remark \ref{bounded.deriv.g}, we obtain
\begin{equation}\label{Rg.bound}
|R_g(t,\omega)|\le C_R\,L_g(\omega)\,\|\delta u\|_{L^2}^3
\qquad\text{for all }t\in[0,T],\ \omega\in\Omega,\ \|\delta u\|_{L^2}\le1,
\end{equation}
with $C_R$ independent of $t$, $\omega$ and $\delta u$.

Returning to  \eqref{gt.taylor1},
from Proposition \ref{1order} we have
$
D g_t(u,\omega)[\delta u]
= \langle \nabla_x g(x(t),\omega), \delta x(t) \rangle.
$ To compute the second-order term, we use the adjoint equation \eqref{adjoint.t}, which gives $p(t)=\nabla_xg(x(t),\omega)$, so that
$$
\begin{aligned}
g_t(v,\omega) &- g_t(u,\omega)
 =  D g_t(u,\omega)[\delta u]
+ \frac{1}{2} \Big(
\langle p(t), \delta^2 x(t) \rangle
+ \langle \nabla_{xx} g(x(t),\omega)\delta x(t), \delta x(t) \rangle
\Big)
+ R_g(t,\omega).
\end{aligned}
$$
It remains to express $\langle p(t), \delta^2 x(t) \rangle$. We have
$$
\begin{aligned}
&\frac{d}{ds} \langle p(s), \delta^2 x(s) \rangle
= \langle -f_x^\top p(s), \delta^2 x(s) \rangle
+ \langle p(s), f_x \delta^2 x(s) \rangle
+ \langle p(s), f_{xx}[\delta x,\delta x] \rangle
+ 2 \langle p(s), f_{xu}[\delta x,\delta u] \rangle,
\end{aligned}
$$
and the first two terms cancel. Integrating in time, since $\delta^2 x(0)=0$, and using the definition of the Hamiltonian $H$,
$$
\begin{aligned}
\langle p(t), \delta^2 x(t) \rangle
&= \int_0^t \Big[
\langle H_{xx}\delta x, \delta x \rangle
+ 2 \langle H_{xu}\delta x, \delta u \rangle
\Big] ds\\
&=
\int_0^t
\bigg\langle
\begin{pmatrix}
H_{xx}(s) & H_{xu}(s)\\
H_{ux}(s) & 0
\end{pmatrix}
\begin{pmatrix}
\delta x(s)\\
\delta u(s)
\end{pmatrix}, \begin{pmatrix}
\delta x(s)\\
\delta u(s)
\end{pmatrix}\bigg\rangle ds.
\end{aligned}
$$
Collecting all second-order terms, we have
\begin{equation}\label{gt.taylor}
g_t(v,\omega) - g_t(u,\omega)
= D g_t(u,\omega)[\delta u] + \frac{1}{2} D^2 g_t(u,\omega)[\delta u,\delta u]
+ R_g(t,\omega),
\end{equation}
where $D^2 g_t(u,\omega)[\delta u,\delta u]$ is given by \eqref{F2order} with $\delta v=\delta u$.
Now, by the variation of constants formula,
$$
\delta x(s)
=
\int_0^s
\Phi^u(s,\tau,\omega)f_u(x(\tau),\omega)\,\delta u(\tau)\,d\tau=:\mathcal S(\delta u)(s), \qquad s\in[0,t],
$$
and from \eqref{dx<du},
$
\|\mathcal S(\delta u)\|_{L^\infty}\leq C_1 \|\delta u\|_{L^2},
$
so that $\mathcal{S}$ is linear and bounded.  Therefore, each term appearing in
\eqref{F2order} is a continuous bilinear form on $\Vspace$: the first one is the composition of the bounded bilinear form
$\nabla_{xx}g(x(t),\omega)$ with $\mathcal S$, and the integral term is continuous since $\|H_{xx}(s)\|\le\|p(s)\|\,L(1+\sqrt m\|u(s)\|)$ is integrable and $\|H_{xu}(s)\|\le\|p(s)\|L\sqrt m$ is bounded (Lemma~\ref{p.lc.w.t.u}(i) and \eqref{growth.f}). By \eqref{gt.taylor} and \eqref{Rg.bound}, \eqref{F2order} is the second Fr\'echet differential of
$g_t(\cdot,\omega)$ at $u$.

\medskip

\noindent(ii)--(iii) Consider $\delta u_1,\delta u_2\in \Vspace.$ From \eqref{gt.taylor1} and polarization, for all $t\in[0,T]$ and  $u\in \mathcal O$,
\begin{equation*}
D^2 g_t(u,\omega)[\delta u_1,\delta u_2]=\left\langle
\nabla_x g(x^u(t),\omega),
\delta^2x_{12}^u(t)
\right\rangle
+\left\langle
\nabla_{xx}g(x^u(t),\omega)\delta x_1^u(t),
\delta x_2^u(t)
\right\rangle,
\end{equation*}
where $\delta x_i^u$, $i=1,2,$ and $\delta^2 x^u_{12} $ are the solutions of the variational equations \eqref{variational} and \eqref{variational2} with control $u$. Thus, from \eqref{dx<du}, \eqref{d2x<du} and Remark \ref{bounded.deriv.g},
\begin{equation}\label{D2g.cont}
|D^2 g_t(u,\omega)[\delta u_1,\delta u_2]|\leq L_g(\omega)(C_1^2+C_2)\|\delta u_1\|_{L^2}\|\delta u_2\|_{L^2}
\end{equation}
for all $\omega \in \Omega,$ and, since $L_g\in L^2(\mu,\Omega)$, (iii) follows.

Finally, for all $u,v\in \mathcal O,$
\begin{align*}
D^2 g_t(u,\omega)[\delta u_1,\delta u_2]
&-
D^2 g_t(v,\omega)[\delta u_1,\delta u_2]\\
&=
\left\langle
\nabla_x g(x^u(t),\omega)
-
\nabla_x g(x^v(t),\omega),
\delta^2x_{12}^u(t)
\right\rangle
\\
&\quad+
\left\langle
\nabla_x g(x^v(t),\omega),
\delta^2x_{12}^u(t)
-
\delta^2x_{12}^v(t)
\right\rangle
\\
&\quad+
\left\langle
\Big(
\nabla_{xx}g(x^u(t),\omega)
-
\nabla_{xx}g(x^v(t),\omega)
\Big)
\delta x_1^u(t),
\delta x_2^u(t)
\right\rangle
\\
&\quad+
\left\langle
\nabla_{xx}g(x^v(t),\omega)
\big(
\delta x_1^u(t)
-
\delta x_1^v(t)
\big),
\delta x_2^u(t)
\right\rangle
\\
&\quad+
\left\langle
\nabla_{xx}g(x^v(t),\omega)
\delta x_1^v(t),
\delta x_2^u(t)
-
\delta x_2^v(t)
\right\rangle .
\end{align*}
By Assumption \ref{A2} and Remark \ref{bounded.deriv.g}, $\nabla_xg(\cdot,\omega)$ and $\nabla_{xx}g(\cdot,\omega)$ are Lipschitz continuous and bounded by $L_g(\omega)$, while, by Propositions \ref{traj.lipschitz} and \ref{Diff.x} and Theorem \ref{Diff.x2}, $x^u(t)$, $\delta x_{i}^u$ and $\delta^2 x_{12}^u$ are bounded and Lipschitz continuous with respect to $u$, with deterministic constants. We deduce that there exists $C_{\mathcal H}=c\,L_g$, with $c>0$ depending only on $T,L,m,K_0$ and $r_\rho$, such that
$$
|D^2 g_t(u,\omega)[\delta u_1,\delta u_2]
-
D^2 g_t(v,\omega)[\delta u_1,\delta u_2]|\leq C_{\mathcal H}(\omega)\|\delta u_1\|_{L^2}\|\delta u_2\|_{L^2}\|u-v\|_{L^2},
$$
which completes the proof.
\end{proof}

\subsection{Proof of Lemma \ref{2orderw} }\label{Proof.lemma.2orderw}

\begin{proof} (i) Let $u\in\mathcal O$ and $\delta u\in\Vspace$ with $u+\delta u\in\mathcal O$ and $\|\delta u\|_{L^2}\le1$. By \eqref{gt.taylor},
$$
g_T(u+\delta u,\omega)
=
g_T(u,\omega)
+
Dg_T(u,\omega)[\delta u]
+
\frac12
D^2g_T(u,\omega)[\delta u,\delta u]
+
R_g(T,\omega).
$$
Integrating with respect to $\mu$ and using Proposition \ref{1order},
$$
J(u+\delta u)-J(u)
= DJ(u)[\delta u] +\frac{1}{2} \int_{\Omega} D^2 g_T(u,\omega)[\delta u,\delta u]\, d\mu(\omega) + \int_{\Omega}R_g(T,\omega)\,d\mu(\omega).
$$
For every
$\delta u,\delta v\in\Vspace$ we define
$
D^2J(u)[\delta u,\delta v]
:=
\int_\Omega
D^2g_T(u,\omega)[\delta u,\delta v]
\,d\mu(\omega),
$
which is a continuous bilinear form on $\Vspace$ by \eqref{D2g.cont}. Finally, by \eqref{Rg.bound},
$$
\bigg|\int_{\Omega}R_g(T,\omega)\,d\mu(\omega)\bigg|\le C_R\|L_g\|_{L^1(\mu,\Omega)}\|\delta u\|_{L^2}^3=o(\|\delta u\|_{L^2}^2),
$$
which proves that $J$ is twice Fr\'echet differentiable. The continuity of $D^2J$ follows from (ii).

\noindent(ii) The Lipschitz continuity of $D^2 J$ is derived from that of $D^2 g_T(\cdot,\omega)$ by integrating \eqref{D2gLIPS} with respect to $\mu$, which is legitimate since $C_{\mathcal H}\in L^2(\mu,\Omega)\subset L^1(\mu,\Omega)$.

\noindent(iii) is immediate since the map
$u\mapsto\tfrac\alpha2\|u\|^2_{L^2}$ is a quadratic form with constant
second derivative $\alpha\langle\cdot,\cdot\rangle_{L^2}$.
\end{proof}

\subsection{Proof of Proposition \ref{hessian.HS}}\label{proof.hessian.HS}

\begin{proof}
Fix $u\in\mathcal O$ and $\omega\in\Omega$. By Proposition \ref{2order}, the
second Fr\'echet derivative of $g_T(\cdot,\omega)$ at $u$ is given by
\eqref{F2order}, where, by the variation of constants formula,
\begin{equation}\label{eq:S.representation}
\mathcal{S}(\delta u)(s)
=
\int_{0}^{T}
G^{s}(\tau)\,
\delta u(\tau)\;d\tau,
\qquad
G^{s}(\tau)
:=
\mathbf{1}_{[0,s]}(\tau)\,\Phi^{u}(s,\tau,\omega)\,f_{u}(x^u(\tau,\omega),\omega).
\end{equation}
The $n\times m$ matrix-valued function $G^s(\tau)$ is, by
\eqref{f.EDO.bound} and Remark \ref{bounded.deriv.f}, uniformly bounded in
norm by $c_G:=L_{f}\,e^{\Lambda_{r_\rho}}.$
Moreover, by Lemma~\ref{p.lc.w.t.u}(i) and \eqref{growth.f}, for a.e. $s\in[0,T]$,
\begin{equation}\label{Hij.bounds}
\begin{aligned}
\|H_{xx}(s)\|&\le\|p(s,\omega)\|\,\|f_{xx}(s)\|\le e^{\Lambda_{r_\rho}}L_g(\omega)\,L\big(1+\sqrt m\|u(s)\|\big),
\\
\|H_{xu}(s)\|&=\|H_{ux}(s)\|\le e^{\Lambda_{r_\rho}}L_g(\omega)\,L_f .
\end{aligned}
\end{equation}
We stress that $H_{xx}(s)$ is in general only integrable, and not bounded, in $s$, since $u\in\mathcal O$ is merely square integrable.

Substituting \eqref{eq:S.representation} into~\eqref{F2order} yields an
explicit expression for
$D^{2}g_{T}(u,\omega)[\delta u,\delta v]$, in which the bilinear form
decomposes into two groups of terms.

\smallskip
\noindent\textit{Terminal term:}
$$
\begin{aligned}
&\Bigl\langle
\nabla_{xx}g(x(T),\omega)\,\mathcal{S}(\delta u)(T),
\;\mathcal{S}(\delta v)(T)
\Bigr\rangle_{\mathbb{R}^{n}}\\
&\hspace{4cm}=
\int_{0}^{T}  \int_{0}^{T}
\delta u(\tau)^{\top}
\underbrace{\Bigl[G^{T}(\tau)^{\top}\,
\nabla_{xx}g(x(T),\omega)\,
G^{T}(\sigma)\Bigr]}_{\displaystyle
=:\;\mathbf{K}_T(\tau,\sigma)}
\delta v(\sigma)\;d\sigma\,d\tau .
\end{aligned}
$$
By Remark \ref{bounded.deriv.g}, $\|\mathbf{K}_T(\tau,\sigma)\|\le c_G^{2}\,L_{g}(\omega)$.

\smallskip
\noindent\textit{Hamiltonian integral term.} Recalling that $H_{uu}\equiv0$ by
\eqref{f.affine} and using~\eqref{eq:S.representation}, this term expands into the three sub-terms
$$
\begin{aligned}
&\int_{0}^{T}
\bigl\langle
H_{xx}(s)\,\mathcal{S}(\delta u)(s),
\mathcal{S}(\delta v)(s)
\bigr\rangle\,ds\\
&\hspace{4cm}=
\int_{0}^{T}  \int_{0}^{T}
\delta u(\tau)^{\top}
\underbrace{\Bigl[\int_{\max(\tau,\sigma)}^{T}
G^{s}(\tau)^{\top}H_{xx}(s)G^{s}(\sigma)\,ds\Bigr]}_{
\displaystyle =:\;\mathbf{K}_{xx}(\tau,\sigma)}
\delta v(\sigma)\;d\sigma\,d\tau ,
\end{aligned}
$$
$$
\int_{0}^{T}
\bigl\langle
H_{xu}(s)\,\delta u(s),
\mathcal{S}(\delta v)(s)
\bigr\rangle\,ds
=
\int_{0}^{T}  \int_{0}^{T}
\delta u(\tau)^{\top}
\underbrace{\Bigl[
G^{\tau}(\sigma)^{\top}H_{xu}(\tau)\Bigr]^{\top}}_{
\displaystyle =:\;\mathbf{K}_{xu}(\tau,\sigma)}
\delta v(\sigma)\;d\sigma\,d\tau ,
$$
$$
\int_{0}^{T}
\bigl\langle
H_{ux}(s)\,\mathcal{S}(\delta u)(s),
\delta v(s)
\bigr\rangle\,ds
=
\int_{0}^{T}  \int_{0}^{T}
\delta u(\tau)^{\top}
\underbrace{\Bigl[
H_{ux}(\sigma)\,G^{\sigma}(\tau)\Bigr]^{\top}}_{
\displaystyle =:\;\mathbf{K}_{ux}(\tau,\sigma)}
\delta v(\sigma)\;d\sigma\,d\tau .
$$
By 
\eqref{Hij.bounds}, for all $\tau,\sigma\in[0,T]$,
$$
\|\mathbf K_{xx}(\tau,\sigma)\|
\le c_G^2e^{\Lambda_{r_\rho}}L_g(\omega)\int_0^TL\big(1+\sqrt m\|u(s)\|\big)\,ds
\le c_G^2\Lambda_{r_\rho}e^{\Lambda_{r_\rho}}L_g(\omega),
$$
where the time integration is essential, and
$\|\mathbf K_{xu}(\tau,\sigma)\|,\|\mathbf K_{ux}(\tau,\sigma)\|\le c_GL_fe^{\Lambda_{r_\rho}}L_g(\omega)$.
Consequently, all the above kernels are bounded pointwise by
$C\,L_g(\omega)$, with $C$ depending only on $T,L_f$ and $\Lambda_{r_\rho}$.

Defining the total kernel
\begin{equation*}
\mathbf{K}^{u,\omega}(\tau,\sigma)
:=
\mathbf{K}_T(\tau,\sigma)
+\mathbf{K}_{xx}(\tau,\sigma)
+\mathbf{K}_{xu}(\tau,\sigma)
+\mathbf{K}_{ux}(\tau,\sigma),
\end{equation*}
and the operator $\mathcal H_g(u,\omega)$ as in \eqref{eq:H.definition}, we
obtain, by construction, that for all $\delta u,\delta v\in\Vspace$,
$$
\langle\mathcal{H}_g(u,\omega)\delta u,\delta v\rangle_{L^2}
=
\int_{0}^{T}  \int_{0}^{T}
\delta u(\tau)^{\top}\mathbf{K}^{u,\omega}(\tau,\sigma)\delta v(\sigma)\,
d\sigma\,d\tau
=
D^{2}g_{T}(u,\omega)[\delta u,\delta v] ,
$$
and the uniqueness of $\mathcal{H}_g(u,\omega)$ follows from the Riesz
representation theorem applied to the continuous bilinear form
$D^{2}g_{T}(u,\omega)$, see \eqref{D2g.cont}. Self-adjointness follows from the
symmetry of the second Fr\'echet differential $D^2g_T(u,\omega)$.

Since all four kernels are bounded by $C L_g(\omega)$ pointwise, we get
$$
\|\mathcal{H}_g(u,\omega)\|_{\HS}^{2}
=
\int_{0}^{T}  \int_{0}^{T}
\|\mathbf{K}^{u,\omega}(\tau,\sigma)\|_{\mathrm{F}}^{2}\;d\tau\,d\sigma
\le
c_{\HS}^2\,L_g(\omega)^2
<
\infty,
$$
for a constant $c_{\HS}>0$ depending only on $T,m,n,L_{f}$ and
$\Lambda_{r_\rho}$, which is \eqref{HS.bound}. Since $L_g\in L^2(\mu,\Omega)$,
we deduce that $\mathcal H_g(u,\cdot)\in L^2(\mu,\Omega;\HS(\Vspace))$ and, in
particular, that it is Bochner integrable in $\HS(\Vspace)$, $\mu$ being a
probability measure. Its integral
$\mathcal H_J(u)$ is then a self-adjoint Hilbert--Schmidt operator and, since
Bochner integrals commute with bounded linear maps (here, $\mathcal T\mapsto\langle\mathcal T\delta u,\delta v\rangle_{L^2}$; see also \cite[Example 4.4.6]{Hsing2013}),
for all $\delta u,\delta v\in\Vspace$,
$$
\langle\mathcal H_J(u)\delta u,\delta v\rangle_{L^2}
=
\int_\Omega\langle \mathcal H_g(u,\omega)\delta u,\delta v\rangle_{L^2}\,d\mu(\omega)
=
\int_\Omega D^2g_T(u,\omega)[\delta u,\delta v]\,d\mu(\omega)
=
D^2J(u)[\delta u,\delta v],
$$
by Lemma \ref{2orderw}(i). Finally, \eqref{H.def} follows from
 Lemma \ref{2orderw}(iii), since the Riesz representative of
$(\delta u,\delta v)\mapsto\alpha\langle\delta u,\delta v\rangle_{L^2}$ is
$\alpha\,\mathrm{Id}$.
\end{proof}

\subsection{Proof of Lemma \ref{weak.continuity}}\label{App.weak.continuity}
\begin{proof}
That $u\in\Uad$ follows since $\Uad$ is convex and closed, hence weakly closed.
By Proposition \ref{wellposedness}, $\|x^{u_k}(t,\omega)\|\le R$ for every $t$
and $k,$ and $\|u_k(t)\|\le\rho$ for a.e. $t\in [0,T],$ so that \eqref{growth.f} gives
$
\|\dot x^{u_k}(t,\omega)\|\le L(1+R)\big(1+\sqrt m\,\rho\big)
$  for a.e. $ t\in[0,T].
$
The family $\{x^{u_k}(\cdot,\omega)\}$ is therefore uniformly bounded and
equi-continuous, and by the Arzel\`a--Ascoli theorem, a subsequence (not
relabelled) converges uniformly on $[0,T]$ to some $y\in C([0,T];\mathbb R^n)$.
Writing the state equation in integral form and using \eqref{f.affine},
$$
x^{u_k}(t,\omega)=\varphi(\omega)
+\int_0^tf_0\big(x^{u_k}(s,\omega),\omega\big)ds
+\int_0^tF\big(x^{u_k}(s,\omega),\omega\big)u_k(s)\,ds .
$$
The first integral converges to $\int_0^tf_0(y(s),\omega)\,ds$, since
$f_0(\cdot,\omega)$ is Lipschitz and the convergence is uniform. For the second
one we split
{\small $$
\int_0^tF\big(x^{u_k}(s,\omega),\omega\big)u_k(s)\,ds
=\int_0^t\Big(F\big(x^{u_k}(s,\omega),\omega\big)-F\big(y(s),\omega\big)\Big)u_k(s)\,ds
+\int_0^tF\big(y(s),\omega\big)u_k(s)\,ds .
$$}
The first term is bounded in norm by
$T\rho\,\|F(x^{u_k}(\cdot,\omega),\omega)-F(y(\cdot),\omega)\|_{L^\infty}$,
which tends to $0$ by Remark \ref{bounded.deriv.f} and uniform convergence. In
the second one, each component is of the form
$\langle\mathbf 1_{[0,t]}F(y(\cdot),\omega)^{\!\top}e_j,\,u_k\rangle_{L^2}$ with
$\mathbf 1_{[0,t]}F(y(\cdot),\omega)^{\!\top}e_j\in\Vspace$, and hence converges
to the corresponding expression with $u$ in place of $u_k$. Therefore $y$
solves \eqref{dynamics} with control $u$, so $y=x^u(\cdot,\omega)$ by
uniqueness. Since the limit is independent of the subsequence, the whole
sequence converges.
\end{proof}

\subsection{Proof of Lemma \ref{key.lemma}}\label{App.key.lemma}
\begin{proof}
Suppose, by contradiction, that the conclusion fails. Then there exist
$t_k\downarrow0$, $G_k\in\mathbb K$, $A_k\to\mathcal H$ in
$\mathcal L(L^2)$, and $h_k\in\Vspace$ with $\|h_k\|_{L^2}\le\bar R$ such that
$
\hat u+t_kh_k\in\Uad
$
and \eqref{VI.t} holds with $(t,G,A,h)=(t_k,G_k,A_k,h_k)$, while
$$
\|h_k-\Theta(G_k)\|_{L^2}>\varepsilon
\qquad\text{for all }k.
$$
Since $\mathbb K$ is compact and $\{h_k\}$ is bounded in $\Vspace$, after passing
to a subsequence we may assume
$$
G_k\to G \quad\text{strongly in }L^2,
\qquad
h_k\rightharpoonup\bar h \quad\text{weakly in }\Vspace.
$$
We shall prove that $h_k\to\Theta(G)$ strongly in $\Vspace$, which, together
with the continuity of $\Theta$ from Lemma~\ref{stampacchia}, yields
$\Theta(G_k)\to\Theta(G)$ and contradicts the preceding inequality. Since $\hat u+t_kh_k\in\Uad$, we have $h_k\in\mathcal R(\hat u).$ Hence
$
\bar h\in\overline{\mathcal R(\hat u)}^{\,L^2}.
$
The latter set 
is closed and convex, and
therefore weakly closed, so $\bar h$ belongs to the tangent cone
$\mathfrak T(\hat u)$.

It remains to verify the criticality condition. Testing \eqref{VI} at $\hat u$
with $u=\hat u+t_kh_k$ gives
$
\langle d,h_k\rangle_{L^2}\ge0.
$
On the other hand, testing \eqref{VI.t} with $v=\hat u$ yields
$$
\langle d,h_k\rangle_{L^2}
\le
-t_k\langle G_k+A_kh_k,h_k\rangle_{L^2}.
$$
Since $\mathbb K$ is compact, $\{A_k\}$ is bounded in
$\mathcal L(L^2)$, and $\|h_k\|_{L^2}\le\bar R$, we obtain
$$
0\le\langle d,h_k\rangle_{L^2}
\le
t_k
\left(
\sup_{G'\in\mathbb K}\|G'\|_{L^2}
+
\sup_j\|A_j\|_{\mathcal L(L^2)}\bar R
\right)\bar R
\longrightarrow0.
$$
Passing to the limit gives
$
\langle d,\bar h\rangle_{L^2}=0.
$
Consequently,
$
\bar h\in\mathcal C(\hat u).
$ We next identify $\bar h$. Let $w\in\mathcal D(\hat u)$. By definition,
there exist $\lambda>0$ and $u\in\Uad$ such that
$
w=\lambda(u-\hat u)
$ and $
\langle d,w\rangle_{L^2}=0.
$
For $t_k\le1/\lambda$, convexity of $\Uad$ implies
$\hat u+t_kw\in\Uad$. We may therefore use
$v=\hat u+t_kw$ as a test function in \eqref{VI.t}. After dividing by $t_k$,
we obtain
$$
\langle d,w-h_k\rangle_{L^2}
+t_k\langle G_k+A_kh_k,w-h_k\rangle_{L^2}\ge0.
$$
Since $\langle d,w\rangle_{L^2}=0$ and
$\langle d,h_k\rangle_{L^2}\ge0$, it follows that
$$
\langle G_k+A_kh_k,w-h_k\rangle_{L^2}\ge0,
$$
or equivalently,
\begin{equation}\label{key.ineq}
\langle A_kh_k,h_k\rangle_{L^2}
\le
\langle G_k,w-h_k\rangle_{L^2}
+\langle A_kh_k,w\rangle_{L^2}.
\end{equation}

We now pass to the limit. Since $G_k\to G$ strongly,
$h_k\rightharpoonup\bar h$ weakly, and
$A_k^*w\to\mathcal H w$ strongly, while $\mathcal H$ is self-adjoint,
$$
\langle G_k,w-h_k\rangle_{L^2}
+\langle A_kh_k,w\rangle_{L^2}
\longrightarrow
\langle G,w-\bar h\rangle_{L^2}
+\langle\mathcal H\bar h,w\rangle_{L^2}.
$$
Moreover, since $A_k\to\mathcal H$ in operator norm and $\{h_k\}$ is bounded,
$$
\langle A_kh_k,h_k\rangle_{L^2}
=
Q(h_k)+o(1).
$$
Thus \eqref{key.ineq} gives
\begin{equation}\label{key.limsup}
\limsup_{k\to\infty}Q(h_k)
\le
\langle G,w-\bar h\rangle_{L^2}
+\langle\mathcal H\bar h,w\rangle_{L^2}
\qquad
\forall\,w\in\mathcal D(\hat u).
\end{equation}

By Lemma~\ref{weak.cont.lemma}, the quadratic form $Q$ is weakly lower
semicontinuous, and hence
$$
Q(\bar h)\le\liminf_{k\to\infty}Q(h_k).
$$
Combining this with \eqref{key.limsup}, we obtain
$$
\langle\mathcal H\bar h+G,w-\bar h\rangle_{L^2}\ge0
\qquad
\forall\,w\in\mathcal D(\hat u).
$$
The density property \eqref{poly.density} allows us to extend this inequality
to every $w\in\mathcal C(\hat u)$. Since $\bar h\in\mathcal C(\hat u)$, the
characterization of $\Theta(G)$ in Lemma~\ref{stampacchia} implies
$$
\bar h=\Theta(G).
$$

It remains to upgrade weak convergence to strong convergence. By the density
of $\mathcal D(\hat u)$ in $\mathcal C(\hat u)$, choose
$w_j\in\mathcal D(\hat u)$ such that $w_j\to\bar h$ strongly in $L^2$. For
each fixed $j$, \eqref{key.limsup} gives
$$
\limsup_{k\to\infty}Q(h_k)
\le
\langle G,w_j-\bar h\rangle_{L^2}
+\langle\mathcal H\bar h,w_j\rangle_{L^2}.
$$
Letting $j\to\infty$ yields
$$
\limsup_{k\to\infty}Q(h_k)
\le
\langle\mathcal H\bar h,\bar h\rangle_{L^2}
=Q(\bar h).
$$
Together with weak lower semicontinuity,
$
Q(h_k)\longrightarrow Q(\bar h).
$ The Legendre property of $Q$, established in Lemma~\ref{weak.cont.lemma},
then implies
$$
h_k\longrightarrow\bar h=\Theta(G)
\ \text{ strongly in }L^2.
$$
Finally, Lemma~\ref{stampacchia} gives
$\Theta(G_k)\to\Theta(G)$, and therefore
$
\|h_k-\Theta(G_k)\|_{L^2}\rightarrow0,
$
contradicting $\|h_k-\Theta(G_k)\|_{L^2}>\varepsilon$ for every $k$.
\end{proof}
\section{Second-order optimality conditions in Banach spaces}\label{appendixB}

We collect here the abstract results of
\cite{casas2012,CasTro2012} that are used in Section  \ref{sec:sosc}, stated with the
notation of that section. Let $U_\infty$ be a Banach space and $U_2$ a Hilbert space, with norms $\|\cdot\|_{U_\infty}$ and $\|\cdot\|_{U_2}$, and assume that
$U_\infty\subset U_2$ with continuous embedding; the case $U_\infty=U_2$ is
allowed. Let $\mathcal U\subset U_\infty$ be a nonempty convex set, let
$A\subset U_\infty$ be an open set with $\mathcal U\subset A$, and let
$\mathfrak J:A\to\mathbb R$. Consider the abstract optimization problem
\begin{equation}\label{abstract.P}\tag{$\mathfrak P$}
\min_{u\in \mathcal U}\mathfrak J (u).
\end{equation}
All notions of differentiability of $\mathfrak J$ are understood in the sense of
$U_\infty$, and for $\hat u\in \mathcal U$ and $\varepsilon>0$ we write
$B_{U_2}(\hat u,\varepsilon):=\{u\in U_2:\|u-\hat u\|_{U_2}<\varepsilon\}$.
We say that $\tilde u\in \mathcal U$ is a {\it stationary point} of \eqref{abstract.P}
if $D\mathfrak J(\tilde u)[u-\tilde u]\ge0$ for all $u\in \mathcal U$. We fix $\hat u\in \mathcal U$ and
consider the following assumptions.

\vspace{0.2cm}
\begin{itemize}
\item[{\bf (A1)}] The functional $\mathfrak J:A\to\mathbb R$ is of class $C^2$ and, for
every $u\in \mathcal U$, the derivatives admit continuous extensions
$
D\mathfrak J(u)\in\mathcal L(U_2,\mathbb R),
$ and $
D^2\mathfrak J(u)\in\mathcal B(U_2\times U_2;\mathbb R).
$

\item[{\bf (A2)}] For every sequence $\{u_k\}_{k\in\mathbb N}\subset \mathcal U$ with
$\|u_k-\hat u\|_{U_2}\to0$ and every sequence
$\{\delta u_k\}_{k\in\mathbb N}\subset U_2$ with
$\delta u_k\rightharpoonup\delta u$ weakly in $U_2$, one has
\begin{align}
&D\mathfrak J(\hat u)[\delta u]=\lim_{k\to\infty}D\mathfrak J(u_k)[\delta u_k],
\label{CT.2.3}
\\
&D^2\mathfrak J(\hat u)[\delta u,\delta u]\le\liminf_{k\to\infty}D^2\mathfrak J(u_k)[\delta u_k,\delta u_k],
\label{CT.2.4}
\\
&\text{if } \delta u=0, \text{ then }
\vartheta\liminf_{k\to\infty}\|\delta u_k\|_{U_2}^2
\le
\liminf_{k\to\infty}D^2\mathfrak J(u_k)[\delta u_k,\delta u_k],
\label{CT.2.5}
\end{align}
for some constant $\vartheta>0$.
\end{itemize}
\vspace{0.2cm}

Associated with $\hat u$, one defines the {\it radial cone}, the
{\it critical cone} and the set of {\it critical radial directions} by
\begin{equation}\label{CT.cones}
\begin{aligned}
\mathcal R(\hat u)
&:=
\bigl\{\delta u\in U_\infty:\ \delta u=\lambda(u-\hat u)\ \text{for some }\lambda>0,\ u\in \mathcal U\bigr\},
\\  
\mathcal C(\hat u)
&:=
\overline{\mathcal R(\hat u)}^{\,U_2}\cap\bigl\{\delta u\in U_2:\ D\mathfrak J(\hat u)[\delta u]=0\bigr\},
\ \text{ and }
\mathcal D(\hat u)
:=
\bigl\{\delta u\in\mathcal R(\hat u):\ D\mathfrak J(\hat u)[\delta u]=0\bigr\},
\end{aligned}
\end{equation}
respectively, where $\overline{\,\cdot\,}^{\,U_2}$ denotes the closure in $U_2$.

The following result provides second-order \emph{necessary} optimality
conditions.

\begin{theorem}[{\cite[Theorem 2.3]{CasTro2012}}]\label{CT.necessary}
Let $\hat u$ be a local solution of \eqref{abstract.P} in $U_\infty$. Assume that
{\bf (A1)} holds and that
$\mathcal C(\hat u)=\overline{\mathcal D(\hat u)}^{\,U_2}$. Then
$
D^2\mathfrak J(\hat u)[\delta u,\delta u]\ge0$ for all $\delta u\in\mathcal C(\hat u).
$
\end{theorem}

\noindent The condition $\mathcal C(\hat u)=\overline{\mathcal D(\hat u)}^{\,U_2}$
is equivalent to the polyhedricity of $\mathcal U$ at $\hat u$; see \cite[\S2]{CasTro2012} and
\cite[\S3.2]{BonnansShapiro2000}.

The following result provides second-order \emph{sufficient} optimality
conditions, in the form of a quadratic growth property.

\begin{theorem}[{\cite[Theorem 2.3]{casas2012}}]\label{CT.growth}
Suppose that {\bf (A1)} and {\bf (A2)} hold, and let $\hat u\in \mathcal U$ be a
stationary point of \eqref{abstract.P} such that
$ D^2\mathfrak J(\hat u)[\delta u,\delta u]>0
 \text{ for all }\delta u\in\mathcal C(\hat u)\setminus\{0\}.
$
Then there exist $\varepsilon>0$ and $\delta>0$ such that
$$
\mathfrak J(\hat u)+\frac{\delta}{2}\|u-\hat u\|_{U_2}^2\le \mathfrak J(u)
\qquad\text{for all } u\in \mathcal U\cap B_{U_2}(\hat u,\varepsilon).
$$
\end{theorem}
\subsection*{Statement of LLM usage}

During the preparation of this manuscript, the authors used Claude Sonnet 5.5 (Anthropic) to improve the language, refine the structure and presentation of the text, and assist in developing the Python code used in the numerical experiments. All AI-assisted content was carefully reviewed and revised by the authors, who assume full responsibility for the final manuscript.
\bibliographystyle{plain}
\bibliography{generalsoledad}

@article{Vinter2005,
	Author = {Vinter, R.},
	Journal = {SIAM journal on control and optimization},
	Number = {3},
	Pages = {939--968},
	Publisher = {SIAM},
	Title = {Minimax optimal control},
	Volume = {44},
	Year = {2005}}

@article{Ross2015,
	Author = {Ross, I. and Proulx, R.J. and Karpenko, M. and Gong, Q.},
	Journal = {Journal of Guidance, Control, and Dynamics},
	Number = {7},
	Pages = {1251--1263},
	Publisher = {American Institute of Aeronautics and Astronautics},
	Title = {{Riemann--Stieltjes} optimal control problems for uncertain dynamic systems},
	Volume = {38},
	Year = {2015}}

@article{casas2012,
author = {Casas, Eduardo and Tröltzsch, Fredi},
year = {2012},
month = {03},
pages = {261-279},
title = {Second Order Analysis for Optimal Control Problems: Improving Results Expected From Abstract Theory},
volume = {22},
journal = {SIAM Journal on Optimization},
doi = {10.1137/110840406}
}

@book{brezis2010functional,
	Author = {Brezis, H.},
	Publisher = {Springer Science \& Business Media},
	Title = {Functional analysis, {Sobolev} spaces and partial differential equations},
	Year = {2010}}

@inproceedings{I}

@article{palladino2025convergence,
  title={Convergence results for control problems with unknown dynamic and applications to reinforcement learning},
  author={Palladino, M. and Pesare, A. and Scarinci, T.},
  journal={Mathematical Control and Related Fields},
  volume={17},
  pages={99--119},
  year={2026}
}

@book{rudin1973functional,
	Author = {Rudin, W.},
	Publisher = {McGraw-Hill, New York},
	Title = {Functional analysis},
	Year = {1973}}

@book{ BonnansShapiro2000,
	Author = {J.F. Bonnans and A. Shapiro},
	Publisher = {Springer-Verlag, New York},
	Title = {Perturbation analysis of optimization problems},
	Year = {2000}}

@article{BettiolKhalil,
      title={Necessary Optimality Conditions For Average Cost Minimization Problems},
      author={P. Bettiol and N. Khalil},
      journal = {Discrete and Continuous Dynamical Systems - B},
      volume = {24},
      number = {5},
      pages = {2093--2124},
      year={2019}
}

@article{zuazua2014,
title = {Averaged control},
journal = {Automatica},
volume = {50},
number = {12},
pages = {3077-3087},
year = {2014},
issn = {0005-1098},
author = {E. Zuazua}
}

@article{RossKarpenko2015,
author = {Ross, I. and Karpenko, Mark and Proulx, Ronald},
year = {2015},
month = {07},
pages = {5599-5604},
title = {A {Lebesgue-Stieltjes} framework for optimal control and allocation},
volume = {2015},
journal = {Proceedings of the American Control Conference},
doi = {10.1109/ACC.2015.7172216}
}

@article{Phelps2013,
  author={Phelps, C. and Royset, J.O. and Gong, Q.},
    title={Sample average approximations in optimal control of uncertain systems}, 
   volume={},
  number={},
pages={1958-1965},
 year={2013},
  journal = {52nd IEEE Conference on Decision and Control, Firenze, Italy},
  doi={10.1109/CDC.2013.6760168}}

@article{PesarePalladinoFalcone.LQR.2021,
author = {Palladino, M. and Pesare, A. and Falcone, M.},
year = {2021},
month = {09},
pages = {},
title = {Convergence results for an averaged {LQR} problem with applications to reinforcement learning},
volume = {33},
journal = {Mathematics of Control Signals and Systems},
doi = {10.1007/s00498-021-00294-y}
}

@article{Scagliotti2023,
   title={Optimal control of ensembles of dynamical systems},
   volume={29},
   ISSN={1262-3377},
   url={http://dx.doi.org/10.1051/cocv/2023011},
   DOI={10.1051/cocv/2023011},
   journal={ESAIM: Control, Optimisation and Calculus of Variations},
   publisher={EDP Sciences},
   author={Scagliotti, A.},
   year={2023},
   pages={22} }

@article{Zuzua-Loheac,
  TITLE = {{From averaged to simultaneous controllability}},
  AUTHOR = {Loh{\'e}ac, J. and Zuazua, E.},
  URL = {https://hal.science/hal-01164525},
  JOURNAL = {{Annales de la Facult{\'e} des Sciences de Toulouse. Math{\'e}matiques.}},
  PUBLISHER = {{Universit{\'e} Paul Sabatier \_ Cellule Mathdoc }},
  VOLUME = {25},
  NUMBER = {4},
  PAGES = {785 - 828},
  YEAR = {2016},
  MONTH = Oct,
  DOI = {10.5802/afst.1511},
  HAL_ID = {hal-01164525},
  HAL_VERSION = {v3},
}

@book{polak2012,
  title={Optimization: Algorithms and Consistent Approximations},
  author={Polak, E.},
  isbn={9781461206637},
  lccn={97002158},
  series={Applied Mathematical Sciences},
  url={https://books.google.com.ec/books?id=uoXuBwAAQBAJ},
  year={2012},
  publisher={Springer New York}
}

@article{AttouchWets1983,
  author    = {Attouch, H{\'e}dy and Wets, Roger J.-B.},
  title     = {A convergence theory for saddle functions},
  journal   = {Transactions of the American Mathematical Society},
  volume    = {280},
  year      = {1983},
  pages     = {1--41}
}

@book{shapiro2021,
  title={Lectures on Stochastic Programming: Modeling and Theory},
  author={Shapiro, A. and Dentcheva, D. and Ruszczy{\'n}ski, A.P.},
  isbn={9781611976588},
  lccn={2014005486},
  series={MOS-SIAM series on optimization},
  url={https://books.google.com.br/books?id=vaY8zgEACAAJ},
  year={2021},
  publisher={Society for Industrial and Applied Mathematics}
}

@article{ArtsteinWets1995,
  author    = {Zvi Artstein and Roger J.-B. Wets},
  title     = {Consistency of minimizers and the SLLN for stochastic programs},
  journal   = {Journal of Convex Analysis},
  volume    = {2},
  number    = {1-2},
  year      = {1995},
  pages     = {1--17}
}

@book{AraujoGine1980,
  author    = {Aloisio Araujo and Evarist Gin{\'e}},
  title     = {The Central Limit Theorem for Real and Banach Valued Random Variables},
  publisher = {John Wiley \& Sons},
  address   = {New York},
  year      = {1980},
  series    = {Wiley Series in Probability and Mathematical Statistics},
  isbn      = {978-0471053040}
}

@phdthesis{Phelps2015,
  author = {C. Phelps},
  title = {Computational Optimal Control of Nonlinear Systems with Parameter Uncertainty},
  school = {University of California, Santa Cruz},
  address = {Santa Cruz, CA},
  year = {2015},
  type = {Ph.D. thesis}
}

@article{WEINAN2018,
  author = {E, Weinan and Han, Jiequn and Li, Qianxiao},
  year = {2018},
  date = {2018-12-13},
  title = {A mean-field optimal control formulation of deep learning},
  journal = {Research in the Mathematical Sciences},
  pages = {10},
  volume = {6},
  number = {1},
  doi = {10.1007/s40687-018-0172-y},
  url = {https://doi.org/10.1007/s40687-018-0172-y},
  issn = {2197-9847}
}

@INPROCEEDINGS{PesarePalladinoFalcone2021b,
  author={Palladino, M. and Pesare, A. and Falcone, M.},
  booktitle={2021 European Control Conference (ECC)}, 
  title={Convergence of the Value Function in Optimal Control Problems with Unknown Dynamics}, 
  year={2021},
  volume={},
  number={},
  pages={2426-2431},
  doi={10.23919/ECC54610.2021.9655079}}

@book{Hsing2013,
  author    = {Tailen Hsing and Randall Eubank},
  title     = {Theoretical Foundations of Functional Data Analysis, with an Introduction to Linear Operators},
  publisher = {John Wiley \& Sons},
  year      = {2013},
  doi       = {10.1002/9781118762547}
}

@book{bill99,
  author    = {P. Billingsley},
  title     = {Convergence of Probability Measures},
  edition   = {2},
  publisher = {John Wiley \& Sons},
  address   = {New York},
  year      = {1999},
  series    = {Wiley Series in Probability and Statistics}
}

@article{melnikov2026convergence,
  title={Convergence rates for ensemble-based solutions to optimal control of uncertain dynamical systems},
  author={Melnikov, Olena and Milz, Johannes},
  journal={SIAM Journal on Control and Optimization},
  volume={64},
  number={3},
  pages={1631--1650},
  year={2026},
  publisher={SIAM}
}

@article{romisch2024asymptotic,
  title={Asymptotic properties of {M}onte {C}arlo methods in elliptic {PDE}-constrained optimization under uncertainty},
  author={R{\"o}misch, Werner and Surowiec, Thomas M},
  journal={Numerische Mathematik},
  volume={156},
  number={5},
  pages={1887--1914},
  year={2024},
  publisher={Springer}
}

@article{milz2024sample,
  title={Sample size estimates for risk-neutral semilinear PDE-constrained optimization},
  author={Milz, Johannes and Ulbrich, Michael},
  journal={SIAM Journal on Optimization},
  volume={34},
  number={1},
  pages={844--869},
  year={2024},
  publisher={SIAM}
}

@book{Kurzweil1986,
  author    = {Kurzweil, Jaroslav},
  title     = {Ordinary Differential Equations: Introduction to the Theory of Ordinary Differential Equations in the Real Domain},
  series    = {Studies in Applied Mechanics},
  volume    = {13},
  publisher = {Elsevier},
  address   = {Amsterdam},
  year      = {1986},
  isbn      = {0-444-99509-9}
}

@article{Ioffe1979III,
    author  = {A. D. Ioffe},
    title   = {Necessary and Sufficient Conditions for a Local Minimum. III: Second Order Conditions and Augmented Duality},
    journal = {SIAM Journal on Control and Optimization},
    volume  = {17},
    number  = {2},
    pages   = {266--288},
    year    = {1979},
    month   = {mar},
    doi     = {10.1137/0317021}
}

@book{florescu2021lebesgue,
  author    = {Florescu, Liviu C.},
  title     = {Lebesgue Integral},
  series    = {Compact Textbooks in Mathematics},
  publisher = {Birkh{\"a}user},
  address   = {Cham},
  edition   = {1},
  year      = {2021},
  isbn      = {978-3-030-60162-1},
  doi       = {10.1007/978-3-030-60163-8}
}

@article{AronnaMonteiroSierra2025,
  author  = {Aronna, M. S. and de Lima Monteiro, Gabriel and
             Sierra Fonseca, O.},
  title   = {Average Optimal Control of Uncertain Control-Affine Systems},
  journal = {Set-Valued and Variational Analysis},
  volume  = {33},
  pages   = {45},
  year    = {2025},
  doi     = {10.1007/s11228-025-00784-0}
}

@article{CasTro2012,
  author  = {Casas, Eduardo and Tr{\"o}ltzsch, Fredi},
  title   = {A general theorem on error estimates with application to a quasilinear elliptic optimal control problem},
  journal = {Computational Optimization and Applications},
  volume  = {53},
  number  = {1},
  pages   = {173--206},
  year    = {2012},
  month   = sep,
  publisher = {Springer},
  doi     = {10.1007/s10589-011-9453-8}
}

@article{KingRockafellar1993,
 ISSN = {0364765X, 15265471},
 URL = {http://www.jstor.org/stable/3690158},
 author = {Alan J. King and R. Tyrrell Rockafellar},
 journal = {Mathematics of Operations Research},
 number = {1},
 pages = {148--162},
 publisher = {INFORMS},
 title = {Asymptotic Theory for Solutions in Statistical Estimation and Stochastic Programming},
 urldate = {2026-09-29},
 volume = {18},
 year = {1993}
}

@article{Shapiro1993,
 ISSN = {0364765X, 15265471},
 URL = {http://www.jstor.org/stable/3690125},
 author = {Alexander Shapiro},
 journal = {Mathematics of Operations Research},
 number = {4},
 pages = {829--845},
 publisher = {INFORMS},
 title = {Asymptotic Behavior of Optimal Solutions in Stochastic Programming},
 urldate = {2026-09-29},
 volume = {18},
 year = {1993}
}

@article{PhelpsRoysetGong2016,
  author  = {Phelps, C. and Royset, J. O. and Gong, Q.},
  title   = {Optimal control of uncertain systems using sample average approximations},
  journal = {SIAM J. Control Optim.},
  volume  = {54},
  number  = {1},
  pages   = {1--29},
  year    = {2016}
}

@book{aliprantis2006infinite,
  author    = {Aliprantis, C. D. and Border, K. C.},
  title     = {Infinite Dimensional Analysis. A Hitchhiker's Guide},
  edition   = {Third},
  publisher = {Springer, Berlin},
  year      = {2006}
}

@article{fournier2015rate,
  author  = {Fournier, N. and Guillin, A.},
  title   = {On the rate of convergence in {W}asserstein distance of the empirical measure},
  journal = {Probab. Theory Related Fields},
  volume  = {162},
  number  = {3-4},
  pages   = {707--738},
  year    = {2015}
}

@article{LiKhaneja2006,
  author  = {Li, J.-S. and Khaneja, N.},
  title   = {Control of inhomogeneous quantum ensembles},
  journal = {Phys. Rev. A},
  volume  = {73},
  number  = {3},
  year    = {2006}
}

@article{RuthsLi2012,
  author  = {Ruths, J. and Li, J.-S.},
  title   = {Optimal control of inhomogeneous ensembles},
  journal = {IEEE Trans. Automat. Control},
  volume  = {57},
  number  = {8},
  pages   = {2021--2032},
  year    = {2012}
}

@article{MarigondaQuincampoix2018,
  author  = {Marigonda, A. and Quincampoix, M.},
  title   = {Mayer control problem with probabilistic uncertainty on initial positions},
  journal = {J. Differential Equations},
  volume  = {264},
  number  = {5},
  pages   = {3212--3252},
  year    = {2018}
}

@article{BettiolKhalil2024,
  author  = {Bettiol, P. and Khalil, N.},
  title   = {Average cost minimization problems subject to state constraints},
  journal = {SIAM J. Control Optim.},
  volume  = {62},
  number  = {3},
  pages   = {1884--1907},
  year    = {2024}
}

@article{AronnaPalladinoSierra2026,
  author  = {Aronna, M. S. and Palladino, M. and Sierra Fonseca, O.},
  title   = {Dynamic programming principle and {H}amilton--{J}acobi--{B}ellman equation for optimal control problems with uncertainty},
  journal = {ESAIM Control Optim. Calc. Var.},
  volume  = {32},
  pages   = {5},
  year    = {2026}
}

@article{Scagliotti2025minimax,
  author  = {Scagliotti, A.},
  title   = {Minimax problems for ensembles of control-affine systems},
  journal = {SIAM J. Control Optim.},
  volume  = {63},
  number  = {1},
  pages   = {502--523},
  year    = {2025}
}

@misc{ScagliottiSurowiec2026risk,
  author       = {Scagliotti, A. and Surowiec, T. M.},
  title        = {Risk-averse ensemble control for control-affine systems},
  howpublished = {arXiv preprint arXiv:2605.02791},
  year         = {2026}
}

@article{BonalliBonnet2023,
  author  = {Bonalli, R. and Bonnet, B.},
  title   = {First-order {P}ontryagin maximum principle for risk-averse stochastic optimal control problems},
  journal = {SIAM J. Control Optim.},
  volume  = {61},
  number  = {3},
  pages   = {1881--1909},
  year    = {2023}
}

@article{KouriSurowiec2016,
  author  = {Kouri, D. P. and Surowiec, T. M.},
  title   = {Risk-averse {PDE}-constrained optimization using the conditional value-at-risk},
  journal = {SIAM J. Optim.},
  volume  = {26},
  number  = {1},
  pages   = {365--396},
  year    = {2016}
}

@article{HoffhuesRomischSurowiec2021,
  author  = {Hoffhues, M. and R\"omisch, W. and Surowiec, T. M.},
  title   = {On quantitative stability in infinite-dimensional optimization under uncertainty},
  journal = {Optim. Lett.},
  volume  = {15},
  number  = {8},
  pages   = {2733--2756},
  year    = {2021}
}

@article{MilzSurowiec2024,
  author  = {Milz, J. and Surowiec, T. M.},
  title   = {Asymptotic consistency for nonconvex risk-averse stochastic optimization with infinite dimensional decision spaces},
  journal = {Math. Oper. Res.},
  volume  = {49},
  number  = {3},
  pages   = {1403--1418},
  year    = {2024}
}

@article{Milz2025stationary,
  author  = {Milz, J.},
  title   = {Consistency of sample-based stationary points for infinite-dimensional stochastic optimization},
  journal = {SIAM J. Optim.},
  volume  = {35},
  number  = {1},
  pages   = {42--61},
  year    = {2025}
}

@article{Milz2023SAA,
  author  = {Milz, J.},
  title   = {Sample average approximations of strongly convex stochastic programs in {H}ilbert spaces},
  journal = {Optim. Lett.},
  volume  = {17},
  pages   = {471--492},
  year    = {2023}
}

@article{LewBonalliPavone2024,
  author  = {Lew, T. and Bonalli, R. and Pavone, M.},
  title   = {Sample average approximation for stochastic programming with equality constraints},
  journal = {SIAM J. Optim.},
  volume  = {34},
  pages   = {3506--3533},
  year    = {2024}
}

@article{DupacovaWets1988,
  author  = {Dupa\v{c}ov\'a, J. and Wets, R. J.-B.},
  title   = {Asymptotic behavior of statistical estimators and of optimal solutions of stochastic optimization problems},
  journal = {Ann. Statist.},
  volume  = {16},
  pages   = {1517--1549},
  year    = {1988}
}

@misc{MilzShapiro2026,
  author       = {Milz, J. and Shapiro, A.},
  title        = {Central limit theorems for sample average approximations in stochastic optimal control},
  howpublished = {arXiv preprint arXiv:2508.01942},
  year         = {2026}
}

@misc{JaveedMilz2026,
  author       = {Javeed, A. and Milz, J.},
  title        = {Statistical inference for scenario-based dynamic optimization under uncertainty},
  howpublished = {arXiv preprint arXiv:2607.14965},
  year         = {2026}
}

@article{FanLiuBelabbas2019,
  author  = {Fan, Y. and Liu, S. and Belabbas, M.-A.},
  title   = {Mid-air motion planning of floating robot using heat flow method},
  journal = {IFAC-PapersOnLine},
  volume  = {52},
  number  = {22},
  pages   = {19--24},
  year    = {2019},
  doi     = {10.1016/j.ifacol.2019.11.041}
}

@article{ScagliottiEtAl2026cancer,
  author  = {Scagliotti, A. and Scagliotti, F. and Locati, L. D. and Sottotetti, F.},
  title   = {Ensemble optimal control for managing drug resistance in cancer therapies},
  journal = {Bull. Math. Biol.},
  volume  = {88},
  pages   = {92},
  year    = {2026},
  doi     = {10.1007/s11538-026-01661-z}
}

@book{vanderVaart1998,
  author    = {van der Vaart, A. W.},
  title     = {Asymptotic Statistics},
  series    = {Cambridge Series in Statistical and Probabilistic Mathematics},
  volume    = {3},
  publisher = {Cambridge University Press},
  year      = {1998}
}

@book{vanderVaartWellner1996,
  author    = {van der Vaart, A. W. and Wellner, J. A.},
  title     = {Weak Convergence and Empirical Processes},
  series    = {Springer Series in Statistics},
  publisher = {Springer, New York},
  year      = {1996}
}

@article{JainMarcus1975,
  author  = {Jain, N. C. and Marcus, M. B.},
  title   = {Central limit theorems for {$C(S)$}-valued random variables},
  journal = {J. Funct. Anal.},
  volume  = {19},
  pages   = {216--231},
  year    = {1975}
}

@article{AronnaDumont,
author = {Aronna, M. S. and Dumont, Y.},
year = {2020},
month = {08},
pages = {},
title = {On Nonlinear Pest/Vector Control via the Sterile Insect Technique: Impact of Residual Fertility},
volume = {82},
journal = {Bulletin of Mathematical Biology},
doi = {10.1007/s11538-020-00790-3}
}

@article{BidiCoronHayat,
      title={Reinforcement Learning in Control Theory: A New Approach to Mathematical Problem Solving}, 
      author={Bidi, K.  and  Coron, J.M. and Hayat, A. and Lichtlé N.},
      year={2023},
      eprint={2310.13072},
      archivePrefix={arXiv},
      primaryClass={math.OC},
      journal ={3rd Workshop on Mathematical
Reasoning and AI at NeurIPS’23},
url={https://arxiv.org/abs/2310.13072}, 
}

@INPROCEEDINGS{AMS.SIT,
  author={Aronna, M.S. and de Lima Monteiro, G. and Sierra, O.},
  booktitle={2025 IEEE 64th Conference on Decision and Control (CDC)}, 
  title={Singular Arcs in Average Optimal Control-Affine Problems}, 
  year={2025},
  volume={},
  number={},
  pages={3748-3753},
  doi={10.1109/CDC57313.2025.11312784}}

@article{machado2025nash,
  title={From {N}ash to {C}ournot--{N}ash equilibria via {$\Gamma$}-convergence},
  author={Machado, Jo{\~a}o Miguel and Mazanti, Guilherme and Pfeiffer, Laurent},
  journal={arXiv preprint arXiv:2509.02205},
  year={2025}
}

@book{BauschkeCombettes2017,
  author    = {Bauschke, Heinz H. and Combettes, Patrick L.},
  title     = {Convex Analysis and Monotone Operator Theory in Hilbert Spaces},
  edition   = {2},
  series    = {CMS Books in Mathematics},
  publisher = {Springer, Cham},
  year      = {2017},
  doi       = {10.1007/978-3-319-48311-5}
}

@article{KolmogorovTikhomirov1959,
  author  = {Kolmogorov, A. N. and Tikhomirov, V. M.},
  title   = {$\varepsilon$-entropy and $\varepsilon$-capacity of sets in function spaces},
  journal = {Uspekhi Matematicheskikh Nauk},
  volume  = {14},
  number  = {2(86)},
  pages   = {3--86},
  year    = {1959},
  note    = {English transl.: Amer. Math. Soc. Transl. (2) 17 (1961), 277--364}
}

\end{document}